\documentclass{amsart}
\usepackage{amscd,amssymb,amsmath,amsbsy,amsthm}
\usepackage[all]{xy}
\xyoption{arc}
\usepackage[colorlinks,plainpages,urlcolor=blue]{hyperref}

\newtheorem{theorem}{Theorem}[section]
\newtheorem{lemma}[theorem]{Lemma}
\newtheorem{corollary}[theorem]{Corollary}
\newtheorem{prop}[theorem]{Proposition}

\theoremstyle{definition}
\newtheorem{definition}[theorem]{Definition}
\newtheorem{example}[theorem]{Example}
\newtheorem{remark}[theorem]{Remark}

\newcommand{\Z}{\mathbb{Z}}
\newcommand{\Q}{\mathbb{Q}}
\newcommand{\R}{\mathbb{R}}
\newcommand{\C}{\mathbb{C}}

\newcommand{\QP}{\mathbb{QP}}
\newcommand{\CP}{\mathbb{CP}}

\newcommand{\bP}{\mathbb{P}}

\newcommand{\cS}{\mathsf{p}}

\newcommand{\vG}{\Gamma}
\newcommand{\T}{T}

\renewcommand{\k}{\Bbbk}
\newcommand{\RR}{\mathcal{R}}
\newcommand{\VV}{\mathcal{V}}
\newcommand{\A}{{\mathcal{A}}}
\newcommand{\B}{{\mathcal{B}}}
\renewcommand{\L}{{\C}}
\newcommand{\CC}{{\mathcal{C}}}
\newcommand{\WW}{\mathcal{W}}

\DeclareMathOperator{\rank}{rank}

\DeclareMathOperator{\im}{im}

\DeclareMathOperator{\codim}{codim}

\DeclareMathOperator{\id}{id}
\DeclareMathOperator{\ab}{{ab}}

\DeclareMathOperator{\Hom}{{Hom}}
\DeclareMathOperator{\spn}{span}
\DeclareMathOperator{\ann}{{ann}}

\DeclareMathOperator{\FF}{{F}}
\DeclareMathOperator{\FP}{{FP}}

\DeclareMathOperator{\lcm}{{lcm}}

\DeclareMathOperator{\Tors}{Tors}

\DeclareMathOperator{\Grass}{Gr}
\DeclareMathOperator{\TC}{TC}
\DeclareMathOperator{\pr}{pr}

\DeclareMathOperator{\orb}{orb}

\newcommand{\wX}{\widetilde{X}}
\newcommand{\wG}{\widehat{G}}

\newcommand{\same}{\Longleftrightarrow}
\newcommand{\surj}{\twoheadrightarrow}
\newcommand{\inj}{\hookrightarrow}
\newcommand{\isom}{\xrightarrow{\,\simeq\:}}
\newcommand{\compl}{\scriptscriptstyle{\complement}}

\newcommand{\abs}[1]{\left| #1 \right|}

\def\set#1{{\{ #1\}}}

\newcommand{\mat}[1]{\left(\begin{smallmatrix} #1 
\end{smallmatrix}\right)}

\newcommand{\bigmid}{\:\big|  \big.\:}
\newcommand{\pt}{\{\text{pt}\}}
\def\dot{\mathchar"013A}  
\newcommand{\hdot}{{\hskip-0.1em\raise1pt\hbox 
to0.35em{\huge $\dot$}}}

\title[Characteristic varieties and Betti numbers of 
free abelian covers]%
{Characteristic varieties and Betti numbers of\\ 
free abelian covers}

\author[Alexander~I.~Suciu]{Alexander~I.~Suciu}

\address{Department of Mathematics,
Northeastern University,
Boston, MA 02115, USA}
\email{\href{mailto:a.suciu@neu.edu}{a.suciu@neu.edu}}
\urladdr{\href{http://www.math.neu.edu/~suciu/}%
{http://www.math.neu.edu/\~{}suciu}}

\thanks{Partially supported by NSA grant H98230-09-1-0021 
and NSF grant DMS--1010298}

\subjclass[2010]{Primary 
14F35, 
55N25; 
Secondary 
20J05,  
32S22,  
57M07,  
57M12.   
}

\keywords{Free abelian cover, characteristic variety,  
exponential tangent cone, Dwyer--Fried set, special 
Schubert variety, translated subtorus, K\"{a}hler manifold, 
quasi-K\"{a}hler manifold, hyperplane arrangement, 
property $\FP_n$.}

\begin{document}

\begin{abstract}
The regular $\Z^r$-covers of a finite cell complex 
$X$ are parameterized by the Grassmannian of 
$r$-planes in $H^1(X,\Q)$. Moving about this variety, 
and recording when the Betti numbers $b_1,\dots, b_i$ 
of the corresponding covers are finite carves out certain 
subsets $\Omega^i_r(X)$ of the Grassmannian.  

We present here a method, essentially going back to 
Dwyer and Fried, for computing these sets in terms  
of the jump loci for homology with coefficients in rank~$1$ 
local systems on $X$.  Using the exponential tangent 
cones to these jump loci, we show that each 
$\Omega$-invariant is contained in the complement 
of a union of Schubert varieties associated to an 
arrangement of linear subspaces in $H^1(X,\Q)$. 

The theory can be made very explicit in the case 
when the characteristic varieties of $X$ are unions 
of translated tori.  But even in this setting, the 
$\Omega$-invariants are not necessarily open, not even 
when $X$ is a smooth complex projective variety.  
As an application, we discuss the geometric finiteness 
properties of some classes of groups.
\end{abstract}

\maketitle
\tableofcontents

\section{Introduction}
\label{sect:intro}

In a short yet insightful paper \cite{DF}, W.~Dwyer and D.~Fried 
showed that the support loci of the Alexander invariants 
of a finite CW-complex completely determine the (rational) 
homological finiteness properties of its regular, free abelian 
covers.  This result was recast in \cite{PS-plms} in terms 
of the characteristic varieties of the given CW-complex,  
as well as their exponential tangent cones. 

Our goal here is to lay the foundations of this theory in 
some detail, develop the machinery to a fuller extent, present 
a number of new results and applications, especially in the setting 
of smooth, complex quasi-projective varieties, and indicate 
some further directions of study. 

\subsection{Characteristic varieties}
\label{subsec:intro cv}

The origins of the subject go back to the 1920s, when  
J.W. Alexander introduced his eponymous knot polynomial. 
Let $K$ be a knot in $S^3$, let $X=S^3\setminus K$ 
be its complement, and let $X^{\ab}\to X$ be the 
universal abelian cover. Then, the first homology group 
$H_1(X^{\ab},\Z)$ is a finitely presented, torsion module 
over the group-ring $\Z[H_1(X,\Z)]\cong \Z[t^{\pm 1}]$;  
the Alexander polynomial of the knot, $\Delta_K(t)$, is simply 
the order of this module. More generally, if $L$ is an $n$-component 
link in $S^3$, there is a multi-variable Alexander polynomial,  
$\Delta_L(t_1, \dots , t_n)$, depending only on the link 
group, and a choice of meridians. 

Even more generally, let $X$ be a connected CW-complex 
with finite $k$-skeleton, for some $k\ge 1$.  Each 
homology group $H_j(X^{\ab},\C)$ with $j\le k$ is a finitely 
generated module over the Noetherian ring $\C[G_{\ab}]$, 
where $G=\pi_1(X,x_0)$. Consider the support loci for the 
direct sum of these modules, up to a fixed degree $i\le k$: 
\begin{equation}
\label{eq:intro  vann}
\VV^i(X) = V\Big(\!\ann\Big(\bigoplus\nolimits_{j\le i} 
H_j\big(X^{\ab}, \C\big)\Big)\!\Big).
\end{equation}

By construction, these loci are Zariski closed subsets of 
the character group  $\wG=\Hom(G,\C^{\times})$. 
As shown in \cite{PS-plms}, the varieties defined 
in this fashion may be reinterpreted as the jump 
loci for the homology of $X$ with coefficients in 
rank~$1$ local systems:
\begin{equation}
\label{eq:intro cvs}
\VV^i(X)=\big\{\rho \in \wG \bigmid H_j(X,\C_{\rho})\ne 0, 
\text{ for some $j\le i$}\big\}.
\end{equation}

It turns out that the geometry of these {\em characteristic varieties}\/ 
is intimately related to the homological properties of regular, 
abelian covers of $X$.  For the purpose of studying free 
abelian covers (like we do here), it will be enough to consider 
the subvarieties $\WW^i(X)$, obtained by intersecting 
$\VV^i(X)$ with the identity component of the character 
group, $\widehat{G}^{\circ}$. 

\subsection{The Dwyer--Fried sets}
\label{subsec:intro df}

For a fixed positive integer $r$, the connected, regular covers 
$Y\to X$ with group of deck-transformations free abelian 
of rank $r$ are parameterized by the Grassmannian of $r$-planes 
in the rational vector space $H^1(X,\Q)$. Moving about 
this variety, and recording when all the Betti numbers 
$b_1(Y),\dots, b_i(Y)$ are finite defines subsets 
\begin{equation}
\label{eq:intro df}
\Omega^i_r(X)\subseteq \Grass_r(H^1(X,\Q)),
\end{equation}
which we call the {\em Dwyer--Fried invariants}\/ of $X$. 
As with the characteristic varieties, 
these sets depend only on the homotopy type of the 
given CW-complex $X$. Consequently, if $G$ is a finitely 
generated group, we may define $\Omega^i_r(G):=
\Omega^i_r(K(G,1))$. 

In \cite{DF}, Dwyer and Fried showed that the $\Omega$-sets 
\eqref{eq:intro df} are completely determined by the 
support varieties \eqref{eq:intro  vann}.  Building on this 
work, and on the refinements from \cite{PS-plms}, we 
recast this result in Theorem \ref{thm:df cv}, as follows.  
Identify the character group $\wG$ with $H^1(X,\C^{\times})$, 
and let $\exp\colon H^1(X,\C)\to H^1(X,\C^{\times})$ be the 
coefficient homomorphism induced by the exponential map 
$\exp\colon \C\to \C^{\times}$. Then,
\begin{equation}
\label{eq:intro grass cv}
\Omega^i_r(X)=\big\{P\in \Grass_r(H^1(X,\Q)) \bigmid
\dim_{\C} ( \exp(P \otimes \C) \cap \WW^i(X)) = 0 \big\}.
\end{equation}
In other words, if $P$ is an $r$-plane in $H^1(X,\Q)$ 
and $Y\to X$ is the corresponding $\Z^r$-cover, then 
the first $i$ Betti numbers of $Y$  are finite if and only if 
the intersection of the algebraic 
torus $\exp(P \otimes \C)$ with the characteristic variety 
$\WW^i(X)$ is finite. 

\subsection{An upper bound for the $\Omega$-sets}
\label{subsec:intro schubert}

Each characteristic variety $\WW^i(X)$ determines an   
arrangement $\CC_i(X)$ of rational subspaces in $H^1(X,\Q)$.  
Indeed, let $\tau_1(\WW^i(X))$ be the set of points $z\in H^1(X,\C)$ 
such that $\exp(\lambda z)$ belongs to $\WW^i(X)$, for all 
$\lambda\in \C$.  As shown in \cite{DPS-duke}, and as 
reproved in detail here in Theorem \ref{thm:tau1}, this set 
is a finite union of rationally defined subspaces. We then 
simply declare that $\tau_1( \WW^i(X)) =  
\bigcup_{L\in \CC_i(X)} L\otimes \C$. 

Using this notion, we establish in Theorems \ref{thm:tau bound} 
and \ref{thm:sch bound} the following upper bound for the 
Dwyer--Fried sets of our space $X$: 
\begin{equation}
\label{eq:intro schubert}
\Omega^i_r(X) \subseteq  \Grass_r(H^1(X,\Q)) \setminus 
\bigcup_{L\in \CC_i(X)} \sigma_r(L ),
\end{equation}
where $\sigma_r(L)$ is the variety of incident 
$r$-planes to $L$.  Thus, each 
set $\Omega^i_r(X)$ is contained in the complement of 
a Zariski closed subset of $\Grass_r(H^1(X,\Q))$, namely, 
the union of the special Schubert varieties $\sigma_r(L)$ 
corresponding to the subspaces $L$ in $\CC_i(X)$. 

If $r=1$, inclusion \eqref{eq:intro schubert} holds as 
equality; thus, the sets $\Omega^i_1(X)$ are Zariski 
open subsets of $\Grass_1(H^1(X,\Q))$. For $r>1$, though, 
the sets $\Omega^i_r(X)$ need not be open, not even in the 
usual topology on the Grassmannian.  This rather surprising 
fact was first noticed by Dwyer and Fried, who constructed 
in \cite{DF} a $3$-dimensional cell complex for which 
$\Omega^2_2(X)$ is a finite set.   We take here a different 
approach, by analyzing the openness of the $\Omega$-sets 
in a particularly simple, yet still intriguing case.  

\subsection{Translated tori}
\label{subsec:intro tt}
This is the case when all positive-dimensional components 
of the characteristic variety $\WW^i(X)$ are translated 
subtori of the character torus  $\widehat{G}^{\circ}$.  Write 
each such component as $\rho_{\alpha} T_{\alpha}$, where 
$T_{\alpha}=\exp(L_{\alpha}\otimes \C)$, for some 
linear subspace $L_{\alpha}\subset H^1(X,\Q)$ and 
some $\rho_{\alpha} \in  \widehat{G}^{\circ}$. Using 
the intersection theory of translated subgroups in a 
complex algebraic torus developed in \cite{Hi, SYZ2}, 
we prove in Theorem \ref{thm:df translated} that 
\begin{equation}
\label{eq:intro df translated}
\Omega^i_r(X) =
\Grass_r(H^1(X,\Q)) \setminus 
\bigcup_{\alpha}  \sigma_r(L_\alpha, \rho_\alpha), 
\end{equation}
where $\sigma_r(L_\alpha,\rho_\alpha)$ is the set of 
rational $r$-planes $P$ for which 
$\rho_\alpha\in \exp((P+L_\alpha)\otimes \C)$ 
and $P\cap L_\alpha\ne \{0\}$.  
These sets, which generalize the usual Schubert varieties 
$\sigma_r(L_\alpha)$, need not be Zariski closed. 

One instance when this happens 
is described in Theorem  \ref{thm:df tors translated}.  
Suppose there is a component $\rho_{\beta} T_{\beta}$ 
of dimension $d\ge 2$ such that it, and all other components parallel 
to it have non-trivial, finite-order translation factors, while the 
other components $\rho_{\alpha} T_{\alpha}$ satisfy 
$\tau_1(T_{\alpha})\cap \tau_1(T_\beta)=\{0\}$.  Then, 
for each $2\le r\le d$, the set $\Omega^i_r(X)$ 
is {\em not}\/ an open subset of $\Grass_r(H^1(X,\Q))$.

We illustrate this phenomenon in Example \ref{ex:omega closed}, 
where we construct a finitely presented  group $G$ for which 
$\Omega^{1}_{2}(G)$ consists of a single point.  As far as we know, 
this is the first example of a group for which one of the Dwyer--Fried 
sets is not open. 

\subsection{The Green--Lazarsfeld sets and the $\Omega$-sets}
\label{subsec:intro proj}

Perhaps the best-known class of spaces for which the 
characteristic varieties consist only of translated tori is 
that of compact K\"{a}hler manifolds.  The basic 
structure of the characteristic varieties of these manifolds  
was determined by Green and Lazarsfeld \cite{GL87, GL91}, 
building on work of Beauville \cite{Be92} and Catanese \cite{Cat91}. 
The theory was further amplified in \cite{Si93, EL97, Ar, Cm01}, 
with some of the latest developments appearing in 
\cite{Di07, Dz08, DPS-duke, Bu09, Cm11, ACM}. 

In particular, suppose $M$ is a compact K\"{a}hler manifold, and $T$ is 
a positive-dimensional component of $\VV^1(M)$.  There is then  
an orbifold fibration $f\colon M\to \Sigma_{g,\mathbf{m}}$  
with either base genus $g\ge 2$, or with $g=1$ and non-trivial 
multiplicity vector $\mathbf{m}$, so that $T$ is a component 
of the pullback along $f_{\sharp} \colon \pi_1(M) \to \Gamma$ 
of the first characteristic variety of the orbifold fundamental group 
$\Gamma=\pi_1^{\orb}(\Sigma_{g,\mathbf{m}})$.  Now, it is 
readily seen that $\VV^1(\Gamma)=\widehat{\Gamma}$ or 
$\widehat{\Gamma} \setminus \widehat{\Gamma}^{\circ}$, 
depending on whether the base genus is at least $2$ or not, 
and this finishes the description of the positive-dimensional 
components of $\VV^1(M)$. 

This description allows us to either estimate or compute 
explicitly the degree one Dwyer--Fried sets of $M$.  
For instance, we show in Theorem \ref{thm:df kahler} that 
\begin{equation}
\label{eq:intro df kahler}
\Omega^1_r(M) \subseteq \Grass_r(H^1(M,\Q))\setminus \bigcup_{\alpha} 
\sigma_r(f_{\alpha}^* (H^1(\Sigma_{g_{\alpha}},\Q)),
\end{equation}
where the union runs over the set of 
orbifold fibrations $f_{\alpha}\colon M\to \Sigma_{g_{\alpha}}$ 
with $g_{\alpha}\ge 2$.  Moreover, if $r=1$, or if there are  
no orbifold fibrations with multiple fibers, then 
\eqref{eq:intro df kahler} holds as equality.

In general, though, the above inclusion is strict.  For instance, 
suppose $M$ is a smooth, complex projective variety, 
admitting an elliptic pencil with multiple fibers.  
We then show in Proposition \ref{prop:dfk not open} that 
$\Omega^1_2(M)$ is {\em not}\/ an open subset of 
the Grassmannian.  A concrete example of this 
phenomenon is provided by the 
Catanese--Ciliberto--Mendes Lopes surface (the quotient 
of $\Sigma_2\times \Sigma_3$ by a certain free involution), 
which fibers over $\Sigma_1$ with two multiple fibers, 
each of multiplicity $2$.

\subsection{Quasi-projective varieties and hyperplane arrangements}
\label{subsec:intro qp arr}

Much of this theory generalizes to quasi-K\"{a}hler manifolds.  
Indeed, if $X=\overline{X}\setminus D$ is such a manifold, 
obtained from a compact K\"{a}hler manifold $\overline{X}$ 
with $b_1(\overline{X})=0$ by removing a normal-crossings 
divisor $D$, then a theorem of Arapura \cite{Ar} insures that 
all components of $\VV^i(X)$ are unitary translates of algebraic 
subtori.

Now, if $X$ is a smooth, quasi-projective variety, 
then all such subtori in $\VV^1(X)$ arise by pullback along 
certain orbifold fibrations $f\colon X\to (\Sigma_{g,s},\mathbf{m})$ 
to Riemann surfaces of genus $g\ge 0$ with $s\ge 0$ points removed. 
Using this fact, we show in Proposition \ref{prop:df quasiproj} 
that 
\begin{equation}
\label{eq:intro omega qp}
\Omega^1_r(X) \subseteq \Grass_r(H^1(X,\Q))\setminus \bigcup_{\alpha} 
\sigma_r(f_{\alpha}^* (H^1(\Sigma_{g_{\alpha},s_{\alpha}},\Q)),
\end{equation}
where the union runs over the set of orbifold fibrations 
$f_{\alpha}\colon X\to \Sigma_{g_{\alpha},s_{\alpha}}$ 
with $2g_{\alpha}+s_{\alpha} \ge 3$.  As before, if $r=1$, 
or if there are no orbifold fibrations with multiple fibers, 
then \eqref{eq:intro omega qp} holds as equality.
On the other hand, we prove in Corollary \ref{cor:tt}
the following: If there is a pencil $X\to \C^{\times}$ 
with multiple fibers, but there is no pencil 
$X\to \Sigma_{g,s}$ with $b_1(X)-b_1(\Sigma_{g,s})\le 1$, 
then inclusion  \eqref{eq:intro omega qp} is strict for $r=2$.

In the case of hyperplane arrangements, the whole theory 
becomes much more specific and combinatorial in nature.  
If $\A$ is an arrangement of $n$ hyperplanes in $\C^{d}$,   
then its complement, $X(\A)$, is a smooth, 
quasi-projective variety.  Its first characteristic variety, 
$\VV^1(X(\A))$, is a union of torsion-translated subtori 
in $(\C^{\times})^n$. As shown by Falk and Yuzvinsky in \cite{FY}, 
the components passing through the origin correspond to 
multinets on the intersection lattice of $\A$. The other 
components are either torsion-translates of the former, 
isolated torsion points, or $1$-dimensional translated 
subtori arising from orbifold fibrations $X(\A)\to \C^{\times}$ 
with at least one multiple fiber.  

This description of the first characteristic variety of an arrangement 
leads to the following combinatorial upper bound for the degree one  
Dwyer--Fried sets:
\begin{equation}
\label{eq:intro dfarr}
\Omega^1_r(X(\A)) \subseteq \Grass_r(\Q^n)\setminus \bigcup_{\alpha} 
\sigma_r(L_{\alpha}),
\end{equation}
where the union runs over the set of orbifold fibrations 
$f_{\alpha}\colon X(\A) \to \Sigma_{0,s_{\alpha}}$, $s_{\alpha}\ge 3$,  
corresponding to multinets on $L(\A)$, and $L_{\alpha}=
f_{\alpha}^* (H^1(\Sigma_{0,s_{\alpha}},\Q))$.  
For certain classes of arrangements (e.g., line arrangements 
for which one or two lines contain all the intersection points 
of multiplicity $3$ and higher), equality holds in \eqref{eq:intro dfarr}.  
In the case of the deleted $\operatorname{B}_3$ arrangement, 
though, the presence of a $1$-dimensional translated subtorus 
in the first characteristic variety forces a strict inclusion 
for $r=2$.

\subsection{Geometric finiteness}
\label{subsec:intro finiteness}

To conclude, we investigate the relationship between the 
Dwyer--Fried invariants and other, well-known finiteness 
properties of groups, namely, Wall's property $\FF_k$ and 
Serre and Bieri's property $\FP_k$.  Using an idea that 
goes back to J.~Stallings' pioneering work on the subject, \cite{Sta}, 
we show that the $\Omega$-sets provide some useful information 
about those properties, too.  

For instance, let $\vG$ be a group of type $\FF_{k-1}$, i.e., a
group admitting a classifying space with finite $(k-1)$-skeleton, 
and suppose $\vG$ arises as the kernel of an epimorphism $G\surj \Z^r$,   
where $G$ is a finitely generated group with $\Omega^{k}_{r}(G)=\emptyset$. 
We then show in Theorem \ref{thm:fpk} that $H_k(\vG,\Z)$ is 
not finitely generated; in particular, $\vG$ is not of type $\FP_k$. 

Motivated by the long-standing Shafarevich conjecture 
on holomorphic convexity, J. Koll\'{a}r asked in \cite{Ko} 
the following question: 
Is the fundamental group of a smooth, complex projective variety 
commensurable (up to finite kernels) with a group admitting 
a smooth, quasi-projective variety as classifying space? 
In \cite{DPS-crelle}, we gave a negative answer to 
this question, by proving the following result (recorded 
here as Theorem \ref{thm:kollar answer}): 
For each $k\ge 3$, there is a smooth, complex 
projective variety $M$ of complex dimension $k-1$, such 
that the group $\Gamma=\pi_{1}(M)$ is of type $\FF_{k-1}$, 
but not of type $\FP_k$.   Using the machinery developed 
here---in particular, Theorem \ref{thm:fpk}---we end the paper 
with a shorter and more transparent proof of the main 
result from \cite{DPS-crelle}.

\subsection{Further directions}
\label{subsec:intro further}
In a companion paper \cite{Su-aspm}, we connect the Dwyer--Fried 
invariants of a space $X$ to the resonance varieties of the 
cohomology algebra $H^*(X,\C)$, under certain 
``straightness" assumptions. More precisely, suppose  
that, for each $i\le k$, all components of $\WW^{i}(X)$ 
passing through the origin of $H^1(X,\C^{\times})$ are 
algebraic subtori, and the tangent cone at $1$ to 
$\WW^i(X)$ equals the $i$th resonance variety, 
\begin{equation}
\label{eq:intro res}
\RR^{i}(X,\C)=\{a \in H^1(X,\C) \mid 
 H^j(H^*(X,\C),\cdot a) \ne 0, \text{ for some $j\le i$}\}.
\end{equation}
Then, for all $i\le k$ and $r\ge 1$, the set 
$\Omega^i_r(X)\subseteq \Grass_r(H^1(X,\Q))$ is contained in the 
complement to the incidence variety $\sigma_r(\RR^i(X,\Q))$. 
If, moreover, all positive-dimensional components 
of $\WW^{i}(X)$ pass through the origin, then this inclusion  
holds as equality.

In joint work with Y.~Yang and G.~Zhao, \cite{SYZ1}, we generalize the 
theory presented here to regular covers $Y\to X$ for which the 
group of deck-transformations is a fixed, finitely generated 
(not necessarily torsion-free)  abelian group $A$.  
Such covers are parametrized by epimorphisms from $\pi_1(X)$ 
to $A$, modulo the action of the automorphism group of $A$ on the 
target. Inside this parameter space, we single out 
the subsets $\Omega_A^i(X)$ consisting of all $A$-covers with finite 
Betti numbers up to degree~$i$.  
These sets can again be computed 
in terms of intersections of algebraic subtori with the characteristic 
varieties. For many spaces of interest, the homological finiteness 
of abelian covers can be tested through the corresponding free 
abelian covers; arbitrary abelian covers, though, exhibit different 
homological finiteness properties than their free abelian counterparts. 

Finally, in \cite{Su-pisa} we explore some of the connections between 
the Dwyer--Fried invariants $\Omega^i_r(X)\subseteq \Grass_r(H^1(X,\Q))$ 
and the Bieri--Neumann--Strebel--Renz invariants 
$\Sigma^i(X,\Z)\subseteq H^1(X,\R)$.  
In \cite{PS-plms}, we showed that 
\begin{equation}
\label{eq:bnsr bound}
\Sigma^i(X,\Z)\subseteq  H^1(X,\R) \setminus 
\bigcup_{L\in \CC_i(X)} L\otimes \R.
\end{equation}

This (sharp) upper bound for the $\Sigma$-invariants may be 
viewed as the analogue of the upper bound \eqref{eq:intro schubert}  
for the $\Omega$-invariants.  Building on the work presented here, 
we prove in \cite{Su-pisa} the following theorem:  If equality 
holds in \eqref{eq:bnsr bound}, then equality must 
also hold in \eqref{eq:intro schubert}.  This result allows us to 
derive useful information about the notoriously intricate  
$\Sigma$-invariants from concrete knowledge of the more 
accessible $\Omega$-invariants.  For instance, a basic 
question about the nature of the $\Sigma$-invariants of 
arrangements was posed at an Oberwolfach Mini-Workshop 
in 2007, and was later revisited in \cite{PS-plms, Su11}. 
Drawing on Example \ref{ex:df deleted B3}, 
we provide in \cite{Su-pisa} an answer to this question, 
thereby laying the ground for further study of the homological 
and geometric finiteness properties of arrangements. 

\section{Equivariant chain complex and characteristic varieties}
\label{sect:cvs}

We start by reviewing the definition and basic properties of the 
characteristic varieties attached to a space.

\subsection{Homology in rank~$1$ local systems}
\label{subsec:loc syst}   

Let $X$ be a connected CW-complex, with finite 
$1$-skeleton.  Without loss of generality, we may assume 
our space has a single $0$-cell, call it $x_0$.  Let $G=\pi_1(X,x_0)$ 
be the fundamental group of $X$, based at $x_0$. Clearly,  
the group $G$ is finitely generated (by the homotopy classes of the 
$1$-cells of $X$).  

Denote by $\C^{\times}$ the multiplicative group of non-zero 
complex numbers. The set of complex-valued characters, 
$\wG=\Hom(G,\C^{\times})$, is a commutative, affine algebraic 
group, with multiplication $\rho_1\cdot \rho_2 (g) 
= \rho_1(g)\rho_2(g)$, and identity the homomorphism 
taking constant value $1\in \C^{\times}$.  

Let $G_{\ab}=G/G'=H_1(X,\Z)$ be the maximal abelian 
quotient of $G$. Since the group $\C^{\times}$ is commutative, 
every character $\rho\colon G\to \C^{\times}$ factors through 
$G_{\ab}$. Thus, the abelianization map $\ab\colon G\to G_{\ab}$ 
induces an isomorphism 
$\widehat{\ab}\colon \wG_{\ab}\isom \wG$,   
which allows us to identify the coordinate ring of $\wG$ 
with the group algebra $\C[G_{\ab}]$. 
Write  $G_{\ab}=\Z^n \oplus A$, where $n=b_1(G)$ and 
$A$ is a finite abelian group.  The identity component 
of the character group, $\wG^{\circ}$, is isomorphic to the 
complex algebraic torus $(\C^{\times})^n$, while $\wG$ is 
isomorphic to the product $(\C^{\times})^n\times \widehat{A}$, 
where $\widehat{A}\cong A$ is a subgroup of $\C^{\times}$, 
consisting of roots of unity.  In particular, all components 
of $\wG$ are $n$-dimensional algebraic tori. 

The set $\wG$ parametrizes rank~$1$ local systems 
on $X$: given a character $\rho$, denote by $\L_{\rho}$ 
the complex vector space $\C$, viewed as a 
right module over the group ring $\Z{G}$ via 
$a \cdot g = \rho(g)a$, for $g\in G$ and $a\in \C$. 

Let $p\colon \wX \to X$ be the universal cover.  The cell 
structure on $X$ lifts in a natural fashion to a cell structure 
on $\wX$.  Fixing a lift $\tilde{x}_0\in p^{-1}(x_0)$ identifies 
the group $G=\pi_1(X,x_0)$ with the group of deck 
transformations of $\wX$. Therefore, we may view 
the cellular chain complex of $\wX$, 
\begin{equation}
\label{eq:equiv cc}
\xymatrixcolsep{22pt}
\xymatrix{\cdots \ar[r]& 
C_{i+1}(\wX,\Z) \ar^(.53){\tilde{\partial}_{i+1}}[r] & 
 C_{i}(\wX,\Z) \ar^(.45){\tilde{\partial}_{i}}[r] & 
  C_{i-1}(\wX,\Z) \ar[r] & \cdots  
},
\end{equation}
as a chain complex of left $\Z{G}$-modules.
The homology of $X$ with coefficients 
in $\L_{\rho}$ is defined as the homology of the 
chain complex $\L_{\rho} \otimes_{\Z{G}}  C_{\hdot}(\wX,\Z)$. 
In concrete terms, $H_*(X,\L_{\rho})$ may be computed 
from the chain complex of $\C$-vector spaces,
\begin{equation}
\label{eq:eval cc}
\xymatrixcolsep{14pt}
\xymatrix{\cdots \ar[r]& C_{i+1}(X,\C) 
\ar^{\tilde\partial_{i+1}(\rho)}[rr] &&
C_{i}(X,\C) \ar^(.47){\tilde\partial_{i}(\rho)}[rr] &&
C_{i-1}(X,\C) \ar[r] & \cdots  
},
\end{equation}
where the evaluation of $\tilde\partial_i$ at $\rho$ 
is obtained by applying the ring homomorphism 
$\Z{G}\to \C$, $ g\mapsto \rho(g)$ 
to each entry of $\tilde\partial_i$.

Alternatively, we may consider the universal abelian 
cover, $X^{\ab}$, and its equivariant chain 
complex, $C_{\hdot}(X^{\ab},\Z)=
\Z{G_{\ab}} \otimes_{\Z{G}} C_{\hdot}(\wX,\Z)$, 
with differentials $\partial^{\ab}_i=\id\otimes\, 
\widetilde{\partial}_i$.  The homology of $X$ with 
coefficients in the rank~$1$ local system given by 
a character $\rho\in \wG_{\ab}=\wG$ is then computed from the 
chain complex \eqref{eq:eval cc}, with differentials 
$\partial^{\ab}_{i}(\rho)=\tilde\partial_{i}(\rho)$.

We will write $b_i(X,\rho)=\dim_{\C} H_i(X,\C_{\rho})$. 
Evidently, the identity $1\in \wG$ yields the trivial 
local system, $\L_{1}=\C$; thus, $H_*(X,\C)$ is the usual 
homology of $X$ with coefficients in $\C$, and 
$b_i(X)=b_1(X,1)$ is the $i$th Betti number of $X$.

\subsection{Jump loci for twisted homology}
\label{subsec:jumps} 

Computing the homology groups of $X$ with coefficients 
in rank~$1$ local systems carves out a notable collection 
of subsets of the character group $\wG$. 

\begin{definition}
\label{def:cvs}
Let $X$ be a connected CW-complex with finite $k$-skeleton, 
for some $k\ge 1$. 
The {\em characteristic varieties}\/ of $X$ are the sets 
\begin{equation*}
\label{eq:cvs}
\VV^i_d(X)=\{\rho \in \wG 
\mid \dim_{\C} H_i(X,\L_{\rho})\ge d\},
\end{equation*}
defined for all degrees $0\le i\le k$ and all depths $d>0$. 
\end{definition}

\begin{remark}
\label{rem:finite cv}
For the purpose of computing the characteristic varieties 
up to degree $i=k$, we may assume without loss of generality 
that $X$ is a finite CW-complex of dimension $k+1$ (see 
\cite[Lemma 2.1]{PS-plms}).
\end{remark}

The terminology from Definition \ref{def:cvs}, due to 
Libgober, is justified by the following lemma.  For 
concreteness, we sketch a proof, based on an 
argument of Green and Lazarsfeld \cite{GL87} 
(see also \cite{Li02}).
Given a commutative ring $R$, and a matrix 
$\varphi\colon R^b\to R^a$, denote by $E_q(\varphi)$ 
the ideal generated by the minors of size $a-q$ of $\varphi$.

\begin{lemma}
\label{lem:cv ideal}
The jump loci $\VV^i_d(X)$ are Zariski closed 
subsets of the algebraic group $\wG$.
\end{lemma}

\begin{proof}
Let $R=\C[G_{\ab}]$ be the coordinate ring of $\wG=\wG_{\ab}$.  
By definition, a character $\rho$ belongs to $\VV^i_d(X)$ 
if and only if $\rank \partial^{\ab}_{i+1}(\rho) + 
\rank \partial^{\ab}_{i}(\rho) \le c_i -d$, 
where $c_i=c_i(X)$ is the number of $i$-cells of $X$.  
Using this description, we may rewrite our jump locus as 
the zero-set of a sum of products of determinantal ideals,
\begin{equation}
\label{eq:ideal cv}
\VV^i_d(X) = V\Big( \sum 
E_p(\partial^{\ab}_i) \cdot E_q(\partial^{\ab}_{i+1}) \Big), 
\end{equation}
with the sum running over all $p,q\ge 0$ with 
$p+q=c_{i-1}+d-1$.
\end{proof}

\subsection{Homotopy invariance}
\label{subsec:cv he} 

Let $X$ and $Y$ be connected CW-complexes with finite $k$-skeleta,   
and denote by $G$ and $H$ the respective fundamental groups.

\begin{lemma}
\label{lem:cv inv}
Suppose $X\simeq Y$. There is then an 
isomorphism $\widehat{H} \cong \wG$, which restricts 
to isomorphisms $\VV^i_d(Y) \cong \VV^i_d(X)$, for 
all $i\le k$, and all $d>0$. 
\end{lemma}

\begin{proof}
Let $f\colon X\to Y$ be a homotopy equivalence;  
without loss of generality, we may assume $f$ 
preserves skeleta.  The induced homomorphism 
on fundamental groups, $f_{\sharp} \colon G \to H$, 
yields an isomorphism of algebraic groups, 
$\hat{f_{\sharp}}\colon \widehat{H} \to \wG$.  
Lifting $f$ to a cellular homotopy 
equivalence, $\tilde{f}\colon \widetilde{X} \to  \widetilde{Y}$, 
defines isomorphisms $H_i(X, \L_{\rho\circ f_{\sharp}}) 
\to  H_i(Y, \L_{\rho}) $, for each character $\rho\in \widehat{H}$. 
Hence, $\hat{f_{\sharp}}$
restricts to isomorphisms $\VV^i_d(Y) \to \VV^i_d(X)$ 
between the respective characteristic varieties.
\end{proof}

In each fixed degree $i$, the characteristic varieties of 
a space $X$ define a descending filtration on the character 
group, 
\begin{equation}
\label{eq:filt cvs}
\wG =\VV^i_0(X) \supseteq \VV^i_1(X) 
\supseteq \VV^i_2(X) \supseteq \cdots.
\end{equation}
Clearly, $1\in \VV^i_d(X)$ if and only if $d\le b_i(X)$. 
Moreover, if $c_i(X)=0$, then $\VV^i_d(X)=\emptyset$, 
for all $d$.  If particular, if $X$ has dimension $k$, then 
$\VV^i_d(X)=\emptyset$, for all $i>k$, and all $d>0$.

In degree $0$, we have $\VV^0_1(X)= \{ 1\}$ 
and $\VV^0_d (X)=\emptyset$, for $d>1$. 
In degree $1$, the characteristic varieties  of $X$ depend only 
on the fundamental group $G=\pi_1(X,x_0)$---in fact, only on 
its maximal metabelian quotient, $G/G''$ (see \S\ref{subsec:cvs one}  
below). Accordingly, we will sometimes write $\VV^1_d(G)$ for 
$\VV^1_d(X)$.  

One may define in the same fashion the characteristic 
varieties $\VV^i_d(X,\k)$ over an arbitrary field $\k$.  
Let us just note here that $\VV^i_d(X,\k)=
\VV^i_d(X,\mathbb{K}) \cap \Hom(G,\k^{\times})$, 
for every field extension $\k\subseteq \mathbb{K}$. 

\subsection{Some basic examples}
\label{subsec:ex cv}
There are a few known instances where all the characteristic 
varieties of a given space can be computed explicitly. 
We record below some of these computations, which 
will be needed later on.

\begin{example}
\label{ex:circle}
We start with the circle $S^1$. Identify  
$\pi_1(S^1,*)=\Z$, and fix a generator $t$; 
the cellular chain complex of 
$\widetilde{S^1}$ then takes the form 
$ C_1 \xrightarrow{\partial_1} C_0 $, 
where $C_1=C_0=\Z\Z$ and $\partial_1(1)=t-1$.  
By specializing at a character $\rho\in \widehat{\Z}=\C^{\times}$, 
we get a chain complex  of $\C$-vector spaces,  
$ \C \xrightarrow{\rho-1} \C $, 
which is acyclic, except for $\rho=1$, when 
$H_0(S^1,\C)=H_1(S^1,\C)=\C$.   
Hence, $\VV^0_1(S^1)=\VV^1_1(S^1)=\{1\}$ 
and $\VV^i_d(S^1)  = \emptyset$, otherwise.
\end{example}

\begin{example}
\label{ex:cv wedge}
Let $\bigvee^{n} S^1$ be a wedge of $n$ circles, $n>1$. 
Then $\pi_1(\bigvee^{n} S^1)=F_n$, the free group of 
rank $n$.  It is readily seen that 
\begin{equation}
\label{eq:cv bouquet}
\VV^i_d\Big(\bigvee^{n} S^1\Big)=
\begin{cases}
(\C^{\times})^n & \text{if $i=1$ and $d< n$},\\
\{1\}& \text{if $i=1$ and $d= n$, or $i=0$ and $d=1$},\\
 \emptyset & \text{otherwise}.
\end{cases}
\end{equation} 
\end{example}

\begin{example}
\label{ex:cv surf}
Let $\Sigma_g$ be a compact, connected, orientable 
surface of genus $g\ge 1$, and identify 
$\Hom(\pi_1(\Sigma_g),\C^{\times})=(\C^{\times})^{2g}$.
We then have
\begin{equation}
\label{eq:cv surf}
\VV^i_d(\Sigma_g)=
\begin{cases}
(\C^{\times})^{2g} & \text{if $i=1$ and $d< 2g-1$},\\
\{1\}& \text{if $i=1$ and $d\in \{2g-1, 2g\}$, or $i\in\{0, 2\}$ and $d=1$},\\
\emptyset & \text{otherwise}.
\end{cases}
\end{equation} 
\end{example}

\subsection{Depth one characteristic varieties}
\label{subsec:cvs one} 

Most important for us will be the depth-$1$ characteristic 
varieties, $\VV^i_1(X)$, and their unions up to a fixed degree, 
$\VV^i(X)=\bigcup_{j=0}^{i} \VV^j_1(X)$.  Let $G=\pi_1(X,x_0)$.  
Clearly, 
\begin{equation}
\label{eq:depth1 cv}
\VV^i(X)=\{\rho \in \wG \mid 
b_j(X,{\rho}) \ne 0, \text{ for some $j\le i$}\},
\end{equation}
These varieties yield an ascending filtration of the character group, 
$\{1\}=\VV^0(X) \subseteq \VV^1(X) \subseteq  \cdots \subseteq \VV^k(X) 
 \subseteq \wG$.

It follows from \cite[Corollary 3.7]{PS-plms}, that the sets 
\eqref{eq:depth1 cv} are the support varieties for the Alexander 
invariants of $X$.  More precisely, view $H_*(X^{\ab},\C)$ as 
a module over the group-ring $\C{G_{\ab}}$.  Then,
\begin{equation}
\label{eq:v ann}
\VV^i(X) = V\Big(\!\ann\Big(\bigoplus_{j\le i} H_j\big(X^{\ab},
\C\big)\Big)\!\Big).
\end{equation}

In particular, $\VV^1(G)= V(\ann (H_1(G',\Z)\otimes \C)$, 
where the $\Z{G_{\ab}}$-module structure on the group 
$H_1(G',\Z)=G'/G''$ arises from the extension  
$0\to  G'/G'' \to G/G'' \to G'/G''\to 0$.
This shows that the characteristic variety $\VV^1(G)$ does indeed 
depend only on the maximal metabelian quotient $G/G''$.

We will also consider the varieties $\WW^i(X)=\VV^i(X)\cap \wG^{\circ}$  
inside the complex algebraic torus  $\wG^{\circ}$. 
An alternate description of these varieties is as follows. 
Let $X^{\alpha}\to X$ be the maximal torsion-free abelian cover 
of $X$, corresponding to the projection $\alpha\colon  
G\surj G_{\ab}/\Tors(G_{\ab})=\Z^n$, where $n=b_1(G)$.  
Identify $\wG^{\circ}=(\C^{\times})^n$, and the group ring 
$\C{\Z^n}$ with the Laurent polynomial ring 
$\Lambda_n=\C[t_1^{\pm 1},\dots ,t_n^{\pm 1}]$.  
Then,
\begin{equation}
\label{eq:w ann}
\WW^i(X) = V\Big(\!\ann\Big(\bigoplus_{j\le i} 
H_j\big(X^{\alpha},\C\big)\Big)\!\Big).
\end{equation}

If $X$ has finite $2$-skeleton, the set $\WW^1(X)$ 
can be computed from a finite presentation for the 
fundamental group, by means of the Fox free differential calculus. 
Suppose $G= \langle x_1, \dots, x_q \mid r_1, \dots, r_m \rangle$, 
and let $\Phi_{G}$ be the corresponding $m$ by $q$ 
Alexander matrix, with entries 
$\alpha(\partial r_i/\partial x_j )$ obtained by 
applying the ring morphism $\alpha\colon \Z{G}\to \Lambda_n$ 
to the Fox derivatives of the relators.  
The variety $\WW^1(G)=\WW^1(X)$, then, 
is defined by the vanishing of the codimension~$1$ minors 
of $\Phi_G$, at least away from the trivial character $1$.  

\begin{example}
\label{ex:cv link}
Let $L=(L_1,\dots ,L_n)$ be a link of smoothly embedded 
circles in $S^3$, with complement 
$X=S^3 \setminus \bigcup_{i=1}^n L_i$.  Choosing orientations 
on the link components yields  a preferred basis for 
$H_1(X,\Z)=\Z^n$ consisting of oriented meridians.   
Using this basis, identify $H^1(X,\C^{\times})=(\C^{\times})^n$.  
Then
\begin{equation}
\label{eq:cv link}
\WW^1(X)= 
\{\zeta\in (\C^{\times})^n \mid \Delta_L(\zeta)=0\} \cup \{1\},
\end{equation}
where $\Delta_L=\Delta_L(t_1,\dots ,t_n)$ is 
the (multi-variable) Alexander polynomial of the link. 
For details and references on this, see \cite{Su11}.
\end{example}

\subsection{Products and wedges}
\label{subsec:prod wedge}

The depth-$1$ characteristic varieties behave 
well with respect to products and wedges. To make 
this statement precise, let $X_1$ and $X_2$ be 
connected CW-complexes with finite $k$-skeleta, 
and with fundamental groups $G_1$ and $G_2$.

We start with the product $X=X_1\times X_2$.  Identify 
$G=\pi_1(X)$ with $G_1\times G_2$; also, 
$\wG=\wG_1\times \wG_2$ 
and $\wG^{\circ}=\wG^{\circ}_1\times \wG^{\circ}_2$. 

\begin{prop}[\cite{PS-plms}]
\label{prop:cv prod}
$\VV^i_1(X_1\times X_2)= 
\bigcup_{p+q=i} \VV^{p}_1(X_1) \times \VV^{q}_1(X_2)$, 
for all $i\le k$.
\end{prop}

The idea of the proof is very simple.   For each 
character $\rho=(\rho_1,\rho_2)\in \wG$, the chain complex   
$C_{\hdot}(X,\C_{\rho})$ decomposes as $C_{\hdot}(X_1,\C_{\rho_1})
\otimes_{\C} C_{\hdot}(X_2,\C_{\rho_2})$. 
Taking homology, we get $H_i(X,\C_{\rho})=\bigoplus_{s+t=i} 
H_{s}(X_1,\C_{\rho_1}) \otimes_{\C} H_{t}(X_2,\C_{\rho_2})$, 
and the claim follows. 

\begin{corollary}
\label{cor:cv prod}
For all $i\le k$, we have that 
$\VV^i(X_1\times X_2)= 
\bigcup_{p+q=i} \VV^{p}(X_1) \times \VV^{q}(X_2)$ and 
$\WW^i(X_1\times X_2)= 
\bigcup_{p+q=i} \WW^{p}(X_1) \times \WW^{q}(X_2)$. 
\end{corollary}

Now consider the wedge, $X=X_1\vee X_2$, 
taken at the unique $0$-cells.  Identify 
$G=\pi_1(X)$ with $G_1* G_2$; 
also, $\wG=\wG_1\times \wG_2$ 
and $\wG^{\circ}=\wG^{\circ}_1\times \wG^{\circ}_2$. 

\begin{prop}[\cite{PS-plms}]
\label{prop:cv wedge}
Suppose $X_1$ and $X_2$ have positive 
first Betti numbers. Then, for all $1\le i\le k$,
\[
\VV^i_1(X_1\vee X_2) = 
\begin{cases}
\wG_1 \times \wG_2 
& \text{if $i=1$,}\\[2pt]
\VV^i_1(X_1)\times \wG_2 
\cup \wG_1 \times \VV^i_1(X_2)
&\text{if $i>1$.}
\end{cases}
\]
\end{prop}

The proof uses the decomposition 
$C_{+}(X,\C_{\rho}) = C_{+}(X_1,\C_{\rho_1})
\oplus C_{+}(X_2,\C_{\rho_2})$.  Taking homology,  
we get $b_i(X,\rho)=b_i(X_1,{\rho_1}) + 
b_i(X_1,{\rho_2})+ \epsilon$, where $\epsilon=1$ 
if $i=1$, $\rho_1\ne 1$,  and $\rho_2\ne 1$, and 
$\epsilon=0$, otherwise.   The claim follows. 

\begin{corollary}
\label{cor:cv wedge}
If $X_1$ and $X_2$ have positive 
first Betti numbers, then $\VV^i(X_1\vee X_2)=\wG$ 
and $\WW^i(X_1\vee X_2)=\wG^{\circ}$, 
for all $i\le k$.
\end{corollary}

The condition that $b_1(X_1)$ and $b_1(X_2)$ be positive 
may be dropped if $i>1$, but not if $i=1$.  For instance, 
take $X_1=S^1$ and $X_2=S^2$.  Then $G_1=\Z$ and 
$G_2=\{1\}$; thus, $\wG=\C^{\times}$, yet 
$\VV^1(S^1\vee S^2)=\{1\}$. 

\subsection{A functoriality property}
\label{subsec:cv funct} 

Every group homomorphism $\varphi\colon G\to H$ 
induces a morphism between character groups, 
$\hat\varphi\colon \widehat{H} \to \wG$, 
given by $\hat\varphi (\rho)(g)=\varphi(\rho(g))$.   
Clearly, the morphism $\hat\varphi$ sends algebraic subgroups 
of $\widehat{H}$ to algebraic subgroups of $\wG$.  
Moreover, if $\varphi$ is surjective, then $\hat\varphi$ 
is injective. The next lemma indicates a partial  
functoriality property for the characteristic varieties of groups.

\begin{lemma}
\label{lem:epi cv}
Let $\varphi\colon G \surj Q$ be an epimorphism 
from a finitely generated group $G$ to a group $Q$. 
Then the induced monomorphism between character groups, 
$\hat\varphi\colon \widehat{Q} \inj \wG$, restricts to an embedding 
$\VV^1_d(Q) \inj \VV^1_d(G)$, for each $d\ge 1$.
\end{lemma}

\begin{proof}
Let $\rho\colon Q\to \C^{\times}$ be a character. 
The $5$-term exact sequence associated to the extension 
$1\to K \to G \to Q \to 1$ and the $\Z{G}$-module 
$M=\L_{\rho \circ \varphi}$ ends in $H_1(G,M) \to H_1(Q,M_K)\to 0$, 
where $M_K$ denotes the module of coinvariants under the $K$-action, 
see \cite[p.~171]{Br}.  Clearly, $M_K=\L_{\rho}$, as 
$\Z{Q}$-modules.  Hence, 
$\dim_{\C} H_1(G,\L_{\rho \circ \varphi})$ is bounded 
below by $\dim_{\C} H_1(Q,\L_{\rho})$.
Thus, if $\rho\in \VV^1_d(Q)$, then $\hat\varphi(\rho) \in \VV^1_d(G)$, 
and we are done.
\end{proof}

\section{A stratification of the rational Grassmannian}
\label{sect:df sets}

In this section, we stratify the Grassmannian of $r$-planes 
in $H^1(X,\Q)$ by certain subsets $\Omega^i_r(X)$
which keep track of those regular $\Z^r$-covers of $X$ 
having finite Betti numbers up to degree $i$.

\subsection{Free abelian covers}
\label{subsec:free abel}

As before, let $X$ be a connected CW-complex 
with finite $1$-skeleton, and let $G=\pi_1(X, x_0)$ 
be the fundamental group, based at the unique $0$-cell $x_0$.
Following Dwyer and Fried \cite{DF}, we start by parameterizing 
in a convenient way the set of all connected, regular covering 
spaces of $X$ (up to equivalence of covers), with group of 
deck transformations a free abelian group of fixed rank $r$. 

The model situation is the $r$-dimensional torus $T^r=K(\Z^r,1)$  
and its universal cover, $\R^r \to T^r$, with group of 
deck transformations $\Z^r$. 
Any epimorphism $\nu\colon G\surj \Z^r$ gives 
rise to a $\Z^r$-cover, $X^{\nu}\to X$, by pulling back 
the universal cover of $T^r$ along a map $f\colon X\to T^r$ 
realizing $\nu$ at the level of fundamental groups,
\begin{equation}
\label{eq:pull}
\xymatrix{
X^{\nu}  \ar[r] \ar[d] &\, \R^r \ar[d] \\
X \ar^{f}[r] & \:T^r. 
}
\end{equation}
Since this is a pull-back diagram, and $\R^r$ is contractible, 
the homotopy fiber of $f$ is homotopy equivalent to $X^{\nu}$. 
By covering space theory, all connected, regular 
$\Z^r$-covers of $X$ arise in the fashion described above.

Now, the homomorphism $\nu\colon G\surj \Z^r$ factors 
as $\nu_*\circ \ab$, where $\ab\colon G\to G_{\ab}$ 
is the abelianization map, and  
$\nu_{*}=\nu_{\ab}\colon G_{\ab} \to \Z^r$ 
may be identified with the induced homomorphism 
$f_*\colon H_1(X,\Z)\to H_1(T^r,\Z)$. 
Passing to the corresponding homomorphism  
in rational homology, we see that the cover 
$X^{\nu}\to X$ is determined (up to equivalence) 
by the kernel of the map $\nu_{*}\colon H_1(X,\Q) \to \Q^r$.   
Conversely, every codimension-$r$ linear subspace 
of $H_1(X,\Q)$ can be realized as the kernel of 
$\nu_{*}\colon H_1(X,\Q) \to \Q^r$, 
for some epimorphism $\nu\colon G\surj \Z^r$, 
and thus gives rise to a cover $X^{\nu}\to X$.  
 
Let $\Grass_r(H^1(X,\Q))$ be the Grassmannian of 
$r$-planes in the finite-dimensional, rational vector 
space $H^1(X,\Q)$.   Proceeding as above, but using 
the dual homomorphism $\nu^{*}\colon \Q^r\to H^1(X,\Q)$, 
instead, we obtain a one-to-one correspondence between 
equivalence classes of regular $\Z^r$-covers of $X$ 
and $r$-planes in $H^1(X,\Q)$, which we record below.

\begin{prop}[Dwyer--Fried \cite{DF}]
\label{prop:df corresp}
The connected, regular covers of $X$ whose group of 
deck transformations is free abelian of rank $r$ are 
parameterized by the rational Grassmannian 
$\Grass_r(H^1(X,\Q))$, via the correspondence 
\begin{align*}
\big\{ \text{regular $\Z^r$-covers of $X$} \big\} & \longleftrightarrow 
\big\{ \text{$r$-planes in $H^1(X,\Q)$}\big\}\\
X^{\nu}  \to X\qquad & \longleftrightarrow\qquad  P_{\nu}:=\im(\nu^*)
\end{align*}
\end{prop}

This correspondence enjoys a nice functoriality property.  
Let $f\colon (X,x_0)\to (Y,y_0)$ be a pointed map, and denote 
by $f_{\sharp}\colon \pi_1(X,x_0)\to \pi_1(Y,y_0)$ the 
induced homomorphism on fundamental groups. 

\begin{lemma}
\label{lem:df funct}
Suppose $f_{\sharp}\colon \pi_1(X,x_0)\to \pi_1(Y,y_0)$, 
is surjective.  If  $\nu\colon \pi_1(Y,y_0) \surj \Z^r$ is an 
epimorphism, then 
$P_{\nu\circ f_{\sharp}}=f^*(P_{\nu})$.
\end{lemma}

\begin{proof}
By the lifting criterion, the map $f$ lifts to a map $\bar{f}$ between 
the $\Z^r$-covers defined by the epimorphisms 
$\nu\circ f_{\sharp}$ and $\nu$.  We then have a pullback diagram,
\begin{equation}
\label{eq:lift}
\xymatrix{
X^{\nu\circ f_{\sharp}} \ar^{\bar{f}}[r] \ar[d]& Y^{\nu} \ar[d] \\
X \ar^{f}[r] & Y.}
\end{equation}

Clearly, the morphism $f_*\colon H_1(X,\Q) \to H_1(Y,\Q)$ is 
surjective.  Thus, the morphism $f^*\colon H^1(Y,\Q) \to H^1(X,\Q)$ 
is injective, and takes the $r$-plane $P_{\nu}=\im(\nu^*)$ 
to the $r$-plane $f^*(P_{\nu})=\im(f^*\circ \nu^*)$, which 
coincides with $\im((\nu\circ f_{\sharp})^*)=P_{\nu\circ f_{\sharp}}$.
\end{proof}

\subsection{The Dwyer--Fried invariants}
\label{subsec:df}

Moving about the rational Grassmannian, and keeping track  
of how the Betti numbers of the corresponding covers 
vary leads to the following definition (see \cite{DF}, 
and also \cite{PS-plms}, \cite{Su11}).

\begin{definition}
\label{def:df sets}
The {\em Dwyer--Fried invariants}\/ of $X$ are the subsets 
\begin{equation*}
\label{eq:grass}
\Omega^i_r(X)=\big\{P_{\nu} \in \Grass_r(H^1(X,\Q)) \bigmid 
\text{$b_{j} (X^{\nu}) <\infty$ for $j\le i$}  \big\},
\end{equation*}
defined for all $i\ge 0$ and all $r\ge 1$, with the convention that 
$\Omega^i_r(X)=\Grass_r(H^1(X,\Q))=\emptyset$ if $r>b_1(X)$.
\end{definition}

If $b_1(X)=0$, then all the $\Omega$-invariants 
of $X$ are empty.  
For a fixed $r\ge 1$, the Dwyer--Fried invariants form a 
descending filtration of the Grassmannian of $r$-planes, 
\begin{equation}
\label{eq:df filt}
\Grass_r(H^1(X,\Q)) = \Omega^0_r(X) \supseteq \Omega^1_r(X)  
\supseteq \Omega^2_r(X)  \supseteq \cdots,  
\end{equation}
with intersection  $\Omega_r(X)=\big\{ P_{\nu} \mid  
\dim_{\Q} H_* (X^{\nu},\Q) <\infty \big\}$.

\begin{remark}
\label{rem:df}
Dwyer and Fried only consider finite CW-complexes 
$X$, and the sets $\Omega_r(X)$.  Of course, if 
$\dim X=k$, then $\Omega_r(X) =\Omega^k_r(X)$.  
We prefer to work with 
the filtration \eqref{eq:df filt}, which provides more refined 
information on the $\Z^r$-covers of $X$. 
\end{remark}

The $\Omega$-sets are homotopy-type invariants of $X$.  
More precisely, we have the following result.

\begin{lemma}
\label{lem:df inv}
Suppose $X$ is homotopy equivalent to $Y$.  
For each $r\ge 1$, there is an isomorphism 
$\Grass_r(H^1(Y,\Q)) \cong \Grass_r(H^1(X,\Q))$  
which sends each subset $\Omega^i_r(Y)$ bijectively 
onto $\Omega^i_r(X)$.
\end{lemma}

\begin{proof}
Let $f\colon X\to Y$ be a (cellular) homotopy equivalence.  
The induced isomorphism in rational cohomology 
defines isomorphisms 
$f^*_r\colon \Grass_r(H^1(Y,\Q)) \to  \Grass_r(H^1(X,\Q))$ 
between the corresponding Grassmannians.  It remains 
to show that 
$f^*_r ( \Omega^i_r(Y)) \subseteq \Omega^i_r(X)$.  
For that, let $P\in  \Omega^i_r(Y)$, and write $P=P_{\nu}$, 
for some epimorphism $\nu\colon \pi_1(Y)\surj \Z^r$. The 
map $f$ lifts to a map $\bar{f}\colon X^{\nu\circ f_{\sharp}} \to Y^{\nu}$
as in \eqref{eq:lift}.   Clearly, $\bar{f}$ is a homotopy equivalence.   
Thus, $b_i(X^{\nu\circ f_{\sharp}})=b_i(Y^{\nu})$, and so 
$f^*_r(P_{\nu})=P_{\nu\circ f_{\sharp}}$ belongs to  $\Omega^i_r(X)$. 
\end{proof}

In view of this lemma, we may extend the definition
of the $\Omega$-sets from spaces to groups. 
Let $G$ be a finitely-generated group.  Pick 
a classifying space $K(G,1)$ with finite $k$-skeleton, 
for some $k\ge 1$.  

\begin{definition}
\label{def:df gp}
The {\em Dwyer--Fried invariants}\/ of a group $G$ are 
the subsets $\Omega^i_r(G)=\Omega^i_r(K(G,1))$ 
of $\Grass_r(H^1(G,\Q))$, defined for all $i\ge 0$ 
and $r\ge 1$.
\end{definition}
Since the homotopy type of $K(G,1)$ depends only $G$, 
the sets $\Omega^i_r(G)$ are well-defined group invariants. 

\subsection{Discussion}
\label{subsec:discuss}

Especially manageable is the situation when 
$n=b_1(X)>0$ and $r=n$.  In this case, $\Grass_n(H^1(X,\Q))=\pt$.  
Under the correspondence from Proposition \ref{prop:df corresp}, 
this single point is realized by the maximal free abelian cover, 
$X^{\alpha}\to X$, where  
$\alpha\colon  G\surj G_{\ab}/\Tors(G_{\ab})=\Z^n$ 
is the canonical projection. We then have 
\begin{equation}
\label{eq:df b1}
\Omega^i_n(X)=\begin{cases}
\pt  & \text{if $b_j(X^{\alpha})<\infty$ for all $j\le i$},\\
\emptyset  & \text{otherwise}.
\end{cases}
\end{equation}

Both situations may occur, as illustrated by a very 
simple example.  

\begin{example}
\label{ex:rn}
Let $X=S^1\vee S^k$, for some $k>1$. Then 
$X^{\alpha}$ is homotopic to a countable 
wedge of $k$-spheres.  Thus, 
$\Omega^i_1(X)=\pt$ for $i<k$, yet   
$\Omega^i_1(X)=\emptyset$, for $i\ge k$. 
\end{example}

It should be emphasized that finiteness of the Betti numbers 
of a free abelian cover $X^{\nu}$ does not necessarily imply 
finite-generation of the integral homology groups of $X^{\nu}$. 
Thus, we cannot replace in Definition \ref{def:df sets} the condition 
``$b_i(X^{\nu})<\infty$, for $i\le q$" by the (stronger) condition 
``$H_i(X^{\nu},\Z)$ is a finitely-generated group, for $i\le q$."  
Example \ref{ex:milnor} below (extracted from a paper of 
Milnor \cite{Mil}) explains why. 

\begin{example}
\label{ex:milnor} 
Let $K$ be a knot in $S^3$, with complement $X=S^3\setminus K$, 
and infinite cyclic cover $X^{\ab}$. As is well-known, 
$H_1(X^{\ab},\Z)=\Z[t^{\pm 1}]/(\Delta_K)$, 
where $\Delta_K$ is the Alexander polynomial of $K$. 
Hence, $H_1(X^{\ab},\Q)=\Q^d$, where $d$ is the 
degree of $\Delta_K$, and so $\Omega^1_1(X)=\pt$.  
On the other hand, if the Alexander polynomial is not monic, 
$H_1(X^{\ab},\Z)$ need not be finitely generated as a $\Z$-module. 
For instance, let $K$ be the $5_2$ knot, with Alexander 
polynomial $\Delta_K=2t^2-3t+2$.  
Then $H_1(X^{\ab},\Z)=\Z[1/2]\oplus \Z[1/2]$ is 
not finitely generated, though, of course, 
$H_1(X^{\ab},\Q)=\Q\oplus \Q$.
\end{example}

\textheight8.1in
\section{Dwyer--Fried invariants from characteristic varieties}
\label{sect:df cv}

In this section, we explain how the characteristic 
varieties of a space $X$ determine the $\Omega$-invariants 
of $X$, and therefore control the homological 
finiteness properties of its (regular) free 
abelian covers. 

\subsection{Homological finiteness and characteristic varieties}
\label{subsec:hf cv}

Given an epimorphism $\nu\colon G\surj \Z^r$, 
let $\hat\nu\colon \widehat{\Z^r} \inj \wG$  be 
the induced homomorphism between character groups. 
Its image, $\T_{\nu}=\hat\nu\big(\widehat{\Z^r}\big)$, 
is a complex algebraic subtorus of $\wG$, isomorphic 
to $(\C^{\times})^r$. 

The following theorem was proved by Dwyer and Fried in \cite{DF} 
for finite CW-com\-plexes, using the support loci for the 
Alexander invariants of such spaces.  It was recast in a slightly more 
general context by Papadima and Suciu in \cite{PS-plms}, 
using the degree-$1$ characteristic varieties.

\begin{theorem}[\cite{DF}, \cite{PS-plms}]
\label{thm:df cvar}
Let $X$ be a connected CW-complex with finite $k$-skeleton, 
for some $k\ge 1$, and let $G=\pi_1(X)$.  For an epimorphism 
$\nu\colon G\surj \Z^r$, the following are equivalent:
\begin{enumerate}
\item 
The vector space 
$\bigoplus_{i=0}^{k} H_{i} (X^{\nu}, \C)$ 
is finite-dimensional. 
\item 
The algebraic torus $\T_{\nu}$ intersects the variety $\VV^k(X)$ 
in only finitely many points. 
\end{enumerate}
\end{theorem}

In other words, whether all the Betti numbers 
$b_1(X^{\nu}),\dots , b_k(X^{\nu})$ are finite or not 
is dictated by whether the variety 
$\T_{\nu}\cap \VV^k(X)$ has dimension $0$ or not. 
More generally, it is shown in \cite{SYZ1} that the 
varieties $\VV^k(X)$ control the finiteness of the Betti numbers 
of all abelian (not necessarily torsion-free) regular covers 
of the space $X$.

\subsection{Reinterpreting the $\Omega$-invariants}
\label{subsec:dfcv}
Let $\exp\colon H^1(X, \C) \to H^1(X, \C^{\times})$ 
be the coefficient homomorphism induced by the homomorphism 
$\C\to \C^{\times}$, $z\mapsto e^z$.  

\begin{lemma}
\label{lem:exp}
Let $\nu\colon G\surj \Z^r$ be an epimorphism. 
Under the universal coefficient isomorphism 
$H^1(X,\C^{\times}) \cong\Hom(G,\C^{\times})$, 
the complex $r$-torus $\exp(P_{\nu}\otimes \C)$ 
corresponds to $\T_{\nu}=\hat\nu\big(\widehat{\Z^r}\big)$. 
\end{lemma}

\begin{proof}
By naturality of the coefficient homomorphism and 
the universal coefficients isomorphism, the following 
diagram commutes.
\[
\xymatrixcolsep{7pt}
\xymatrixrowsep{22pt}
\xymatrix{
&\Q^r \ar[rr]  \ar^(.4){\nu^*}[rr] \ar@{^{(}->}[d] 
&& H^1(X,\Q) \ar@{^{(}->}[d] \\
\Hom(\Z^r, \C) \ar@{=}^(.7){\sim}[r] \ar_{\Hom(\_ \, , \exp)}[d] 
& \C^r  \ar^(.4){\nu^*}[rr] \ar^{\exp}[d] \ar@{-->}[rrd]
&&  H^1(X,\C) \ar@{=}^(.47){\sim}[r] \ar^{\exp}[d] 
& \Hom(G,\C) \ar^{\Hom(\_ \, , \exp)}[d] \\
\Hom(\Z^r, \C^{\times}) \ar@{=}^(.65){\sim}[r]  
\ar@/_18pt/[rrrr]_{\hat\nu=\Hom(\nu,\_)}
& (\C^{\times})^r \ar^(.4){\nu^*}[rr] 
&& H^1(X,\C^{\times})\ar@{=}^(.47){\sim}[r]  & \Hom(G,\C^{\times}).
}
\] 

By definition, $P_{\nu}$ is the image of the top $\nu^*$ map. 
Hence, $P_{\nu}\otimes \C$ is the image of the middle $\nu^*$ 
map, while $\exp(P_{\nu}\otimes \C)$ is the image of the shaded 
arrow.  By surjectivity of $\exp\colon \C^r \to (\C^{\times})^r$, 
this last image is the same as the image of the bottom 
$\nu^*$ map, which in turn corresponds to 
$\T_{\nu}=\im(\hat\nu)$ under the bottom-right  isomorphism. 
\end{proof}

Recall we introduced in \S\ref{subsec:cvs one} 
the varieties $\WW^k(X)=\VV^k(X)\cap \wG^{\circ}$, 
lying in the identity component of the character group 
$\wG=H^1(G,\C^{\times})$.  

\begin{theorem}
\label{thm:df cv}
Let $X$ be a CW-complex with finite $k$-skeleton. 
Then,  for all $i\le k$ and $r\ge 1$, 
\begin{equation*}
\label{eq:grass cv}
\Omega^i_r(X)=\big\{P\in \Grass_r(H^1(X,\Q)) \bigmid
\# \big(\!\exp(P \otimes \C) \cap \WW^i(X) 
\big)<\infty \big\}.
\end{equation*}
\end{theorem}

\begin{proof}
Let $P$ be a rational $r$-plane in $H^1(X,\Q)$.  
By Proposition \ref{prop:df corresp}, there is 
an epimorphism $\nu\colon G\surj \Z^r$ such that 
$P=P_{\nu}$.  Clearly, the algebraic torus $\T_{\nu}$ 
lies in $\wG^{\circ}$; thus, 
$T_{\nu}\cap \VV^i(X)=T_{\nu}\cap \WW^i(X)$.  
Applying Theorem \ref{thm:df cvar} 
and Lemma \ref{lem:exp} finishes the proof.
\end{proof}

In other words, an $r$-plane $P_{\nu}\in \Grass_r(H^1(X,\Q))$ belongs 
to $\Omega^i_r(X)$ if and only if there are only finitely many 
characters $\rho\colon \Z^r \to \C^{\times}$ such that 
$b_j(X, \rho\circ \nu)\ne 0$, for some $j\le i$. 

\subsection{Some applications}
\label{subsec:apps}

In the case when $r=b_1(X)$, Theorem \ref{thm:df cv}
allows us to improve on the discussion from \S\ref{subsec:discuss}
regarding the Betti numbers of maximal free abelian covers.

\begin{theorem}
\label{thm:max abel}
Let $X$ be a CW-complex with finite $k$-skeleton, 
and set $n=b_1(X)$.  Let $X^{\alpha}\to X$ be the maximal 
free abelian cover.   
For each $i\le k$, the following are equivalent:
\begin{enumerate}
\item \label{f1}
$b_j(X^{\alpha})<\infty$, for all $j\le i$. 
\item \label{f2}
$\Omega^i_n(X)\ne \emptyset$.
\item \label{f3}
$\WW^i(X)$ is finite.
\end{enumerate}
\end{theorem}

\begin{proof}
We may assume $n>0$, since, otherwise, there is nothing 
to prove. The only $n$-plane in $H=H^1(X,\Q)$ is $H$ itself;  
thus, $\Grass_n(H)$ consists of a single point, which 
corresponds to the cover $X^{\alpha}\to X$  
given by the projection $\alpha\colon \pi_1(X)\surj \Z^n$. 

By formula \eqref{eq:df b1}, conditions \eqref{f1} and 
\eqref{f2} are equivalent.  On the other hand,  
 $\exp(H\otimes \C)\cap \WW^i(X)=\WW^i(X)$;   
thus, Theorem \ref{thm:df cv} shows that 
conditions \eqref{f2} and \eqref{f3} are equivalent. 
\end{proof}

The following particular case is worth singling out.

\begin{corollary}
\label{cor:b1fab}
Assume $X$ has finite $1$-skeleton. Then 
$\WW^1(X)$ is finite if and only if 
$b_1(X^{\alpha})$ is finite.
\end{corollary}

\subsection{Large $\Omega$-invariants}
\label{subsec:large}
We now analyze the situation when the  
Dwyer--Fried sets comprise the whole Grassmannian.  
The next result follows at once from Theorem \ref{thm:df cv}.

\begin{prop}
\label{prop:df smallcv}
If $\WW^i(X)$ is a finite set, then 
$\Omega^i_r(X)= \Grass_r(H^1(X,\Q))$, for all 
 $r\ge 1$.  
\end{prop}

A general class of examples is provided by nilmanifolds.

\begin{example}
\label{cor:df nilmanifold}

Let $G$ be a torsion-free, finitely 
generated nilpotent group.  Then $G$ admits as 
classifying space a compact nilmanifold of the form 
$M=\R^n/G$. Furthermore, $\VV^i_d(M)$ equals $\{1\}$ if 
$d\le b_i(M)$ and is empty, otherwise.  
(This fact was established in \cite{MP}, using the 
Hochschild-Serre spectral sequence and induction 
on the nilpotency class of $G$.) 
It follows that $\Omega^i_r(M)=\Grass_r(H^1(M,\Q))$, 
for all $i\ge 0$ and $r \ge 1$.

In particular, $\Omega^i_r(T^n)=\Grass_r(\Q^n)$, 
reflecting the fact that every connected cover of the 
$n$-torus is homotopy equivalent to a $k$-torus, for 
some $0\le k\le n$.
\end{example}

Further examples are provided by knot complements. 

\begin{prop}
\label{prop:df prod knots}
Let $K_1, \dots, K_n$ be codimension-$2$ spheres smoothly 
embedded in $S^{m}$, for some $m\ge 3$, and 
let $X=X_1\times \cdots \times X_n$ be the product 
of the respective complements. Then $\Omega^i_r(X)=
\Grass_r(\Q^n)$, for all $i\ge 0$.  
\end{prop}

\begin{proof}
First assume $n=1$, and let $X$ be a knot complement. 
The homology groups of the infinite cyclic cover $X^{\ab}$ 
with coefficients in $\C$ are modules over the principal ideal 
domain $\Lambda=\C[t^{\pm 1}]$.  A  spectral sequence 
argument (due to J.~Levine) shows that $t-1$ acts invertibly 
on these Alexander modules.  It follows that 
$\bigoplus_{j\le i} H_{j}(X^{\ab},\C)$ is a torsion 
$\Lambda$-module, for each $i\ge 0$; hence, its 
support $\WW^i(X)$ is a finite subset of $\C^{\times}$. 

For the general case, we know from the above argument that 
the sets $\WW^i(X_1), \dots$,  $\WW^i(X_n)$ are finite.  
From Corollary \ref{cor:cv prod}, 
we deduce that $\WW^i(X)$ is finite, for each $i\ge 0$.  
The desired conclusion follows from Proposition \ref{prop:df smallcv}.  
\end{proof}

\subsection{Empty $\Omega$-invariants}
\label{subsec:empty}

At the other extreme, we have a rather large supply of spaces 
for which the Dwyer--Fried sets are empty. 

\begin{prop}
\label{prop:df bigcv}
Suppose $\VV^j_1(X) = H^1(X, \C^{\times})$, 
for some $j>0$.  Then $\Omega^i_r(X)=\emptyset$, 
for all $i\ge j$ and all $r\ge 1$.  
\end{prop}

\begin{proof}
If $r>b_1(X)$, then of course $\Omega^i_r(X)=\emptyset$.  
Otherwise,  we may use Theorem \ref{thm:df cv}.  
By hypothesis, $\VV^i(X)=H^1(X, \C^{\times})$, for all 
$i\ge j$.  Hence,  for every rational $r$-plane $P\subset H^1(X,\Q)$, 
the intersection $\exp(P\otimes \C)\cap \VV^i(X)$ is isomorphic 
to $(\C^{\times})^r$, which is an infinite set.  
Thus, $\Omega^i_r(X)=\emptyset$.
\end{proof}

From Propositions \ref{prop:cv wedge} and \ref{prop:df bigcv}, 
we obtain the following corollary. 

\begin{corollary}
\label{cor:df wp} 
Let $X=X_1\vee \cdots \vee X_n$, and  
suppose $b_1(X_s)>0$, for all $s$. 
Then $\Omega^i_r(X) = \emptyset$,  
for all $i, r\ge 1$.   
\end{corollary}

Similarly, using Propositions \ref{prop:cv prod} and \ref{prop:df bigcv}, 
we obtain the following corollary. 

\begin{corollary}
\label{cor:df bigprod}
Let $X=X_1\times \cdots \times X_n$, and  suppose that 
$\VV^1_1(X_s) = H^1(X_s, \C^{\times})$, for all $s$. 
Then $\Omega^i_r(X)=\emptyset$, 
for all $i\ge n$ and $r\ge 1$.  
\end{corollary}

\begin{example}
\label{ex:df surf}
Let $\Sigma_g$ be a Riemann surface of genus $g>1$.  
From Example \ref{ex:cv surf}, we know that 
$\VV^1_1(\Sigma_g)=H^1(\Sigma_g, \C^{\times})$. 
Proposition \ref{prop:df bigcv} now gives
$\Omega^i_r(\Sigma_g)=\emptyset$, for all $i, r\ge 1$. 
For a product of surfaces, Corollary \ref{cor:df bigprod} yields
\begin{equation}
\label{eq:df prod surf}
\Omega^n_r(\Sigma_{g_1}\times \cdots \times \Sigma_{g_n})=\emptyset, 
\quad \text{for all $r\ge 1$}.
\end{equation}
\end{example}

\begin{example}
\label{ex:df bouquet}
Let $Y_m=\bigvee^m S^1$ be a wedge of $m$ circles, 
$m>1$.  From Example \ref{ex:cv wedge}, we know that 
$\VV^1_1(Y_m)=H^1(Y_m, \C^{\times})$. 
By Proposition \ref{prop:df bigcv} (or 
Corollary \ref{cor:df wp}), we have
$\Omega^i_r(Y_m)=\emptyset$, for all $i, r\ge 1$. 
For a product of wedges of circles, Corollary \ref{cor:df bigprod} 
yields
\begin{equation}
\label{eq:df prod wedge}
\Omega^n_r(Y_{m_1}\times \cdots \times Y_{m_n})=\emptyset, 
\quad \text{for all $r\ge 1$}.
\end{equation}
\end{example}

\subsection{A naturality property}
\label{subsec:pullback}
We saw in Lemma \ref{lem:epi cv} that the characteristic 
varieties enjoy a nice naturality property with respect to 
epimorphisms between finitely generated groups.  
As we shall see in Example \ref{ex:epi df}, the analogous 
property does not hold for the Dwyer--Fried sets. Nevertheless, 
it does hold for their complements. 

\begin{prop}
\label{prop:nat df}
Let $\varphi\colon G \surj Q$ be an epimorphism 
from a finitely generated group $G$ to a group $Q$, 
and let $\varphi^*\colon H^1(Q,\Q) \inj H^1(G,\Q)$ 
be the induced monomorphism in cohomology.  
Then, for each $r\ge 1$,  the corresponding morphism 
between Grassmannians, 
$\varphi^*_r\colon \Grass_r(H^1(Q,\Q)) \inj \Grass_r(H^1(G,\Q))$,  
restricts to an embedding 
$\Omega^1_r(Q)^{\compl}\inj \Omega^1_r(G)^{\compl}$. 
\end{prop}

\begin{proof}
Let $P=P_{\nu}$ be a plane in 
$\Grass_r(H^1(Q,\Q)) \setminus \Omega^1_r(Q)$. 
By Theorem \ref{thm:df cv}, the intersection of $\exp(P\otimes \C)$ 
with $\WW^1(Q)$ is infinite. By Lemma \ref{lem:epi cv}, the morphism 
$\hat{\varphi} \colon \widehat{Q} \inj \widehat{G}$ restricts to an 
embedding $\WW^1(Q)\inj \WW^1(G)$.   By Lemma \ref{lem:df funct}, 
we have that $\varphi^*_r(P_{\nu})=P_{\hat\varphi(\nu)}$.   
Hence, the intersection of $\exp(\varphi^*_r(P)\otimes \C)$ 
with $\WW^1(G)$ is infinite.  Thus, $\varphi^*_r(P)$ belongs to 
$\Grass_r(H^1(G,\Q)) \setminus \Omega^1_r(G)$.
\end{proof}

\begin{example}
\label{ex:epi df}
Let $F_n$ be the free group of rank $n>1$. The abelianization 
map, $\ab\colon F_n\surj \Z^n$, induces isomorphisms 
$\ab^*_r\colon \Grass_r(H^1(\Z^n,\Q)) \isom\Grass_r(H^1(F_n,\Q))$, 
for all $r\ge 1$. Corollary \ref{cor:df nilmanifold} shows that 
$\Omega^1_r(\Z^n)=\Grass_r(\Q^n)$, whereas Example \ref{ex:df bouquet} 
shows that $\Omega^1_r(F_n)=\emptyset$. Thus, the map 
$\ab^*_r$ does {\em not}\/ restrict to an embedding 
$\Omega^1_r(\Z^n) \inj \Omega^1_r(F_n)$, though 
of course it restricts to an embedding 
$\Omega^1_r(\Z^n)^{\compl} \inj \Omega^1_r(F_n)^{\compl}$.  
\end{example}

\section{Characteristic subspace arrangements}
\label{sect:df subarr}

In this section, we associate to each space $X$ a 
sequence of subspace arrangements, $\CC_i(X)$, 
all lying in the rational vector space $H^1(X,\Q)$, 
and providing a rough linear approximation to 
the characteristic varieties $\WW^i(X)$.

\subsection{The exponential tangent cone to a variety}
\label{subsec:exp tcone}

We start by reviewing a notion introduced by Dimca, 
Papadima and Suciu in \cite{DPS-duke}, a notion that 
will play an essential role for the rest of this paper. 

\begin{definition}
\label{def:exp tc}
Let $W$ be an algebraic subset 
of the complex algebraic torus $(\C^{\times})^n$.  
The {\em exponential tangent cone}\/ of $W$ at the identity 
$1$ is the homogeneous subvariety 
$\tau_1(W)$ of $\C^n$, given by
\begin{equation}
\label{eq:tau1}
\tau_1(W)= \{ z\in \C^n \mid \exp(\lambda z)\in W,\ 
\text{for all $\lambda\in \C$} \}.
\end{equation}
\end{definition}

The terminology reflects certain similarities between 
the exponential tangent cone $\tau_1(W)$ and the 
classical tangent cone $\TC_1(W)$.  As noted 
in \cite{DPS-duke}, the inclusion $\tau_1(W)
\subseteq \TC_1(W)$ always holds, but not 
necessarily as an equality.  For more on this, 
we refer to \cite{PS-plms} and \cite{Su-aspm}. 
We illustrate the construction with a simple example. 

\begin{example}
\label{ex:tau1 toy}
Let $W=\{t\in (\C^{\times})^2 \mid 
t_1+ t_2 =2\}$. The set $\tau_1(W)$ 
consists of all pairs $(z_1,z_2)\in \C^2$ for 
which $e^{\lambda z_1} + e^{\lambda z_2} = 2$, 
for all $\lambda\in \C$.  Expanding in Taylor series around $0$, 
we find that $z_1^2+z_2^2=0$, and thus $\tau_1(W)=\{0\}$.  
\end{example}

It is readily seen that $\tau_1$ commutes with finite unions and 
arbitrary intersections. Clearly, the exponential tangent cone of 
$W$ only depends on the analytic germ of $W$ at the identity 
$1\in (\C^{\times})^n$.  In particular, $\tau_1(W)\ne \emptyset$ 
if and only if $1\in W$.   The following diagram summarizes 
the situation:
\[
\xymatrixcolsep{5pt}
\xymatrix{
0\ar@{|->}[d]   \ar@{}[r]|-{\in} & \tau_1(W)\: \ar@{^{(}->}[rr] 
\ar^{\exp}[d] &&\C^n\, \ar^{\exp}[d] \\
1 \ar@{}[r]|-{\in} & W \:\ar@{^{(}->}[rr] && (\C^{\times})^n. 
}
\] 

\begin{example}
\label{ex:tau1 torus}
Suppose $T$ is an algebraic $r$-subtorus of $(\C^{\times})^n$.  
Then $T=\exp(P\otimes \C)$, for some $r$-dimensional 
subspace $P$ inside $\Q^n$.  It follows that $\tau_1(T)=P\otimes \C$, 
a rationally defined linear subspace of $\C^n$. Hence, 
$\tau_1(T)$ coincides with $T_1(T)$, the tangent space 
at the identity to the Lie group $T$.
\end{example}

In fact, as shown in \cite{DPS-duke}, a much more general 
result holds:   {\em every}\/ exponential tangent cone 
$\tau_1(W)$ is a union of linear subspaces defined 
over $\Q$.  Since the proof given in \cite{DPS-duke} 
is rather sketchy, we will provide full details in 
Lemma \ref{lem:wronski} and Theorem \ref{thm:tau1} below.  
Along the way, we shall outline an explicit computational 
algorithm for determining $\tau_1(W)$, for any algebraic subvariety 
$W\subset (\C^{\times})^n$.

\subsection{The structure of exponential tangent cones}
\label{subsec:exp tc linear}

To start with, fix a non-zero Laurent polynomial 
$f\in \C[t_1^{\pm 1},\dots , t_n^{\pm 1}]$ with $f(1)=0$.  Write 
\begin{equation}
\label{eq:laurent}
f(t_1,\dots,t_n)= \sum_{a\in S} c_a t_1^{a_1}\cdots t_n^{a_n}, 
\end{equation}
where $S$ is a finite subset of $\Z^n$, and $c_a\ne 0$ 
for each $a=(a_1,\dots,a_n)\in S$. We say a partition 
 $\cS=(\cS_1\mid \cdots \mid \cS_q)$ of the support $S$ is 
 {\em admissible}\/ if $\sum_{a\in \cS_i} c_a =0$, for all 
 $1\le i\le q$.  In particular, the trivial partition $\cS=(S)$ 
 is admissible, by our assumption on $f$.
 
To a partition $\cS$ as above we associate a linear 
subspace $L(\cS)\subset \Q^n$, defined as
\begin{equation}
\label{eq:rat sub}
L(\cS)= \{ x\in \Q^n \mid (a-b)\cdot x =0, 
\; \forall a,b\in \cS_i, \; \forall\, 1\le i\le q \}.
\end{equation}

Given a vector $z\in \C^n$, 
define an analytic function $\phi_z\colon \C\to \C$ by 
\begin{equation}
\label{eq:phiz}
\phi_z(\lambda)=f(e^{\lambda z_1},\dots, e^{\lambda z_n})=
\sum_{a\in S} c_a e^{(a\cdot z)\lambda},
\end{equation}
where $a\cdot z=\sum_{i=1}^n a_iz_i$ is the standard dot product. 

\begin{lemma}
\label{lem:wronski}
The function $\phi_z$ vanishes identically if and only if 
$z$ belongs to a linear subspace of the form 
$L(\cS)\otimes \C$, where  $\cS$ is a admissible 
partition of $S$.
\end{lemma}

\begin{proof}
Suppose $z\in L(\cS)\otimes \C$.  Then, for each part $\cS_i$, 
there is a constant $k_i$ such that $a\cdot z=k_i$ for all $a\in \cS_i$. 
Write $f_i=\sum_{a\in \cS_i} c_a t_1^{a_1}\cdots t_n^{a_n}$, 
and let  $\phi^{i}_z$ be the corresponding analytic function.   
We then have 
\begin{equation*}
\label{eq:phi i}
\phi^{i}_z(\lambda)
=\sum_{a\in \cS_i}c_a e^{(a\cdot z) \lambda}=
\sum_{a\in \cS_i}c_ae^{k_i \lambda}=
\Big(\sum_{a\in \cS_i}c_a\Big)e^{k_i\lambda},
\end{equation*}
which equals $0$, since $\cS$ is admissible.  
On the other hand,  $\phi_z=\sum_{i=1}^q \phi^i_z$.  
Thus, $\phi_z\equiv 0$.

Conversely, suppose $\phi_z\equiv 0$.  
Since $f$ is non-zero, we can find a maximal subset 
$\cS_1\subset S$ such that $a\cdot z$ is constant 
for all $a\in\cS_1$.  Replacing $S$ by $S\setminus \cS_1$ 
and proceeding in this fashion, we ultimately arrive at 
a partition $\cS=(\cS_1\mid \cdots \mid \cS_q)$, 
and pairwise distinct constants 
$k_1,\dots, k_q$ such that 
$a\cdot z=k_i$, for all $a \in \cS_i$.   
Thus, $z\in L(\cS)\otimes \C$.  On the 
other hand, the identity  
$\phi_z=\sum_{i=1}^q \phi^i_z\equiv 0$, 
translates into a linear dependence,   
\[
\sum_{i=1}^{q} \Big(\sum_{a\in \cS_i}c_a\Big)e^{k_i\lambda} =0, 
\quad \forall \lambda\in \C.
\]
Computing the Wronskian of 
$e^{k_1 \lambda}, \dots ,e^{k_q \lambda}$ 
reveals that these functions are linearly independent.  
Hence, $\sum_{a\in \cS_i}c_a=0$, for all $i$, 
showing that $\cS$ is an admissible partition. 
\end{proof}

\begin{theorem}[\cite{DPS-duke}]
\label{thm:tau1}
The exponential tangent cone $\tau_1(W)$ is a finite union 
of rationally defined linear subspaces of $\C^n$.
\end{theorem}

\begin{proof}
The variety $W\subset (\C^{\times})^n$ is the common 
zero-locus of finitely many Laurent polynomials in 
$n$ variables, say, $f_1, \dots, f_m$.  
Since $\tau_1$ commutes with intersections, we 
have that $\tau_1(W)= \bigcap_{j=1}^m \tau_1(V(f_j))$.  
Thus, it is enough to consider the case $W=V(f)$, 
where $f$ is a non-zero Laurent polynomial.  Without loss 
of generality, we may assume $f(1)=0$, for otherwise 
$\tau_1(W)=\emptyset$. 

By definition, a point $z\in \C^n$ belongs to $\tau_1(W)$ 
if and only if the function $\phi_z$ defined in \eqref{eq:phiz} 
vanishes identically.  By Lemma \ref{lem:wronski}, this 
happens precisely when $z$ belong to a linear subspace 
of the form $L(\cS)\otimes \C$, with $\cS$ an 
admissible partition.  Hence, $\tau_1(W)$ is the 
union of all such (rationally defined) subspaces.  
\end{proof}

Note that only maximal partitions contribute to the above union.  
Indeed, if $\cS$ is a refinement of $\cS'$, then 
$L(\cS)\supseteq L(\cS')$. 

\subsection{A class of rational subspace arrangements}
\label{subsec:df subarr}

As usual, let $X$ be a connected CW-com\-plex with finite 
$k$-skeleton. Set $n=b_1(X)$, and identify $H^1(X, \C)=\C^n$ 
and  $H^1(X, \C^{\times})^{\circ}=(\C^{\times})^n$.   
Let us apply the exponential tangent cone construction to the 
characteristic varieties $\WW^i(X)\subseteq (\C^{\times})^n$. 
In view of Theorem \ref{thm:tau1}, we may make  
the following definition.  

\begin{definition}
\label{def:df subarr}
For each $i\le k$, the {\em $i$-th characteristic arrangement}\/ 
of $X$, denoted $\CC_i(X)$, is the subspace arrangement in 
$H^1(X,\Q)$ whose complexified union is the exponential 
tangent cone to $\WW^i(X)$:
\begin{equation}
\label{eq:dfxi}
\tau_1(\WW^i(X)) =  
\bigcup_{L\in \CC_i(X)} L\otimes \C.
\end{equation}
\end{definition}

Put differently, the set of rational points on the exponential 
tangent cone to $\WW^i(X)$, 
\begin{equation}
\label{eq:tau1q}
\tau_1^{\Q}(\WW^i(X)) = \tau_1(\WW^i(X))\cap 
H^1(X,\Q),
\end{equation}
equals the union of the subspaces comprising the 
$i$-th characteristic arrangement.    
We thus have a sequence $\CC_0(X)$, $\CC_1(X), \dots , 
\CC_k(X)$ of rational subspace arrangements, all lying in 
the same affine space $H^1(X,\Q)=\Q^n$.  

From Lemma \ref{lem:cv inv}, we easily see that 
the subspace arrangements $\CC_i(X)$ 
depend only on the homotopy type of $X$.  
Using now Lemma \ref{lem:epi cv}, 
we obtain the following consequence. 

\begin{prop}
\label{prop:nat tau}
Let $\varphi\colon G \surj Q$ be an epimorphism 
from a finitely generated group $G$ to a group $Q$.  
The induced morphism in cohomology, 
$\varphi^*\colon H^1(Q,\Q) \inj H^1(G,\Q)$, 
restricts to a map 
$\varphi^*\colon \tau^{\Q}_1(Q)\inj \tau^{\Q}_1(G)$, whose 
image coincides with $\tau_1^{\Q}(\hat{\varphi}(\WW^1(Q)))$.
\end{prop}

\begin{example}
\label{ex:chain}
Let $L$ be the closed three-link chain 
(the link $6^3_1$ from Rolfsen's tables \cite{Ro}), 
with Alexander polynomial
$\Delta_L=t_1+t_2+t_3-t_1t_2-t_1t_3-t_2t_3$. 
Let $X=S^3\setminus L$, and 
fix a meridional basis $e_1,e_2,e_3$ for $H_1(X,\Z)=\Z^3$. 
From Example \ref{ex:cv link}, we know that 
$\WW^1(X)$ is the zero-locus of $\Delta_L$.  A straightforward   
computation show that the support of this polynomial, 
$S=\{e_1,e_2,e_3,e_1+e_2,e_1+e_3,e_2+e_3\}$, 
has precisely three (maximal) admissible partitions, namely, 
$\cS=(e_1, \, e_2+e_3 \mid e_2, \, e_1+e_2 \mid e_3, \, e_1+e_3 )$,
$\cS'=(e_1, \, e_1+e_2 \mid e_2, \, e_1+e_3 \mid e_3, \, e_2+e_3 )$,
and $\cS''=(e_1, \, e_1+e_3 \mid e_2, \, e_2+e_3 \mid e_3, \, e_1+e_2 )$. 
Therefore, the arrangement $\CC_1(X)$ consists of 
three lines in $\Q^3$, to wit, 
$L(\cS)= \{ x_1=x_2+x_3=0\}$, 
$L(\cS')= \{ x_2=x_1+x_3=0\}$, and 
$L(\cS'')= \{ x_3=x_1+x_2=0\}$.
\end{example}

\section{An upper bound for the $\Omega$-invariants}
\label{sect:bound}

In this section, we give a readily computable 
approximation to the Dwyer--Fried sets, based 
on the exponential tangent cone construction 
described above.

\subsection{$\Omega$-invariants and characteristic arrangements}
\label{subsec:bound}

As usual, let $X$ be a connected CW-complex 
with finite $k$-skeleton, for some $k\ge 1$.  
For each $i\le k$, denote by $\CC_i(X)$ the 
rational subspace arrangement defined by \eqref{eq:dfxi}.

\begin{theorem}
\label{thm:tau bound}
For all $i\le k$ and all $r\ge 1$, 
we have: 
\begin{equation}
\label{eq:ubound}
\Omega^i_r(X) \subseteq \bigg(
\bigcup_{L\in \CC_i(X)} \big\{P \in \Grass_r(H^1(X,\Q)) \bigmid 
P\cap L \ne \{0\} \big\} \bigg)^{\compl}.
\end{equation}
\end{theorem}

\begin{proof}
Fix an $r$-plane $P$ inside $H^1(X,\Q)$, 
and let $T=\exp(P\otimes \C)$ be the corresponding 
algebraic subtorus in $H^1(X,\C^{\times})$.   Then:
\[
\begin{aligned}
P\in \Omega^i_r(X) & \underset{\text{(i)}}{\same}
T \cap \WW^i(X) \text{ is finite}\\
& \underset{\text{(ii)}}{\implies}  
\tau_1(T \cap \WW^i(X)) =\{0\}\\
&\underset{\text{(iii)}}{\same}
(P\otimes \C) \cap \tau_1(\WW^i(X))= \{0\}\\
& \underset{\text{(iv)}}{\same}
P\cap L=\{0\},  \text{ for each $L\in \CC_i(X)$},
\end{aligned}
\]
where in (i) we used Theorem \ref{thm:df cv},  
in (ii) we used Definition \ref{def:exp tc}, 
in (iii) we used the fact that $\tau_1$ commutes 
with intersections and $\tau_1(T)=P\otimes \C$, 
and in (iv) we used Definition \ref{def:df subarr}. 
\end{proof}

If $\tau_1(\WW^i(X))= \{0\}$, then the right-hand side 
of \eqref{eq:ubound} is the whole Grassmannian, and 
the upper bound is tautological.  Under appropriate 
hypothesis, though, the bound is sharp. 
The next corollary isolates a class of spaces for 
which this happens.

\begin{corollary}
\label{cor:tau sharp}
Suppose all positive-dimensional components of $\WW^i(X)$ 
are algebraic subtori.  Then, for all $r\ge 1$, 
we have: 
\begin{equation}
\label{eq:bound sharp}
\Omega^i_r(X) = \Big(
\bigcup_{L\in \CC_i(X)} \big\{P \in \Grass_r(H^1(X,\Q)) \bigmid 
P\cap L \ne \{0\} \big\} \Big)^{\compl}.
\end{equation}
\end{corollary}

\begin{proof}
By hypothesis, $\WW^{i}(X)=W \cup Z$, where 
$W=\exp(\tau_1(\WW^i(X)))$ is a finite union of algebraic 
subtori, and $Z$ is a finite  set.  
With notation as in the proof of Theorem \ref{thm:tau bound}, 
suppose $\tau_1(T \cap \WW^i(X)) =\{0\}$. Then $T\cap W=\{1\}$, 
and thus $T\cap \WW^{i}(X)$ is finite. This shows implication 
(ii) can be reversed in this situation, thereby finishing the proof.
\end{proof}

A class of spaces for which the hypothesis of 
Corollary \ref{cor:tau sharp} holds are the ``straight" 
spaces studied in \cite{Su-aspm}.  In general, though, 
the inclusion from Theorem \ref{thm:tau bound} is strict, 
as illustrated in the next example.

\begin{example}
\label{ex:link 421}
Let $L$ be the $2$-component link 
denoted $4^2_1$ in Rolfsen's tables \cite{Ro}, 
and let $X$ be its complement. Then $\Delta_L=t_1+t_2$,   
and so the variety $\WW^1(X)\subset (\C^{\times})^2$ 
consists of the identity $1$, together with the $1$-dimensional 
translated torus $\{t_1t_2^{-1}=-1\}$.  Using Theorem \ref{thm:max abel}, 
we find that $\Omega^1_2(X)=\emptyset$. On the other hand, 
$\tau_1(\WW^1(X))=\{0\}$; thus, the right-hand side 
of \eqref{eq:ubound} equals $\Grass_2(\Q^2)=\pt$. 
\end{example}

\subsection{The case of infinite cyclic covers}
\label{subsec:zcov}

When $r=1$, the upper bound from 
Theorem \ref{thm:tau bound} becomes an equality. 
To see why that is the case, 
set $H=H^1(X,\Q)$, and identify the Grassmannian 
$\Grass_1 (H)$ with the projective space 
$\QP^{n-1}=\bP(H)$, where $n=b_1(X)$. 
Given a linear subspace 
$L\subset H$, denote its projectivization by 
$\bP(L)$.  We then have the 
following theorem, which recovers results from 
Dwyer--Fried \cite{DF} and Papadima--Suciu \cite{PS-plms}.  
For completeness, we include a proof adapted to our setting.

\begin{theorem}
\label{thm:df1}
Let $X$ be a CW-complex with finite $k$-skeleton. 
Set $n=b_1(X)$.  Then $\Omega^i_1(X) = \QP^{n-1} \setminus 
\bigcup_{L\in \CC_i(X)} \bP(L)$, 
for all $i\le k$, 
\end{theorem}

\begin{proof}
Let $P$ be a line in $H=\Q^n$.  Set $T=\exp(P\otimes \C)$ 
and $T'=T\cap \WW^i(X)$. Looking at the proof of 
Theorem \ref{thm:tau bound}, we only need to reverse 
implication (ii); that is, we need to show that 
$\tau_1(T')= \{0\} \implies T' \text{ is finite}$.  
Now, $T'$ is a Zariski closed subset of $T\cong \C^{\times}$. 
So, if $T'$ were not finite, then $T'$ would equal $T$, implying 
$\tau_1(T')\cong \C$. 
\end{proof}

This theorem yields a nice qualitative result about the 
rank-$1$ Dwyer--Fried sets. 

\begin{corollary}[\cite{DF}]
\label{cor:omega1}
Each set $\Omega^i_1(X)$ is the complement of a 
finite union of projective subspaces in $\QP^{n-1}$.  
In particular, $\Omega^i_1(X)$ is a Zariski open set in 
$\QP^{n-1}$.
\end{corollary}

Theorem \ref{thm:df1} also gives an easy-to-use 
homological finiteness test for the regular, infinite cyclic 
covers of a space $X$.  

\begin{corollary}[\cite{PS-plms}]
\label{cor:nutau}
Let $X^{\nu}\to X$ be a connected, regular $\Z$-cover, classified 
by an epimorphism $\nu\colon \pi_1(X)\surj \Z$, and  let 
$\bar\nu\in H^1(X,\Z)\subset H^1(X,\C)$ be the 
corresponding cohomology class.  Then 
$\sum_{i=1}^{k} b_{i} (X^{\nu}) < \infty$ if and only if 
$\bar\nu\not\in  \tau_1( \WW^k(X))$.
\end{corollary}

\begin{example}
\label{ex:chain again}
Let $X$ be the link complement from Example \ref{ex:chain}. 
By Theorem \ref{thm:df1}, the set $\Omega^1_1(X)$ 
is the complement in $\QP^2$ of the points $(0,1,-1)$, $(1,0,-1)$, 
and $(1,-1,0)$.  Thus, every regular $\Z$-cover $X^{\nu}\to X$ 
has $b_1(X^{\nu})<\infty$, except for the three covers 
specified by the those vectors. 
\end{example}

\section{Special Schubert varieties}
\label{sect:schubert}

In this section, we reinterpret the upper bound for the 
Dwyer--Fried sets in terms of a well-known geometric 
construction, by showing that each $\Omega$-invariant  
lies in the complement of an arrangement of (special) 
Schubert varieties.  

\subsection{The incidence correspondence}
\label{subsec:incidence}

We start by recalling a classical construction from 
algebraic geometry (see \cite[p.~69]{Har}).  
Fix a field $\k$, and let $V$ be a homogeneous variety in 
$\k^n$.  Consider  the locus of $r$-planes in $\k^n$ intersecting 
$V$ non-trivially, 
\begin{equation}
\label{eq:incident}
\sigma_r(V) = \big\{ P \in \Grass_r(\k^n) 
\bigmid P \cap  V \ne \{0\} \big\}.
\end{equation}
This set is a Zariski closed subset of the Grassmannian 
$\Grass_r(\k^n)$, called the {\em variety of incident 
$r$-planes}\/ to $V$. In particular, $\sigma_r(\{0\})=\emptyset$. 

Recall $\dim  \Grass_r(\k^n)  =r(n-r)$. 
The following well-known fact 
(cf.~\cite[p.~153]{Har}) will be useful to us.

\begin{lemma}
\label{lem:dim incidence}
Let $V$ be a homogeneous, irreducible variety in $\k^n$, 
of dimension $m>0$.  Then, for all $0<r< n-m$, the incidence 
variety $\sigma_r(V)$ is an irreducible 
subvariety of $\Grass_r(\k^n)$, of dimension $(r-1)(n-r)+m-1$. 
\end{lemma}

Particularly simple is the case when $V$ is 
a linear subspace $L\subset \k^n$.  The corresponding 
incidence variety, $\sigma_r(L)$, is known as the  
{\em special Schubert variety}\/ defined by $L$. 
Clearly, $\sigma_1(L)=\bP(L)$, viewed as a projective 
subspace in $\bP(\k^n)$.

\begin{corollary}
\label{cor:dim sigma}
Let $L$ be a non-zero, codimension $d$ linear subspace in $\k^n$.  
Then $\sigma_r(L)$ has codimension $d-r+1$ in $\Grass_r(\k^n)$. 
\end{corollary}

To write down equations for the special Schubert varieties, 
start by embedding $\Grass_r(\k^n)$ in the projective 
space $\bP(\bigwedge^r \k^{n})$ via the Pl\"{u}cker 
embedding, which sends an $r$-plane $P$ to $\bigwedge^r P$, 
and let $p_{i_1,\dots ,i_r}(P)$, $1\le i_1< \cdots < i_r\le n$, be 
the Pl\"{u}cker coordinates of this plane. 
Let $L$ be an $s$-dimensional plane in $\k^n$, represented 
as the row space of an $s\times n$ matrix. Then $L$ 
meets an $r$-plane $P$ non-trivially if and only if 
all the maximal minors of the matrix 
$\left( \begin{smallmatrix} L \\ P \end{smallmatrix} \right)$ 
vanish.   Laplace expansion of each minor along the rows 
of $P$ yields a linear equation in the Pl\"{u}cker coordinates. 
We illustrate the procedure with an example.

\begin{example}
\label{ex:grass24}
The Grassmannian $\Grass_2(\k^4)$ is the hypersurface 
in $\k\bP^5$ with equation
\begin{equation}
\label{eq:plucker}
p_{12}p_{34}-p_{13}p_{24}+p_{23}p_{14}=0.
\end{equation}
Let $L$ be a plane in $\k^4$.  The corresponding 
Schubert variety, $\sigma_2(L)$, is the $3$-fold in 
$\Grass_2(\k^4)$ cut out by the hyperplane 
\begin{equation}
\label{eq:plucker plane}
L_{34}p_{12}-L_{24}p_{13}+L_{14}p_{23}+L_{23}p_{14}
-L_{13}p_{24}+L_{12}p_{34}=0,
\end{equation}
where $L_{ij}$ is the $(ij)$-minor of a $2\times 4$ matrix 
representing $L$. 
\end{example}

\subsection{Reinterpreting the upper bound for the $\Omega$-sets}
\label{subsec:df schubert}

Using the notions discussed above, we may reformulate 
the exponential tangent cone upper bound from 
Theorem \ref{thm:tau bound}, as follows.

\begin{theorem}
\label{thm:sch bound}
Let $X$ be a CW-complex with finite $k$-skeleton. 
Then
\begin{equation}
\label{eq:omega schubert}
\Omega^i_r(X) \subseteq  \Grass_r(H^1(X,\Q)) \setminus 
\bigcup_{L\in \CC_i(X)} \sigma_r(L ),
\end{equation}
for all $i\le k$ and $r\ge 1$. 
\end{theorem}

In other words, each Dwyer--Fried set $\Omega^i_r(X)$ is 
contained in the complement of a Zariski closed subset 
of $\Grass_r(H^1(X,\Q))$, namely, the union of special Schubert 
varieties corresponding to the subspaces belonging to the 
characteristic arrangement $\CC_i(X)$. 

\begin{remark}
\label{rem:not strict}
If inclusion \eqref{eq:omega schubert} holds as equality, then 
clearly $\Omega^i_r(X)$ is a Zariski open subset of $\Grass_r(H^1(X,\Q))$. 
Nevertheless, as we saw in Example \ref{ex:link 421}, this 
inclusion can be strict even when the set $\Omega^i_r(X)$ is open. 
\end{remark}

\begin{remark}
\label{rem:not open}
As shown in Theorem \ref{thm:df1}, inclusion 
\eqref{eq:omega schubert} holds as equality if $r=1$. 
On the other hand, if $r>1$, the sets $\Omega^i_r(X)$ 
are not necessarily open.  
This rather surprising fact was first noticed by Dwyer and Fried, 
who constructed in \cite{DF} a $3$-dimensional cell complex 
of the form $X=(T^3\vee S^2)\cup_{\varphi} e^3$ 
for which $\Omega^2_2(X)$ is a finite set 
(see also \cite{Su11} for more details).  
We will come back to this topic in \S\ref{sect:tt}, 
where we will give an example of a finitely presented 
group $G$ for which $\Omega^1_2(G)$ is a single point. 
\end{remark}

Corollary \ref{cor:tau sharp} can also be reformulated, 
as follows. 

\begin{corollary}
\label{cor:tau schubert}
Suppose all positive-dimensional components of $\WW^i(X)$ 
are algebraic subtori.  Then, for all $r\ge 1$, 
we have: 
\begin{equation}
\label{eq:sharp schubert}
\Omega^i_r(X) = \Grass_r(H^1(X,\Q)) \setminus 
\sigma_r( \tau_1^{\Q}(\WW^i(X))).
\end{equation}
In particular, $\Omega^i_r(X)$ is a Zariski open subset 
in $\Grass_r(H^1(X,\Q))$.
\end{corollary}

\begin{example}
\label{ex:schubert g}
Let $G=F_2\times F_2$.   By Corollary \ref{cor:cv prod}, 
the first characteristic variety $\VV^1(G)$ consists of two, 
$2$-dimensional algebraic subtori in  $(\C^{\times})^4$,  
namely $T_1=\{ t_1=t_2=1\}$ and $T_2=\{t_3=t_4=1\}$.  
Thus, the characteristic subspace arrangement $\CC_1(G)$ 
consists of two planes in $\Q^4$, namely, 
$L_1=\{ x_1=x_2=0\}$ and $L_2=\{x_3=x_4=0\}$. 

Now identify the Grassmannian $\Grass_2(\Q^4)$ with the quadric  
in $\QP^5$ given by equation \eqref{eq:plucker}.  In view 
of formulas \eqref{eq:sharp schubert} and \eqref{eq:plucker plane}, 
then, the set $\Omega^1_2(X)$ is the complement 
in $\Grass_2(\Q^4)$ of the variety cut out by the 
hyperplanes $\{p_{12}=0\}$ and $\{p_{34}=0\}$.
\end{example}

In favorable situations, the information provided by 
Theorem \ref{thm:sch bound} allows us to identify 
precisely the $\Omega$-sets.

\begin{corollary}
\label{cor:df sch} 
Suppose the characteristic arrangement $\CC_i(X)$ 
contains a non-zero subspace of codimension $d$.  
Then $\Omega^i_r(X) = \emptyset$, for all $r> d$.  
\end{corollary}

\begin{proof}
Let $L$ be a subspace in $\CC_i(X)$, with $L\ne \{0\}$ and $\codim (L)=d$.  
Suppose $r\ge d+1$.  By Corollary \ref{cor:dim sigma}, we 
have $\codim (\sigma_r(L))=d-r+1\le 0$, 
and so $\sigma_r(L)=\Grass_r(H^1(X,\Q))$. In view of 
\eqref{eq:omega schubert}, this forces 
$\Omega^i_r(X)$ to be empty.
\end{proof}

\begin{corollary}
\label{cor:vanish omega} 
Let $X^{\alpha}$ be the maximal free abelian 
cover of $X$. If $\tau_1(\WW^1(X))\ne \{0\}$,   
then $b_1(X^{\alpha})=\infty$.  
\end{corollary}

\begin{proof}
Set $n=b_1(X)$.  By hypothesis, the arrangement 
$\CC_1(X)$ contains a non-zero subspace, say $L$. 
Since $\codim(L)\le n-1$, the previous 
corollary implies $\Omega^1_n(X)=\emptyset$. 
The desired conclusion follows from \eqref{eq:df b1}.
\end{proof}

The converse to the implication from Corollary \ref{cor:vanish omega} 
does not hold.   For instance, if $X$ is the complement to the $2$-component 
link from Example \ref{ex:link 421}, then $\tau_1(\WW^1(X))=\{0\}$, yet 
$\Omega^1_2(X)=\emptyset$, and so $b_1(X^{\alpha})=\infty$. 

\section{Translated tori}
\label{sect:tt}

In this section, we give a fairly complete description of the 
Dwyer--Fried sets $\Omega^i_r(X)$ in the case when all 
the positive-dimensional components of the characteristic 
variety $\WW^i(X)$ are torsion-translated subtori of the 
character torus of $\pi_1(X)$. 

\subsection{Intersections of translated tori}
\label{subsec:int tt}
We start by recalling some basic terminology.  
Fix a complex algebraic torus $(\C^{\times})^n$.  
A {\em subtorus}\/ is a connected algebraic subgroup; 
every such subgroup  
is of the form $T=\exp(L\otimes \C)$, for some linear 
subspace $L\subset \Q^n$.  A {\em translated}\/ subtorus 
is a coset of a subtorus, i.e, a subvariety of $(\C^{\times})^n$ 
of the form $\rho\cdot T$;  
if $\rho \in (\C^{\times})^n$ can be chosen to have finite order, we say $\rho\cdot T$ 
is a {\em torsion-translated}\/ subtorus. 

For our purposes, it is important to understand the way 
translated subtori intersect inside an algebraic torus.  
Much of the relevant theory was worked out by 
E.~Hironaka in \cite{Hi}; the theory was further 
developed by Suciu, Yang, and Zhao in \cite{SYZ2}. 
The specific result we shall need can be formulated as follows.  

\begin{prop}[\cite{SYZ2}]
\label{prop:tt intersect}
Let $T_1=\exp(L_1\otimes \C)$ and $T_2=\exp(L_2\otimes \C)$ 
be two algebraic subtori in $(\C^{\times})^n$, 
and let $\rho_1$ and $\rho_2$ be two elements in $(\C^{\times})^n$.  
Then 
\begin{enumerate}
\item \label{ss1}
The variety $Q=\rho_1 T_1 \cap \rho_2 T_2$ is non-empty if and 
only if $\rho^{}_1 \rho_2^{-1}$ belongs to the subgroup 
$T_1 \cdot T_2=\exp((L_1+L_2)\otimes \C)$. 

\item \label{ss2}
If the above condition is satisfied, then $\dim Q=\dim(T_1\cap T_2)$.
\end{enumerate}
\end{prop} 
  
Given a linear subspace $L\subset \Q^n$ and an element 
$\rho\in (\C^{\times})^n$, it is convenient to introduce the 
following notation:
\begin{equation}
\label{eq:sll}
\sigma_r(L,\rho) = \{P \in \Grass_r(\Q^n) \mid 
\text{$\rho\in \exp((P+L)\otimes \C)$ and $P\cap L\ne \{0\}$} \}.
\end{equation}

Clearly, $\sigma_r(L,\rho)\subseteq \sigma_r(L)$ for all $\rho$, 
with equality if $\rho\in \exp(L\otimes \C)$. If $\dim L=1$, then 
$\sigma_r(L,\rho)$ consists of those $r$-planes $P$ 
which contain the line $L$ and for which  
$\rho\in \exp(P\otimes \C)$; 
in particular, if $\rho\notin \exp(L\otimes \C)$, 
then $\sigma_1(L,\rho)=\emptyset$.

\begin{remark}
\label{rem:sigmalrho}
Unlike the special Schubert varieties $\sigma_r(L)$, 
the sets $\sigma_r(L,\rho)$ need not be Zariski closed 
subsets of the rational Grassmannian.  For instance, 
if $L$ has dimension $2\le d \le n-1$, and $\rho$ is 
an element of finite-order in $(\C^{\times})^n\setminus  
\exp(L\otimes \C)$, then, as we shall see 
in the proof of Theorem \ref{thm:df tors translated}, the 
set $\sigma_{r}(L,\rho)$ is not a closed subset of 
$\Grass_{r}(\Q^n)$, for any $2\le r\le d$.  In fact, 
as we shall see in the proof of Proposition 
\ref{prop:ttdf}, if $d= n-1$, then 
$\sigma_{r}(L,\rho) = \Grass_{r}(\Q^n) \setminus \Grass_{r}(L)$.
\end{remark}

Proposition \ref{prop:tt intersect} has an immediate corollary.

\begin{corollary}
\label{cor:exp intersect}
Let $\rho T$ be a translated subtorus in  $(\C^{\times})^n$.   
Write $T=\exp(L\otimes \C)$, 
for some linear subspace $L\subset \Q^n$.  
Given a subspace $P\in \Grass_r(\Q^n)$, we have 
\begin{equation}
\label{eq:pinlrho}
\dim_{\C} ( \exp(P\otimes \C) \cap \rho T )>0 \same 
P\in \sigma_r(L,\rho).
\end{equation}
\end{corollary} 

\subsection{Translated tori and Dwyer--Fried sets}
\label{subsec:df tt}
We are now ready to state and prove the main result of this section.

\begin{theorem}
\label{thm:df translated}
Let $X$ be a CW-complex with finite $k$-skeleton, 
for some $k\ge 1$.   For a fixed integer $i\le k$, suppose 
all positive-dimensional components of 
$\WW^i(X)$ are translated subtori in 
$H^1(X,\C^{\times})^{\circ}$.   Write 
\begin{equation}
\label{eq:wix}
\WW^i(X) = \Big( \bigcup_{\alpha}  \rho_{\alpha} T_{\alpha} \Big)
\cup Z,
\end{equation}
where $Z$ is a finite set and 
$T_{\alpha}=\exp(L_{\alpha}\otimes \C)$, for some 
linear subspace $L_{\alpha}\subset H^1(X,\Q)$ and 
some $\rho_{\alpha} \in H^1(X,\C^{\times})^{\circ}$. 
Then, for all $r\ge 1$,
\begin{equation}
\label{eq:df translated}
\Omega^i_r(X) =
\Grass_r(H^1(X,\Q)) \setminus 
\bigcup_{\alpha}  \sigma_r(L_\alpha, \rho_\alpha).
\end{equation}
\end{theorem}

\begin{proof}
Fix an $r$-plane $P$ inside $H^1(X,\Q)$, 
and let $T=\exp(P\otimes \C)$ be the corresponding 
algebraic subtorus in $H^1(X,\C^{\times})^{\circ}$.   Then:
\[
\begin{aligned}
P\in \Omega^i_r(X) 
& {\ \same \ }  
T \cap \WW^i(X) \text{ is finite,} 
&& \text{by Theorem \ref{thm:df cv}} \\
&{\ \same\ } 
 \text{$T\cap \rho_{\alpha} T_{\alpha}$ is finite, for all $\alpha$,} 
&& \text{by assumption \eqref{eq:wix}} \\ 
&{\ \same \ } 
 \text{$P\notin \sigma_r(L_{\alpha}, \rho_{\alpha})$, for all $\alpha$,} 
&& \text{by Corollary \ref{cor:exp intersect}.}
\end{aligned}
\]
This ends the proof.
\end{proof}

If all the translation factors $\rho_{\alpha}$ in \eqref{eq:wix} 
are equal to $1$, then we are in the situation of 
Corollary \ref{cor:tau schubert}, and the Dwyer--Fried sets 
$\Omega^i_r(X)$ are Zariski open subsets 
of $\Grass_r(H^1(X,\Q))$.  On the other hand, if 
one of those translation factors is a non-trivial root 
of unity, the corresponding subtorus has dimension 
$r>1$, and that subtorus intersects the other subtori transversely 
at the identity, then, as we show next, the set $\Omega^i_r(X)$ 
is not open, even for the usual topology on the Grassmannian.

\begin{theorem}
\label{thm:df tors translated}
As above, suppose 
all positive-dimensional components of $\WW^i(X)$ 
are translated subtori in $H^1(X,\C^{\times})^{\circ}$, 
of the form $\rho_{\alpha} T_{\alpha}$.  
Furthermore, suppose there is a subtorus $T$ 
of dimension $d\ge 2$ such that:
\begin{enumerate}
\item \label{int1}
There is at least one component $\rho_{\beta} T_{\beta}$ 
with $T_\beta=T$. 

\item \label{int1.5}
If $T_{\alpha}=T$, then $\rho_{\alpha}$ has finite order 
(modulo $T$) and $\rho_{\alpha}\notin T$. 

\item \label{int2} 
If $T_\alpha\ne T$, then 
$\tau_1(T_{\alpha})\cap \tau_1(T)=\{0\}$.
\end{enumerate}
Then, for each $2\le r\le d$, the set $\Omega^i_r(X)$ is {\em not}\/ 
open in $\Grass_r(H^1(X,\Q))$.
\end{theorem}

\begin{proof}
Set $n=b_1(X)$, and identify $H^1(X,\Q)=\Q^n$.  
Write $T_{\alpha}=\exp(L_{\alpha}\otimes \C)$, 
with $L_{\alpha}=\tau_1(T_{\alpha})$ a linear subspace 
of $\Q^n$, and similarly, $T=\exp(L\otimes \C)$.  
For each component $\rho_\alpha T_\alpha$ with 
$L_\alpha= L$, write $\rho_{\alpha}=\exp(2 \pi i \lambda_{\alpha})$; 
by assumption \eqref{int1.5}, then, we have $\lambda_{\alpha}\in \Q^n$ 
but $\lambda_{\alpha}\notin L_\alpha+\Z^n$.  

Fix a basis $\{v_1,\dots ,v_d\}$ for $L$  and an integer 
$2\le r \le d$.  By assumption \eqref{int1}, there is an 
index $\beta$ for which $L_\beta=L$.  Consider the 
$r$-dimensional subspaces 
\[
P=\spn\{v_1,\dots,  v_r \}
\quad\text{and}\quad  
P_q=\spn\{v_1,\dots, v_{r-1}, v_r+\tfrac{\lambda_\beta}{q} \}.
\]  
Clearly, $P_q \to P$.  By Theorem \ref{thm:df translated}, we 
have that $\Omega^i_r(X) =\Grass_r(\Q^n) \setminus 
\bigcup_{\alpha}  \sigma_r(L_\alpha, \rho_\alpha)$.
Thus, to finish the proof, it is enough to show 
that $P\notin \sigma_r(L_\alpha,\rho_\alpha)$, for all $\alpha$, 
yet $P_q\in \sigma_r(L_\beta,\rho_\beta)$, for all $q$. 

To prove the first claim, first note that $P\subseteq L$. 
Next, consider an index $\alpha$ so that $L_\alpha \ne L$.  
By assumption \eqref{int2}, we have that 
$L_\alpha \cap L=\{0\}$;  thus, $P\notin \sigma_r(L_\alpha)$, 
and so $P\notin \sigma_r(L_\alpha,\rho_\alpha)$.  
Finally, consider an index $\alpha$ so that $L_\alpha = L$. 
Then $\rho_{\alpha}\notin \exp((P+L_\alpha)\otimes \C)=T_\alpha$;  
hence, $P\notin \sigma_r(L_{\alpha},\rho_{\alpha})$. 

To prove the second claim, note that  $\dim (P_q\cap L_{\beta})=r-1>0$.  
Furthermore, note that 
\[
\rho_{\beta} =\exp(-2\pi i qv_r)\cdot \exp(2\pi i q(v_r+\tfrac{\lambda_{\beta}}{q} ) )
\in \exp(L_{\beta}\otimes \C)\cdot \exp(P_q\otimes \C).
\]
Corollary \ref{cor:exp intersect} now implies that 
$P_q\in  \sigma_{r}(L_{\beta},\rho_{\beta})$, and we are done. 
\end{proof}

In particular, if the only positive-dimensional component of 
$\WW^i(X)$ is a torsion-trans\-lated subtorus of dimension 
$d\ge 2$, then the sets $\Omega^i_2(X), \dots , \Omega^i_d(X)$ 
are not open.

\subsection{Examples of non-open $\Omega$-invariants}
\label{subsec:omega not open}

We now isolate a situation where we can explicitly 
identify certain Dwyer--Fried sets which are not open subsets 
of the rational Grassmannian. 

\begin{prop}
\label{prop:ttdf}
Let $X$ be a CW-complex with finite $k$-skeleton, 
and set $n=b_1(X)$.  Suppose that, for some $i\le k$, 
there is an $(n-1)$-dimensional subspace $L\subset H^1(X,\Q)$
such that $\WW^i(X)=\big(\bigcup_{\alpha} \rho_{\alpha} T\big) \cup Z$, 
where $Z$ is a finite set, $T=\exp(L\otimes \C)$ and 
$\rho_{\alpha}\notin T$, for all $\alpha$.   
Then 
\[
\Omega^i_r(X)=
\begin{cases} 
\QP^{n-1} & \text{if $r=1$,}\\ 
\Grass_r(L) & \text{if $1< r <n$,}\\
\emptyset & \text{if $r\ge n$.}
\end{cases}
\] 
In particular, $\Omega^i_r(X)$ is not open for $1< r <n$, 
and $\Omega^i_{n-1}(X)=\{L\}$.
\end{prop} 

\begin{proof}
First note that $\tau_1(\WW^i(X))=\{0\}$, since 
$1\notin \rho_{\alpha} T$; thus, by Theorem \ref{thm:df1}, 
$\Omega^i_1(X)=\QP^{n-1}$. The fact  that 
$\Omega^i_n(X)=\emptyset$ follows from 
Theorem \ref{thm:max abel}. 

Now assume $2\le r\le n-1$, and  
suppose $P$ is an $r$-plane in $\Q^n$.  If  $P\subsetneqq L$, 
then $P+L=\Q^n$ and $\dim (P\cap L)\ge r-1>0$; 
thus, $P\in \sigma_r(L,\rho_{\alpha})$. Conversely, if $P\subseteq L$, 
then clearly $P\notin \sigma_r(L,\rho_{\alpha})$.  The desired 
conclusion follows from Theorem \ref{thm:df translated}.
\end{proof}
 
We illustrate this proposition with a concrete example%
---as far as we know, the first of its kind---%
of a finitely presented group $G$ for which the set 
$\Omega^1_2(G)$ is not open.  

\begin{example}
\label{ex:omega closed}
Consider the group $G$ with generators $x_1, x_2, x_3$ 
and relators $r_1= [x_1^2,x_2]$, 
$r_2=[x_1,x_3]$,  
$r_3=x_1[x_2,x_3]x_1^{-1}[x_2,x_3]$.  
Note that $G_{\ab}=\Z^3$. 
Computing Fox derivatives, we find that 
\[
\partial_2^{\ab} = \begin{pmatrix}
(x_2-1) (1+x_1)& (1-x_1)(1+x_1) & 0\\
x_3-1 & 0 & 1-x_1\\
0 & (x_3-1) (1+x_1) & (1-x_2) (1+x_1)
\end{pmatrix}.
\]

Identifying $\wG=(\C^{\times})^3$, we obtain  
$\WW^1(G)= \{1\}\cup \set{ t \in  (\C^{\times})^3 \mid t_1=-1}$. 
By Proposition \ref{prop:ttdf}, the set $\Omega^1_2(G)$ 
consists of a single point in $\Grass_2(\Q^3)$, 
corresponding to the plane  $x_1=0$. 
In particular, $\Omega^1_2(G)$ is not open, 
not even in the usual topology on $\QP^2$.
\end{example}

\section{K\"{a}hler manifolds}
\label{sect:proj}

In this section, we discuss the characteristic varieties  
(also known in this context as the Green--Lazarsfeld sets) 
and the Dwyer--Fried sets of compact K\"{a}hler manifolds, 
highlighting the manner in which orbifold fibrations determine 
these sets in degree $1$.  

\subsection{Maps to $2$-orbifolds}
\label{subsec:pencils} 

Let $\Sigma_g$ be a Riemann surface of genus $g\ge 1$. 
As we saw in Example \ref{ex:cv surf}, the characteristic 
varieties of $\Sigma_g$ either fill the whole character torus 
$H^1(\Sigma_g,\C^{\times})=(\C^{\times})^{2g}$, or consist 
only of the identity.   In particular, if $g>1$, then 
$\VV^1(\Sigma_g)=(\C^{\times})^{2g}$.  

Following the approach of Dimca \cite{Di07}, 
Delzant \cite{Dz08}, Campana \cite{Cm11}, and 
Artal Bartolo, Cogolludo, and Matei \cite{ACM}, 
let us consider the more general situation 
of $2$-orbifolds.  Fix points $q_1,\dots ,q_t$ in $\Sigma_g$, 
and assign to these points integer weights $m_1,\dots , m_t$ 
with $m_i\ge 2$.  The orbifold fundamental group 
$\Gamma=\pi_1^{\orb}(\Sigma_g, \mathbf{m})$ 
associated to these data may be presented as
\begin{equation}
\label{eq:orbipi1}
\Gamma=
\left\langle 
\begin{array}{c}
x_1,\dots, x_g, y_1,\dots , y_g, \\[2pt]
z_1, \dots ,z_t 
\end{array}\left|
\begin{array}{c}
[x_1,y_1]\cdots [x_g,y_g] z_1\cdots z_t =1, \\[2pt]
z_1^{m_1}=\cdots =z_t^{m_t}=1
\end{array}
\right.
\right\rangle. 
\end{equation}

Clearly, we have an epimorphism 
$\kappa\colon\Gamma\surj \pi_1(\Sigma_g)$, 
obtained by sending $z_i\mapsto 1$.  Upon abelianizing, 
we obtain an isomorphism 
$\Gamma_{\ab} \cong \pi_1(\Sigma_g)_{\ab} \oplus A$, 
where 
\begin{equation}
\label{eq:tors gab}
A:=\Tors(\Gamma_{\ab})=
\Z_{m_1}\oplus \cdots \oplus \Z_{m_t}/(1,\dots ,1).
\end{equation}
Evidently, the group $A$ has order 
$m_1\cdots m_t/\lcm(m_1,\dots ,m_t)$; in particular, if 
$t=0$ or $t=1$, then $A=0$, but otherwise, $A\ne 0$. 

Identify $\widehat{\Gamma}=\widehat{\Gamma}^{\circ}  
\times \widehat{A}$, where $\widehat{\Gamma}^{\circ}=
\hat{\kappa}\big(\widehat{\pi_1(\Sigma_g)}\big) \cong 
(\C^{\times})^{2g}$.  A Fox calculus computation as in 
\cite[Proposition 3.11]{ACM} shows that
\begin{equation}
\label{eq:v1piorb}
\VV^1(\Gamma)=\begin{cases}
\widehat{\Gamma} & \text{if $g\ge 2$}, \\[2pt]
\big( \widehat{\Gamma}\setminus 
\widehat{\Gamma}^{\circ}\big) \cup  \{1\} 
& \text{if $g=1$ and $t> 1$}, \\[2pt]
\{1\}  
& \text{if $g=1$ and $t \le 1$}.
\end{cases}
\end{equation}

\begin{prop}
\label{prop:orbi}
Let $G$ be a finitely generated group, and suppose there 
is an epimorphism $\varphi\colon G\surj \Gamma$, 
where $\Gamma=\pi_1^{\orb}(\Sigma_g, (m_1,\dots, m_t))$, 
with either $g\ge 2$, or $g=1$ and $t\ge 2$. 
Set $A=\Tors(\Gamma_{\ab})$. Then 
\begin{equation}
\label{eq:v1gt}
\VV^1(G)\supseteq \bigcup_{\rho} \rho \cdot T,
\end{equation}
where $T$ is the $(2g)$-dimensional subtorus of $\wG$ 
obtained by pulling back $\widehat{\pi_1(\Sigma_g)}$ 
along the map $\varphi_0=\kappa\circ \varphi \colon G\surj \pi_1(\Sigma_g)$,  
and the union is taken over the set $\hat{\varphi}(\widehat{A})$ in 
the first case, and $\hat{\varphi}(\widehat{A})\setminus \{1\}$ 
in the second case. 
\end{prop}

\begin{proof}  
By Lemma \ref{lem:epi cv}, the induced morphism 
on character groups, 
$\hat{\varphi} \colon \widehat{\Gamma} \inj \wG$, embeds the 
characteristic variety $\VV^1(\Gamma)$ into $\VV^1(G)$.  
Applying formula \eqref{eq:v1piorb} and the discussion 
preceding it ends the proof.
\end{proof}

Note that, in either case, the right side of \eqref{eq:v1gt} is a 
finite union of torsion-translated subtori of $\wG$.

\subsection{The Green--Lazarsfeld sets}
\label{subsec:cv kahler}
Let $M$ be a compact, connected, K\"{a}hler 
manifold, for instance, a smooth, complex projective variety.   
The structure of the characteristic varieties of such manifolds 
was determined by Green and Lazarsfeld in \cite{GL87, GL91}, 
building on work of Castelnuovo and de Franchis, 
Beauville \cite{Be92}, and Catanese \cite{Cat91}.   
For this reason, the varieties $\VV^i_d(M)$ are 
also known in this context as the 
{\em Green--Lazarsfeld sets}\/ of $M$. 

The theory was further amplified by Simpson \cite{Si93}, 
Ein and Lazarsfeld \cite{EL97},  Arapura \cite{Ar}, and 
Campana \cite{Cm01}, with some of the latest developments 
appearing in \cite{Di07, Dz08, DPS-duke, Cm11, ACM}.  
The relationship between the present definition of 
characteristic varieties and the original definition from 
\cite{GL87} is discussed at length by Budur in \cite{Bu09}. 
The basic nature of the Green--Lazarsfeld sets 
is summarized in the following theorem. 

\begin{theorem}[\cite{GL91, Be92, Si93, Ar}]
\label{thm:cv kahler}
Each characteristic variety $\VV^i_d(M)$ of a 
compact K\"{a}hler manifold $M$ 
is a finite union of unitary translates of algebraic 
subtori of $H^1(M,\C^{\times})$.  Furthermore, if $M$ is 
projective, then all the translates are by torsion 
characters.
\end{theorem}

In degree $i=1$, the structure of the Green--Lazarsfeld 
set $\VV^1(M)=\VV^1_1(M)\cup \{1\}$ can be made more 
precise. First, we need some background on orbifold fibrations
(also known as orbifold morphisms, or pencils).

As before, let $\Sigma_g$ be a Riemann surface of genus 
$g\ge 1$, with marked points $q_1,\dots, q_t$, and weight 
vector $\mathbf{m}=(m_1,\dots, m_t)$, where $m_i\ge 2$ and 
$\abs{\mathbf{m}}:=t\ge 0$.   
A surjective map $f\colon M\to (\Sigma_g,\mathbf{m})$ 
is called an {\em orbifold fibration}\/ if $f$ is holomorphic, 
the fiber over any non-marked point is connected, 
and, for every point $q_i$, the multiplicity of the fiber 
$f^{-1}(q_i)$ equals $m_i$.  Such a map induces 
an epimorphism $f_{\sharp} \colon \pi_1(M) \surj \Gamma$, 
where $\Gamma=\pi_1^{\orb}(\Sigma_g, \mathbf{m})$ is 
the orbifold fundamental group described in \eqref{eq:orbipi1}.
By Proposition \ref{prop:orbi}, the induced morphism of 
character groups, 
$\widehat{f_{\sharp}}\colon \widehat{\Gamma}\inj \widehat{\pi_1(M)}$, 
sends $\VV^1(\Gamma)$ to a union of (possibly torsion-translated) 
subtori inside $\VV^1(M)$.

Two orbifold fibrations, $f\colon M\to (\Sigma_g,\mathbf{m})$ and 
$f'\colon M\to (\Sigma_{g'},\mathbf{m}')$, are equivalent 
if there is a biholomorphic map $h\colon \Sigma_g\to \Sigma_{g'}$ 
which sends marked points to marked points, while preserving multiplicities. 
Write $\chi^{\orb}(\Sigma_g, \mathbf{m})=2-2g-\sum_{i=1}^{t} (1-1/m_i)$. 
As shown by Delzant \cite{Dz08}, a K\"ahler manifold $M$ admits 
only finitely many equivalence classes of orbifold fibrations  
for which the orbifold Euler characteristic of the base is negative.

The next theorem, which is a distillation of several results 
from the quoted sources, 
shows that {\em all}\/ positive-dimensional components in the 
first characteristic variety of $M$ arise by pullback along this 
finite set of pencils.  

\begin{theorem}[\cite{Be92, Ar, Cm01, Di07, Dz08, Cm11, ACM}]
\label{thm:cv1 kahler}
Let $M$ be a compact K\"{a}hler manifold.  Then 
\begin{equation}
\label{eq:vvk}
\VV^1(M)=\bigcup_{\alpha} \widehat{(f_{\alpha})_{\sharp}}
\big(\VV^1(\pi_1^{\orb}(\Sigma_{g_\alpha},\mathbf{m}_{\alpha}) \big) 
\cup Z, 
\end{equation}
where $Z$ is a finite set of torsion characters, and the union runs 
over the (finite) set of equivalence classes of orbifold fibrations 
$f_{\alpha}\colon M\to (\Sigma_{g_\alpha},\mathbf{m}_{\alpha})$ 
with either $g_{\alpha}\ge 2$, or $g_{\alpha}=1$ and 
$\abs{\mathbf{m}_{\alpha}}\ge 2$. 
\end{theorem}

In particular, every positive-dimensional component of $\VV^1(M)$ 
is of the form $\rho\cdot T$, with  $T$ a $(2g)$-dimensional algebraic 
subtorus in $H^1(M,\C^{\times})$, and $\rho$ of finite order (modulo $T$).  
Clearly,  if $\rho\in T$, then $g\ge 2$.  Moreover, if $\rho\notin T$ 
and $g\ge 2$, then $T$ is also a component of $\VV^1(M)$, but if 
$g=1$, then $T$ is {\em not}\/ a component of $\VV^1(M)$. 
Let us record a particularly simple situation as a corollary.  

\begin{corollary}
\label{cor:kahler nomult}
Let $M$ be a compact K\"{a}hler manifold, and suppose $M$ 
admits no orbifold fibrations with multiple fibers.   Then 
\begin{equation}
\label{eq:vv nomult}
\VV^1(M)=\bigcup_{\alpha} 
f_{\alpha}^* ( H^1(\Sigma_{g_{\alpha}},\C^{\times}))\cup Z, 
\end{equation}
where $Z$ is a finite set of torsion characters, 
and the union runs over the set of equivalence classes of 
orbifold fibrations $f_{\alpha}\colon M\to \Sigma_{g_{\alpha}}$ 
with $g_{\alpha}\ge 2$, and 
$f_{\alpha}^*\colon  H^1(\Sigma_{g_{\alpha}},\C^{\times})\inj 
H^1(M,\C^{\times})$ is the induced homomorphism in cohomology.
\end{corollary}

Another simple situation is the one considered by 
Green and Lazarsfeld \cite{GL87} and Hacon 
and Pardini \cite{HP} in a slightly less general setting.

\begin{prop}
\label{prop:kahler zerores}
Let $M$ be a compact K\"{a}hler manifold, and suppose 
there is no orbifold fibration $M\to \Sigma_g$ with $g\ge 2$.  
Then 
\begin{enumerate}
\item \label{gl1}  $\dim \VV^1(M)=0$.  
\item \label{gl2}  If, moreover, $p_g(M) = q(M) = 3$, 
then $\VV^1(M)=\{1\}$.  
\end{enumerate}
\end{prop}

\begin{proof}
Part \eqref{gl1} follows directly from Theorem \ref{thm:cv1 kahler} 
and the discussion following it. 

Part \eqref{gl2} follows from the proof of \cite[Proposition 2.8]{HP}, 
which crucially relies on Theorem 1.2 (1.2.3) from \cite{EL97}.
\end{proof}

\begin{example}
\label{ex:symm}
Following \cite{CCM}, consider the symmetric product 
$M=(\Sigma_3\times \Sigma_3)/\sigma$, where $\sigma$ 
is the involution interchanging the two factors. 
Then $M$ is a minimal surface of general type with 
$p_g(M) = q(M) = 3$ and $K_M^2=6$.  As noted 
in \cite{HP}, this surface has no irrational pencils 
of genus $g\ge 2$, and thus $\VV^1(M)=\{1\}$. 
\end{example}

\subsection{Characteristic subspace arrangements}
\label{subsec:tau1 kahler}
As an immediate corollary to Theorem \ref{thm:cv1 kahler}, 
we also obtain the following characterization of the exponential 
tangent cone to the first Green--Lazarsfeld set of $M$.

\begin{corollary}
\label{cor:tcone kahler}
Let $M$ be a compact K\"{a}hler manifold.  Then 
\begin{equation}
\label{eq:omega kahler}
\tau_1(\WW^1(M))=
\bigcup_{\alpha} f_{\alpha}^* (H^1(\Sigma_{g_{\alpha}},\C)),
\end{equation}
where the union runs over the set of equivalence classes of 
orbifold fibrations $f_{\alpha}\colon M\to \Sigma_{g_{\alpha}}$ 
with $g_{\alpha}\ge 2$, and $f_{\alpha}^*\colon  H^1(\Sigma_{g_{\alpha}},\C)\inj 
H^1(M,\C)$ is the induced homomorphism in cohomology.
\end{corollary}

Put differently, the characteristic arrangement $\CC_1(M)$ 
consists of the linear subspaces 
$f_{\alpha}^* (H^1(\Sigma_{g_{\alpha}},\Q)) \subseteq H^1(M,\Q)$ 
arising by pullback along orbifold fibrations with base 
genus $g_{\alpha}\ge 2$.  Note that each of these subspaces 
has dimension $2g_{\alpha}\ge 4$.  

It follows from 
\cite{EL97, DPS-duke} that $\tau_1(\WW^1(M))$ 
coincides with the resonance variety $\RR^1(M,\C)$ 
defined in \eqref{eq:intro res}; thus, the arrangement 
$\CC_1(M)$ depends only on the cup-product map 
$H^1(M,\Q)\wedge H^1(M,\Q) \to H^2(M,\Q)$. 
In the projective setting, a bit more can be said. 

\begin{theorem}[\cite{DPS-imrn}]
\label{thm:qp cv}
Let $M$ be a smooth, complex projective variety, and let 
$\rho_{\alpha} T_{\alpha}$ and $\rho_{\beta} T_{\beta}$ 
be two distinct components of $\VV^1(M)$.  Then either 
$T_{\alpha} = T_{\beta}$, or 
$\tau_1 (T_{\alpha}) \cap \tau_1 (T_{\beta})= \{ 0\}$. Furthermore, 
$\rho_{\alpha} T_{\alpha} \cap \rho_{\beta} T_{\beta}$ 
is a finite set of torsion characters. 
\end{theorem}

Therefore, the linear subspaces 
$\tau_1 (T_{\alpha}) =f_{\alpha}^* (H^1(\Sigma_{g_{\alpha}},\Q))$
comprising the arrangement $\CC_1(M)$ intersect pairwise 
transversely. 

\subsection{Dwyer--Fried sets}
\label{subsec:df kahler}

The above results of the Green--Lazarsfeld sets of compact 
K\"{a}hler manifolds $M$ can be used to describe 
the Dwyer--Fried invariants of such manifolds.   
For instance, if $\WW^i(M)$ contains no positive-dimensional 
translated subtori, then Corollary \ref{cor:tau schubert} and 
Theorem \ref{thm:cv kahler} insure that 
$\Omega^i_r(M) = \sigma_r(\tau^{\Q}_1(\WW^i(M)))^{\compl}$, 
for all $r\ge 1$.  When $i=1$, we can be much more concrete.

Given an orbifold fibration $f\colon M\to (\Sigma_{g},\mathbf{m})$, 
write $L=f^* (H^1(\Sigma_{g},\Q))$ and 
$A=\Tors(\pi_1^{\orb}(\Sigma_{g},\mathbf{m}))$. 
Furthermore, denote by 
$\varphi\colon \pi_1(M)\surj  \pi_1^{\orb}(\Sigma_{g},\mathbf{m})$ 
and $\varphi_0\colon \pi_1(M)\to  \pi_1(\Sigma_{g})$ the 
homomorphisms induced by $f$, and identify 
$\widehat{\varphi_0}(\widehat{\pi_1(\Sigma_{g})})=
f^*(H^1(\Sigma_g,\C^{\times}))$ 
with the subtorus $T=\exp(L\otimes \C)$.  

\begin{theorem}
\label{thm:dfk}
Let $M$ be a compact K\"{a}hler manifold. 
For all $r\ge 1$, 
\begin{equation}
\label{eq:om1k}
\Omega^1_r(M) = \Grass_r(H^1(M,\Q)) \setminus 
\bigcup_{\alpha} \bigcup_{\rho} \sigma_r(L_\alpha, \rho),
\end{equation}
where the first union runs over the set of equivalence classes of 
orbifold fibrations 
$f_{\alpha}\colon M\to (\Sigma_{g_\alpha},\mathbf{m}_{\alpha})$ 
with either $g_{\alpha}\ge 2$, or $g_{\alpha}=1$ and 
$\abs{\mathbf{m}_{\alpha}}\ge 2$, while the second union 
runs over $\rho\in \hat{\varphi}_{\alpha}(\widehat{A}_{\alpha})$ and  
$\rho\in\hat{\varphi}_{\alpha}(\widehat{A}_{\alpha})\setminus \{1\}$,  
respectively. 
\end{theorem} 

\begin{proof}
Follows from Theorem  \ref{thm:df translated} on one hand, 
and Proposition \ref{prop:orbi} and Theorem \ref{thm:cv1 kahler} 
on the other hand.
\end{proof}

Next, we give a cohomological upper bound for the Dwyer--Fried 
sets $\Omega^1_r(M)$, and single out a couple of situations   
when that bound is attained.  

\begin{theorem}
\label{thm:df kahler}
Let $M$ be a compact K\"{a}hler manifold.  
Then, for all $r\ge 1$, 
\begin{equation}
\label{eq:omega1 kahler}
\Omega^1_r(M) \subseteq \Grass_r(H^1(M,\Q))\setminus 
\bigcup_{\alpha} \sigma_r(f_{\alpha}^* (H^1(\Sigma_{g_{\alpha}},\Q)),
\end{equation}
where the union runs over the set of equivalence classes of 
orbifold fibrations $f_{\alpha}\colon M\to \Sigma_{g_{\alpha}}$ 
with $g_{\alpha}\ge 2$.  Moreover, inclusion \eqref{eq:omega1 kahler} 
holds as equality, provided either 
\begin{enumerate}
\item \label{c1} $r=1$, or
\item \label{c2} $M$ admits no orbifold fibrations with multiple fibers.
\end{enumerate}
\end{theorem}

\begin{proof}
Formula \eqref{eq:omega1 kahler} is a direct consequence of 
Theorem  \ref{thm:sch bound} and Corollary \ref{cor:tcone kahler}. 

Claim \eqref{c1} now follows from Theorem \ref{thm:df1}, while 
Claim \eqref{c2} follows from Corollaries \ref{cor:tau schubert} 
and \ref{cor:kahler nomult}.
\end{proof}

In the situations  \eqref{c1} and  \eqref{c2} recorded above, 
each Dwyer--Fried set $\Omega^1_r(M)$ is 
the complement of a finite union of special Schubert varieties; 
in particular, an open subset of $\Grass_r(H^1(M,\Q))$.  
Other instances when this happens are recorded in the 
next proposition.

\begin{prop}
\label{prop:dfk empty} 
Let $M$ be a compact K\"{a}hler manifold.  
\begin{enumerate}
\item \label{kg1}
If $M$ admits an orbifold fibration with base genus $g\ge 2$,  
then $\Omega^1_r(M) = \emptyset$, for all $r> b_1(M)-2g$.
\item \label{kg2}
If $M$ does not admit an orbifold fibration with 
base genus $g\ge 2$, then  $\Omega^1_r(M)=\Grass_r(H^1(M,\Q))$, 
for all $r\ge 1$.
\end{enumerate}
\end{prop}

\begin{proof}
For part \eqref{kg1}, use Corollaries \ref{cor:df sch} and 
\ref{cor:tcone kahler}, while for part \eqref{kg2}, use 
Propositions \ref{prop:df smallcv} and \ref{prop:kahler zerores}.
\end{proof}
 
In general, though, the presence of (elliptic) pencils with 
multiple fibers drastically changes the nature of the 
$\Omega$-invariants. 

\begin{prop}
\label{prop:dfk not open}
Let $M$ be a smooth, complex projective variety, and 
suppose $M$ admits an orbifold fibration  with multiple fibers 
and base genus $g=1$. Then $\Omega^1_2(M)$ is {\em not} 
an open subset of $\Grass_2(H^1(M,\Q))$. 
\end{prop}

\begin{proof}
By hypothesis,  there is an orbifold fibration, 
$f\colon M\to (\Sigma_1,\mathbf{m})$,  
with $\abs{\mathbf{m}}\ge 2$. Thus, $\WW^1(M)$ 
contains a component of the form  $\rho T$, where 
$T=f^*(H^1(\Sigma_1,\C^{\times}))$ is a $2$-dimensional 
subtorus, and $\rho$ is a finite-order character (modulo $T$), 
$\rho\notin T$, yet $T$ itself is not a component of $\WW^1(M)$. 

By Theorem \ref{thm:cv1 kahler}, all positive-dimensional 
components of $\WW^1(M)$ are torsion-trans\-lated subtori, 
of the form $\rho_\alpha T_\alpha$.  
If $T_{\alpha}=T$, then the component $\rho_\alpha T_\alpha$ 
must arise from an elliptic pencil; thus, $\rho_{\alpha}\notin T_\alpha$. 
On the other hand, if $T_{\alpha}\ne T$, then 
Theorem \ref{thm:qp cv} insures that 
$\tau_1(T_{\alpha})\cap \tau_1(T)=\{0\}$.  Thus, all conditions 
of Theorem \ref{thm:df tors translated} are satisfied, 
and the desired conclusion follows.
\end{proof}

\subsection{The Catanese--Ciliberto--Mendes Lopes surface}
\label{subsec:kahler not open}
We now give a concrete example of a smooth, complex projective  
variety $M$ for which one of the Dwyer--Fried sets is not open.  
The variety in question is a minimal surface of general type with 
$p_g(M) = q(M) = 3$ and $K_M^2=8$.  The construction 
goes back to Catanese, Ciliberto, and Mendes Lopes \cite{CCM}; 
it was revisited by Hacon and Pardini \cite{HP}, and more 
recently, by Akhmedov and Park \cite{AP}, as well as Akbulut \cite{Ak}. 

\begin{figure}[t]
\setlength{\unitlength}{6pt}
\begin{picture}(28,19)(4,-8.2)
\xy <0mm,0mm>;<7pt,0pt>:<0pt,7pt>::
,*\xycircle(4,8){\dir{-}}
@={(-4.0,0.5),(4.0,0.5)},
s0="prev" @@{;"prev";**@{.}="prev"}
@={(-5.8,0.5),(-4.0,0.5)},
s0="prev" @@{;"prev";**@{-}="prev"}
@={(4.0,0.5),(4.8,0.5)},
s0="prev" @@{;"prev";**@{-}="prev"}
@={(5.3,0.5),(6,0.5)},
s0="prev" @@{;"prev";**@{-}="prev"}
\POS (-4,0.5) ,*\xycircle(0.16,0.16){}  
,*\xycircle(0.12,0.12){}  ,*\xycircle(0.06,0.06){} 
\POS (4,0.5)  ,*\xycircle(0.16,0.16){}  
,*\xycircle(0.12,0.12){}  ,*\xycircle(0.06,0.06){} 
,(5.9,-0.9)*+{}; (7,0.8)*+{}
**\crv{(3.5,2.2)} ?>*\dir{>};
(6,-1.6) *+{\sigma_1}
,(-1.5,3)*+{}; (1.5,3)*+{}
**\crv{(0,4.1)}
,(-1.6,3.95)*+{}; (1.6,3.95)*+{}
**\crv{(0,1.55)}
,(-1.5,-2)*+{}; (1.5,-2)*+{}
**\crv{(0,-0.9)}
,(-1.6,-1.05)*+{}; (1.6,-1.05)*+{}
**\crv{(0,-3.45)}
\endxy
\hspace*{3pc}
\xy <0mm,0mm>;<7pt,0pt>:<0pt,7pt>::
,*\xycircle(4,8){\dir{-}}
@={(1.4,0.5),(3.8,0.5)},
s0="prev" @@{;"prev";**@{.}="prev"}
@={(-5.8,0.5),(0.5,0.5)},
s0="prev" @@{;"prev";**@{-}="prev"}
@={(4.3,0.5),(4.8,0.5)},
s0="prev" @@{;"prev";**@{-}="prev"}
@={(5.3,0.5),(6,0.5)},
s0="prev" @@{;"prev";**@{-}="prev"}
,(5.9,-0.9)*+{}; (7,0.8)*+{}
**\crv{(3.5,2.2)} ?>*\dir{>};
(6,-1.6) *+{\sigma_2}
,(-1.5,4)*+{}; (1.5,4)*+{}
**\crv{(0,5.1)}
,(-1.6,4.95)*+{}; (1.6,4.95)*+{}
**\crv{(0,2.55)}
,(-1.5,0.2)*+{}; (1.5,0.2)*+{}
**\crv{(0,1.3)}
,(-1.6,1.15)*+{}; (1.6,1.15)*+{}
**\crv{(0,-1.25)}
,(-1.5,-3.2)*+{}; (1.5,-3.2)*+{}
**\crv{(0,-2.1)}
,(-1.6,-2.25)*+{}; (1.6,-2.25)*+{}
**\crv{(0,-4.65)}
\endxy
\end{picture}
\caption{The complex curves $C_1$ and $C_2$, with involutions 
$\sigma_1$ and $\sigma_2$}
\label{fig:g2g3}
\end{figure}
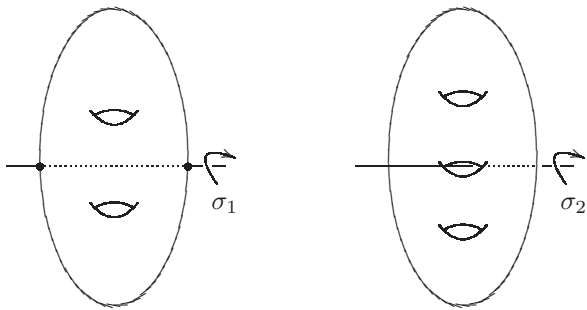

\begin{example}
\label{ex:omega kahler}
Let $C_1$ be a (smooth, complex) curve of genus $2$ 
with an elliptic involution $\sigma_1$ and $C_2$ a curve of 
genus $3$ with a free involution $\sigma_2$ (see Figure \ref{fig:g2g3}). 
Then $\Sigma_1=C_1/\sigma_1$ is a curve of genus $1$, 
and $\Sigma_2=C_2/\sigma_2$ is a curve of genus $2$. 
We let $\Z_2$ act freely on the product $C_1 \times C_2$ 
via the involution $\sigma_1\times \sigma_2$, and denote  
by $M$ the quotient surface.  
Projection onto the first coordinate yields an orbifold fibration  
$f_1 \colon M\to (\Sigma_1,(2,2))$ with two multiple fibers, each of 
multiplicity $2$, while projection onto the second coordinate 
defines a holomorphic fibration $f_2\colon M\to \Sigma_2$: 
\[
\xymatrix{
C_2 \ar^(.45){/\sigma_2}[d]& C_1 \times C_2 \ar_(.55){\pr_2}[l]
\ar^(.45){/\sigma_1\times\sigma_2}[d] 
\ar^(.55){\pr_1}[r] & C_1\ar^(.45){/\sigma_1}[d]\\
\Sigma_2 &  M \ar^{f_1}[r] \ar_{f_2}[l] & \Sigma_1
}
\]

Using the procedure described in \cite{BCGP}, we find that 
the fundamental group of $M$ has presentation 
with generators $x_1, \dots,  x_6$ and relators 
\begin{align*}
&[x_3^2,x_1],  \quad [x_3^2,x_2], \quad 
[x_2,x_1] [x_2^{x_3},x_1^{x_3}], \quad 
[x_3,x_4] [x_5,x_6], 
\\
&[x_1,x_4],  \quad [x_2,x_4], \quad
[x_1,x_5],  \quad [x_2,x_5],  \quad 
[x_1,x_6],  \quad [x_2,x_6], 
\\
&[x_1^{x_3}, x_4], \quad  [ x_2^{x_3}, x_4],\quad 
[x_1^{x_3}, x_5 ],  \quad [ x_2^{x_3}, x_5 ], \quad 
[ x_1^{x_3}, x_6 ], \quad [x_2^{x_3}, x_6 ],
\end{align*}
where $z^w:=w^{-1}zw$.  
Identify the character group $\Hom(\pi_1(M),\C^{\times})$ with 
$(\C^{\times})^6$.   A straightforward computation shows 
that the characteristic variety $\VV^1(M$) has two  
components, corresponding to the above two pencils.    
The first component is the subtorus $T_1=\exp(L_1\otimes \C)$, 
with $L_1=\{x\in \Q^6 \mid x_1=x_2=0\}$, while the second 
component is the translated subtorus $\rho T_2$, where 
$\rho = (1,1,-1,1,1,1)$ and $T_2=\exp(L_2\otimes \C)$, with  
$L_2=\{x\in \Q^6 \mid x_3=x_4=x_5=x_6=0\}$.  

By Theorem \ref{thm:dfk}, we have that  
$\Omega^1_2(M)=\Grass_2(\Q^6)\setminus 
(\sigma_2(L_1)\cup \sigma_2(L_2,\rho))$.  
Clearly, the plane $L_2$ belongs to $\Omega^1_2(M)$. 
By Proposition \ref{prop:dfk not open}, the set 
$\Omega^1_2(M)$ is not open, not even 
in the usual topology on the Grassmannian. 
\end{example}

\begin{remark}
\label{rem:hp}
In \cite{HP}, Hacon and Pardini detect the component $T_1$ 
of the characteristic variety $\VV^1(M)$ 
by means of the infinitesimal description of this 
variety, due to Green and Lazarsfeld \cite{GL87}. 
Although their method does not detect the translated 
component $\rho T_2$, it does distinguish topologically 
the surface $M$ from the symmetric product considered 
in Example \ref{ex:symm}.
\end{remark}

\section{Quasi-K\"{a}hler manifolds}
\label{sect:quasi proj}

In this section, we consider the more general setting 
of quasi-K\"{a}hler manifolds, which includes smooth, 
quasi-projective varieties such as complements 
of plane curves or hyperplane arrangements. 

\subsection{Non-compact $2$-orbifolds}
\label{subsec:orbi}

Let $\Sigma_{g,s}=\Sigma_g \setminus \{p_1,\dots , p_s\}$ be 
a Riemann surface of genus $g\ge 0$ with $s$ points removed 
($s\ge 1$).   
Then $\Sigma_{g,s}$ is a connected, smooth, quasi-projective 
variety, with $\pi_1(\Sigma_{g,s})=F_n$, where $n=b_1(\Sigma_{g,s})=2g+s-1$.  
Thus, $\VV^1(\Sigma_{g,s})=H^1(\Sigma_{g,s},\C^{\times})$, 
unless $g=0$ and $s\le 2$.  

More generally, let us consider as in \cite{Di07, ACM} 
a non-compact $2$-orbifold, i.e., a manifold $\Sigma_{g,s}$ 
with marked points $q_1,\dots, q_t$  
and weight vector $\mathbf{m}=(m_1,\dots, m_t)$.  
The orbifold fundamental group 
$\Gamma=\pi_1^{\orb}(\Sigma_{g,s}, \mathbf{m})$ 
associated to these data is the free product 
\begin{equation}
\label{eq:orbipi1 noncompact}
\Gamma=
F_{n} * \Z_{m_1} * \cdots * \Z_{m_t}.
\end{equation}
Thus, $\Gamma_{\ab}=\Z^{n}\oplus A$, 
where $A=\Z_{m_1}\oplus \cdots \oplus \Z_{m_t}$.  
Identify $\widehat{\Gamma}= (\C^{\times})^{n} \times \widehat{A}$.  
A Fox calculus computation as in \cite[Proposition 3.10]{ACM} 
shows that
\begin{equation}
\label{eq:v1pi}
\VV^1(\Gamma)=\begin{cases}
\widehat{\Gamma} & \text{if $n\ge 2$}, \\
\big( \widehat{\Gamma}\setminus \widehat{\Gamma}^{\circ}\big) \cup \{1\}
& \text{if $n=1$ and $t>0$}, \\
\{1\} & \text{if $n=1$ and $t=0$}.
\end{cases}
\end{equation}

In the above, note that $n=1$ if and only if $g=0$ and $s=2$, i.e., 
$\Sigma_{g,s}=\C^{\times}$. 

Proceeding as in the proof of Proposition \ref{prop:orbi}, we 
obtain the following result (see also \cite{Di07, ACM}). 

\begin{prop}
\label{prop:orbi noncompact}
Let $G$ be a finitely generated group, and suppose 
there is an epimorphism $\varphi\colon G\surj \Gamma$, 
where $\Gamma=\pi_1^{\orb}(\Sigma_{g,s}, (m_1,\dots, m_t))$, 
with $g\ge 0$ and $s\ge 1$.
Set $A=\Tors(\Gamma_{\ab})$ and $n=2g+s-1$, 
and suppose further that $n\ge 2$, or $n=1$ and $t>0$. 
Then $\VV^1(G)$ contains all the translated subtori in $\wG$ of 
the form $\rho \cdot T$, where $T$ is obtained by pulling back the torus 
$\widehat{F_n}=(\C^{\times})^n$ along the homomorphism 
$G\xrightarrow{\varphi}\Gamma \surj F_n$,  
and $\rho$ belongs to $\hat{\varphi}(\widehat{A})$ 
in the first case, and 
$\hat{\varphi}(\widehat{A})\setminus \{1\}$ in the second case. 
\end{prop}

A particular case is worth singling out. 

\begin{corollary}
\label{prop:za}
Suppose $G$ has a factor group of the form $\Z*\Z_m$, with $m\ge 2$.  
There is then a $1$-dimensional algebraic subtorus $T$ 
and a torsion character $\rho\in \wG \setminus T$ such that 
$\VV^1(G)$ contains the translated tori $\rho T,\dots, \rho^{m-1} T$.  
\end{corollary}

\begin{example}
\label{ex:1torus}
Let $G=\langle x_1, x_2 \mid x_1 x_2^2 = x_2^2 x_1 \rangle$, and  
identify $\wG=(\C^{\times})^2$, with coordinates $t_1$ and $t_2$. 
Taking the quotient of $G$ by the normal subgroup generated by 
$x_2^2$, we obtain an epimorphism $\varphi\colon G\surj \Z*\Z_2$. 
The morphism induced by $\varphi$ on character groups sends $\widehat{\Z}$ to 
$T=\{t_2=1\}$, and $\widehat{\Z}_2\setminus \{1\}$ to $\rho=(1,-1)$.  
By the above proposition, $\VV^1(G)$ contains the translated 
torus $\rho T=\{ t_2=-1\}$.  Direct computation with Fox derivatives 
shows that, in fact, $\VV^1(G)=\{1\}\cup \rho T$. 
\end{example}

\subsection{Arapura's theorem}
\label{subsec:arapura}

A smooth, connected manifold $X$ is said 
to be a {\em quasi-K\"{a}hler manifold}\/ 
if there is a compact K\"{a}hler manifold $\overline{X}$ and a 
normal-crossings divisor $D$ such that $X=\overline{X}\setminus D$. 
The most familiar examples of quasi-K\"{a}hler manifolds are 
smooth, quasi-projective complex varieties, 
such as complements of projective hypersurfaces, or 
configuration spaces of smooth projective varieties. 
As is well-known, every quasi-projective variety has 
the homotopy type of a finite CW-complex. 

The basic structure of the characteristic varieties of 
quasi-K\"{a}hler manifolds was described by D.~Arapura. 

\begin{theorem}[\cite{Ar}]
\label{thm:arapura}
Let $X=\overline{X}\setminus D$ be a quasi-K\"{a}hler manifold.  
If either $D=\emptyset$ or $b_1(\overline{X})=0$, 
then each characteristic variety $\VV^i_d(X)$ 
is a finite union of unitary translates of algebraic 
subtori of $\widehat{\pi_1(X)}$.
\end{theorem}

In the above, each component of $\VV^i_d(X)$ is 
of the form $\rho\cdot  \widehat{f_{\sharp}} (\widehat{\pi_1(T)})$, 
for some character $\rho\colon \pi_1(X)\to S^1$, 
and some holomorphic map $f\colon M\to T$ to 
a complex Lie group $T$ which decomposes as 
a product of factors of the form $\C^{\times}$ and  
$S^1\times S^1$. 

When $X$ is quasi-projective and $i=1$, a more precise 
structure theorem holds.   Since the compact case was 
already treated in \S\ref{sect:proj}, we will focus  
for the rest of this section on the non-compact case.

Let $(\Sigma_{g,s},\mathbf{m})$ be a Riemann surface of genus 
$g\ge 0$, with $s\ge 0$ points removed (so that 
$\Sigma_{g,0}=\Sigma_g$), and with marked points 
$(q_1,m_1),\dots, (q_t,m_t)$.   A surjective, holomorphic map 
$f\colon X\to \Sigma_{g,s}$ is called an {\em orbifold fibration}\/ 
(or, a pencil) if the generic fiber is connected,
the multiplicity of the fiber over each $q_i$ equals $m_i$,  
and $f$ has an extension to the respective compactifications, 
$\bar{f}\colon \overline{X}\to \Sigma_{g}$, which is 
also a surjective, holomorphic map with connected 
generic fibers. 

An orbifold fibration as above induces an 
epimorphism $f_{\sharp} \colon \pi_1(X) \surj \Gamma$, 
where $\Gamma=\pi_1^{\orb}(\Sigma_{g,s}, \mathbf{m})$ is the 
orbifold fundamental group described in \eqref{eq:orbipi1} 
if $s=0$ and \eqref{eq:orbipi1 noncompact} if $s>0$.  
By Propositions \ref{prop:orbi} and \ref{prop:orbi noncompact}, 
the induced morphism, 
$\widehat{f_{\sharp}}\colon \widehat{\Gamma}\inj \widehat{\pi_1(X)}$, 
sends $\VV^1(\Gamma)$ to a union of (possibly torsion-translated) 
subtori inside $\VV^1(X)$.

\begin{theorem}[\cite{Ar, Di07, ACM}]
\label{thm:arapura bis}
Let $X$ be a smooth, quasi-projective variety.  Then
\begin{equation}
\label{eq:abis}
\VV^1(X)=\bigcup_{\alpha} \widehat{(f_{\alpha})_{\sharp}}
\big(\VV^1(\pi_1^{\orb}(\Sigma_{g_\alpha, s_{\alpha}},
\mathbf{m}_{\alpha}) \big) 
\cup Z, 
\end{equation}
where $Z$ is a finite set of torsion characters, and 
the union runs over a finite set of orbifold fibrations 
$f_{\alpha}\colon X\to (\Sigma_{g_\alpha, s_{\alpha}},\mathbf{m}_{\alpha})$. 
\end{theorem}

In particular, every positive-dimensional component of $\VV^1(X)$ 
is of the form $\rho\cdot T$, where $T$ is an algebraic subtorus 
in $H^1(X,\C^{\times})$, and $\rho$ is of finite order (modulo $T$).  
If this component arises from an orbifold fibration with base 
$\Sigma_{g,s}$, then $T$ has dimension $n:=b_1(\Sigma_{g,s})$.  
Moreover, if $\rho\in T$, then $n=2g\ge 4$ or $n=2g+s-1\ge 2$, 
according to whether $s=0$ or not.   On the other hand, if $\rho\notin T$,  
then the direction torus $T$ is {\em not}\/ a component of $\VV^1(X)$ 
precisely when $n=2$ (if $s=0$) or $n=1$ (if $s>0$).  

More information on the pencils $f_{\alpha}$ occurring in \eqref{eq:abis} 
can be found in \cite[Proposition~7.2]{DPS-duke}.  In particular, 
if $b_1(\overline{X})=0$, then all those pencils  have base genus 
$g_{\alpha}=0$.

Next, we derive some immediate corollaries of the above theorem.

\begin{corollary}
\label{cor:qp nomult}
Let $X$ be a smooth, quasi-projective variety, and suppose $X$ 
admits no orbifold fibrations with multiple fibers.   Then 
\begin{equation}
\label{eq:vv qp nomult}
\VV^1(X)=\bigcup_{\alpha} 
f_{\alpha}^* ( H^1(\Sigma_{g_{\alpha},s_{\alpha}},\C^{\times}))\cup Z, 
\end{equation}
where $Z$ is a finite set of torsion characters, 
the union runs over the set of 
orbifold fibrations $f_{\alpha}\colon X\to \Sigma_{g_{\alpha},s_{\alpha}}$ 
with $2g_{\alpha}+s_{\alpha} \ge 3$, and 
$f_{\alpha}^*\colon  H^1(\Sigma_{g_{\alpha},s_{\alpha}},\C^{\times})\inj 
H^1(X,\C^{\times})$ is the induced homomorphism in cohomology.
\end{corollary}

\begin{corollary}
\label{cor:tcone qp}
Let $X$ be a smooth, quasi-projective variety. Then 
\begin{equation}
\label{eq:tau qp}
\tau_1(\WW^1(X))=
\bigcup_{\alpha} f_{\alpha}^* (H^1(\Sigma_{g_{\alpha},s_{\alpha}},\C)),
\end{equation}
where the union runs over all  
orbifold fibrations $f_{\alpha}\colon M\to \Sigma_{g_{\alpha},s_{\alpha}}$ 
with $2g_{\alpha}+s_{\alpha} \ge 3$, and 
$f_{\alpha}^*\colon  H^1(\Sigma_{g_{\alpha},s_{\alpha}},\C)\inj 
H^1(X,\C)$ is the induced homomorphism in cohomology.
\end{corollary}

The exact analogue of Theorem \ref{thm:qp cv} also holds 
for quasi-projective varieties. 

\subsection{Dwyer--Fried invariants}
\label{subsec:df qp}

Much as in the compact case, the above structural results 
inform on the $\Omega$-sets of (non-compact) quasi-K\"{a}hler manifolds.   

For instance, let $X=\overline{X}\setminus D$, where 
$\overline{X}$ is a compact K\"{a}hler manifold, 
$D$ is a normal-crossings divisor, and 
$b_1(\overline{X})=0$ if $D\ne \emptyset$. 
Furthermore, suppose $\WW^i(X)$ contains no positive-dimensional 
translated subtori. Then Corollary \ref{cor:tau schubert} and 
Theorem \ref{thm:arapura} insure that the inclusion 
$\Omega^i_r(X) \subseteq \sigma_r(\tau^{\Q}_1(\WW^i(X)))^{\compl}$ 
holds as equality, for all $r\ge 1$.  In general, though, the inclusion 
can very well be strict, as the next example illustrates.

\begin{example}
\label{ex:qp product}
Let $\mathcal{C}$ be the affine plane curve 
with equation $x^p-y^q=0$, where $p$ and $q$ are 
coprime positive integers.  Its complement, 
$Y=\C^2 \setminus \mathcal{C}$, is 
homotopy equivalent to the complement in $S^3$ 
of a $(p,q)$-torus knot $K$, whose Alexander polynomial 
is $\Delta_K=(t^{pq}-1)(t-1)(t^p-1)^{-1}(t^q-1)^{-1}$.

Now let $X=Y\times (\CP^1 \setminus \{\text{$n$ points}\})$, 
for some $n\ge 3$.  Clearly, $X$ is a smooth, quasi-projective 
variety, with $H_1(X,\Z)=\Z^{n}$.  Identify the character group of 
$\pi_1(X)$ with $(\C^{\times})^{n}$.  
Using Corollary \ref{cor:cv prod}, as well as 
Examples \ref{ex:cv wedge} and \ref{ex:cv link}, 
we see that 
\[
\WW^2(X)=\{1\}\cup 
\bigcup_{\zeta :  \Delta_K(\zeta)=0} \{\zeta\} \times (\C^{\times})^{n-1}. 
\]

Since $\Delta_K(1)\ne 0$,  we have that $\tau_1(\WW^2(X))=\{0\}$.  
On the other hand, Proposition \ref{prop:ttdf} gives that 
$\Omega^2_r(X)=\Grass_r(\{0\}\times \Q^{n-1})$ if $1< r<n$ and 
$\Omega^2_r(X)=\emptyset$ if $r=n$.  In the first case, 
$\Omega^2_r(X)$ is {\em not}\/ an open 
subset of $\Grass_r(\Q^{n})$, and in either case, 
$\Omega^2_r(X)$ is strictly included in 
$\sigma_r(\tau^{\Q}_1(\WW^2(X)))^{\compl}=\Grass_r(\Q^n)$.
\end{example} 

In degree $i=1$, we can say more.  As before, let $X$ be a smooth 
(connected, non-compact) quasi-projective variety. 
Given an orbifold fibration $f\colon X\to (\Sigma_{g,s},\mathbf{m})$, 
set $L=f^* (H^1(\Sigma_{g,s},\Q))$ and 
$A=\Tors(\pi_1^{\orb}(\Sigma_{g,s},\mathbf{m}))$, 
and denote by 
$\varphi\colon \pi_1(X)\surj  
\pi_1^{\orb}(\Sigma_{g,s},\mathbf{m})$ 
the homomorphism induced by $f$.  

\begin{prop}
\label{prop:df qp}
Let $X$ be a smooth, quasi-projective variety. 
For all $r\ge 1$, 
\begin{equation}
\label{eq:om1 qp}
\Omega^1_r(X) = \Grass_r(H^1(X,\Q)) \setminus 
\bigcup_{\alpha} \bigcup_{\rho} \sigma_r(L_\alpha, \rho),
\end{equation}
where the first union runs over the set of equivalence classes of 
orbifold fibrations 
$f_{\alpha}\colon X\to (\Sigma_{g_\alpha,s_\alpha},\mathbf{m}_{\alpha})$ 
with either 
\begin{enumerate}
\item \label{i1} $2g_{\alpha}+s_{\alpha}\ge 3$, or 
\item \label{i2} $g_{\alpha}=1$, $s_{\alpha}=0$ and $\abs{\mathbf{m}_{\alpha}}\ge 2$, or 
\item \label{i3} $g_{\alpha}=0$, $s_{\alpha}=2$ and $\abs{\mathbf{m}_{\alpha}}\ge 1$, 
\end{enumerate}
while the second union runs over all  
$\rho\in \hat{\varphi}_{\alpha}(\widehat{A}_{\alpha})$ in case \eqref{i1} or 
$\rho\in \hat{\varphi}_{\alpha}(\widehat{A}_{\alpha})\setminus \{1\}$  
in cases \eqref{i2} and \eqref{i3}.
\end{prop} 

\begin{proof}
In view of Propositions  \ref{prop:orbi} 
and \ref{prop:orbi noncompact}, as well as 
Theorem \ref{thm:arapura bis}, the conclusion follows 
from Theorem \ref{thm:df translated}.
\end{proof}

\begin{prop}
\label{prop:df quasiproj}
Let $X$ be a smooth, quasi-projective variety.   
Then, for all $r\ge 1$, 
\begin{equation}
\label{eq:omega qp}
\Omega^1_r(X) \subseteq \Grass_r(H^1(X,\Q))\setminus \bigcup_{\alpha} 
\sigma_r(f_{\alpha}^* (H^1(\Sigma_{g_{\alpha},s_{\alpha}},\Q)),
\end{equation}
where the union runs over the set of orbifold fibrations 
$f_{\alpha}\colon X\to \Sigma_{g_{\alpha},s_{\alpha}}$ 
with $2g_{\alpha}+s_{\alpha} \ge 3$.  Moreover, if $r=1$, 
or there are no orbifold fibrations with multiple fibers, 
then \eqref{eq:omega qp} holds as equality.
\end{prop}

\begin{proof}
Formula \eqref{eq:omega qp} is a consequence of 
Theorem  \ref{thm:sch bound} and Corollary \ref{cor:tcone qp}.  
The last two claims  follow from Theorem \ref{thm:df1} 
and Corollaries \ref{cor:tau schubert} and \ref{cor:qp nomult}, 
respectively.
\end{proof}

Using Corollary \ref{cor:df sch}, we obtain the following immediate 
consequence.

\begin{corollary}
\label{cor:dfqp empty} 
Suppose there is an orbifold fibration $f\colon X\to \Sigma_{g,s}$  
with $2g+s\ge 3$.  Then $\Omega^1_r(X) = \emptyset$, for all 
$r>b_1(X)-b_1(\Sigma_{g,s})$.
\end{corollary}

On the other hand, as the next result shows, 
if the characteristic variety $\WW^1(X)$ contains positive-dimensional 
translated components, its exponential tangent cone 
 may fail to determine the 
Dwyer--Fried sets $\Omega^1_r(X)$ with $r>1$. 

\begin{prop}
\label{prop:transtorus}
Let $X$ be a smooth, quasi-projective variety.   
Suppose that 
$\WW^1(X)$ has a $1$-dimensional component not 
passing through $1$, and 
$\tau_1(\WW^1(X))$ has codimension greater than $1$. 
Then $\Omega^1_2(X)$ is strictly contained in 
$\sigma_2(\tau^{\Q}_1(\WW^1(X))^{\compl}$.  
\end{prop}

A proof of this result is given in \cite[Theorem 10.11]{Su-aspm}, 
in the particular case when the quasi-projective variety $X$ is 
$1$-formal.  The argument given there extends without 
any essential modifications to this more general setting, 
upon replacing the resonance variety $\RR^1(X,\C)$ by 
the exponential tangent cone $\tau_1(\WW^1(X))$.

\begin{corollary}
\label{cor:tt}
Let $X$ be a smooth, quasi-projective variety.  
Suppose there is a pencil $X\to \C^{\times}$ 
with multiple fibers, but there is no pencil 
$X\to \Sigma_{g,s}$ with $b_1(X)-b_1(\Sigma_{g,s})\le 1$. 
Then $\Omega^1_2(X) \subsetneqq \Grass_2(H^1(X,\Q)) \setminus 
\sigma_2(\tau^{\Q}_1(\WW^1(X)))$.
\end{corollary}

\begin{proof}
Follows from Proposition \ref{prop:transtorus}, 
formulas \eqref{eq:v1piorb} and \eqref{eq:v1pi}, 
and Theorem \ref{thm:arapura bis}.
\end{proof}

\section{Hyperplane arrangements}
\label{sect:arr}

Among quasi-pro\-jective varieties,  
complements of complex hyperplane arrangements 
stand out as a particularly fascinating class of examples. 

\subsection{Characteristic varieties of arrangements}
\label{subsec:hyp arr}

Let $\A=\{H_1,\dots, H_n\}$ be a finite set of hyperplanes 
in some finite-dimensional complex vector space $\C^{d}$.  
Most of the time, we will assume the arrangement is central, 
i.e., all hyperplanes pass through the origin.  In this case, 
a defining polynomial for the arrangement is the product 
$Q=\prod_{i=1}^n \alpha_i$, where $\alpha_i\colon \C^{d} \to \C$ 
are linear forms with $\ker(\alpha_i)=H_i$.  

Let $X(\A)=\C^{d}\setminus \bigcup_{i=1}^n H_i$ be the 
complement of the arrangement.  This is a smooth, quasi-projective 
variety, with the homotopy type of a connected, finite CW-complex 
of dimension $d$.   The cohomology ring $H^*(X(\A),\Z)$ 
is the quotient of the exterior algebra on classes dual 
to the meridians around the hyperplanes, modulo a certain 
ideal generated in degrees $\ge 2$.  As shown by 
Orlik and Solomon in the 1980s, this ideal is determined 
by the intersection lattice, that is, the poset $L(\A)$ 
of all non-empty intersections of $\A$, ordered by 
reverse inclusion. 

Using the meridian basis of $H_1(X(\A),\Z)=\Z^n$, 
we may identify the character group of $\pi_1(X(\A))$ 
with the complex torus $(\C^{\times})^n$. 
By Theorem \ref{thm:arapura}, the characteristic varieties   
$\VV^i_1(X(\A))$ are unions of unitary translates of algebraic 
subtori in $(\C^{\times})^n$.   If $\A$ is central, then its 
complement is diffeomorphic to the product of $\C^{\times}$ with 
the complement in $\CP^{d-1}$ of the projectivization of $\A$.  
From Proposition \ref{prop:cv prod}, we see that 
$\VV^i_1(X(\A))$ is a subvariety of the complex torus 
$\{t \in (\C^{\times})^n \mid t_1\cdots t_n=1\}\cong (\C^{\times})^{n-1}$.

We will be only interested here in the degree-one characteristic 
variety, $\VV^1(\A):=\VV^1(X(\A))$.  By Theorem \ref{thm:arapura bis}, 
this variety is a union of torsion-translated subtori. Results  
from \cite{CS99} and \cite{DPS-duke} imply that the exponential 
tangent cone to $\VV^1(\A)$ coincides with the resonance 
variety $\RR^1(X(\A))$.  Thus, the components of $\VV^1(\A)$ 
passing through the origin are determined by the intersection 
lattice of $\A$.

If $\B\subset \A$ is a sub-arrangement, the inclusion 
$X(\A) \inj X(\B)$ induces an epimorphism 
$\pi_1(X(\A)) \surj \pi_1(X(\B))$. 
By Lemma \ref{lem:epi cv}, the resulting monomorphism 
between character groups restricts to an embedding 
$\VV^1(\B) \inj \VV^1(\A)$.  Components of $\VV^1(\A)$ 
which are not supported on any proper sub-arrangement 
are said to be essential. 

\subsection{Pencils, multinets, and translated tori}
\label{subsec:arr pencils}

Let us describe in finer detail the geometry of the 
characteristic variety $\VV^1(\A)$.  We start with a simple, 
yet basic example.

\begin{example}
\label{ex:pencil}
Let $\A$ be an arrangement of $n\ge 3$ lines through the 
origin of $\C^2$, defined by the polynomial $Q=x^n-y^n$. 
Then $X(\A)$ is diffeomorphic to $\C^{\times} \times \Sigma_{0,n}$.     
Using Example \ref{ex:cv wedge} and Proposition \ref{prop:cv prod}, 
we see that $\VV^1(\A)=(\C^{\times})^{n-1}$. Clearly, this single component  
arises as pull-back along the projection map $X(\A)\to \Sigma_{0,n}$. 
\end{example}

More generally, let $\A$ be a central arrangement in $\C^d$.  
Since $X(\A)$ admits a non-singular compactification with  
$b_1=0$, every orbifold fibration $f\colon 
X(\A) \to \Sigma_{g,s}$ must have base genus $g=0$. 

Taking a generic $2$-section, we obtain an arrangement $\A'$ 
of affine lines in $\C^2$, with $\pi_1(X(\A))=\pi_1(X(\A'))$; therefore, 
$\VV^1(\A)=\VV^1(\A')$.   Each intersection point of multiplicity 
$s\ge 3$ in $\A'$ gives rise to a ``local" component of $\VV^1(\A)$ 
of dimension $s-1$. Such components arise as pull-backs along 
pencils $X(\A)\to \Sigma_{0,s}$ similar to the ones from 
Example \ref{ex:pencil}.  In general, though, there are non-local 
components in  $\VV^1(\A)$.  

\begin{example}
\label{ex:braid arr}
Let $\A$ be the braid arrangement in $\C^3$, defined by the 
polynomial $Q=(x^2-y^2)(x^2-z^2)(y^2-z^2)$.  The variety 
$\VV^1(\A)\subset (\C^{\times})^6$ has $4$ local components 
of dimension $2$, corresponding to $4$ triple points in a generic 
section.  Additionally, the rational map 
$f\colon \CP^2 \dashrightarrow \CP^1$ 
given by $f(x,y,z)=(x^2-y^2,x^2-z^2)$ restricts to 
a holomorphic fibration $f\colon X(\A)\to \Sigma_{0,3}$, 
where  $\Sigma_{0,3}=\CP^1 \setminus \{ (1,0), (0,1), (1,1) \}$. 
This yields an essential, $2$-dimensional component in $\VV^1(\A)$. 
\end{example}

For an arbitrary central arrangement $\A$ in $\C^3$ with 
defining polynomial $Q=\prod_{i=1}^n \alpha_i$, all the 
irreducible components of $\VV^1(\A)$ passing through 
the origin can be described in terms of ``multinets" on 
the intersection lattice of $\A$.  As shown by Falk, Pereira, 
and Yuzvinsky in \cite{FY, PeY, Yu}, 
every such multinet determines a completely reducible curve 
$\{Q^{\mu}=0\}$ in $\CP^2$, where 
$Q^{\mu}=\prod_{i=1}^n \alpha_i^{\mu_i}$, 
for some vector $\mu\in \Z^n$.  In turn, this curve 
defines a pencil $X(\A)\to \Sigma_{0,s}$ for some 
$s\ge 3$, and this pencil produces an $(s-1)$-dimensional 
component of $\VV^1(\A)$, passing through $1$. 
Moreover, if this component is non-local, then $s=3$ or $4$. 

In general, though, the characteristic variety $\VV^1(\A)$ has  
irreducible components not passing through the origin. 

\begin{example}
\label{ex:deleted B3}
Let $\A$ be the deleted $\operatorname{B}_3$ 
arrangement, with defining polynomial 
$Q=xyz(x-y)(x-z)(y-z)(x-y-z)(x-y+z)$; a 
generic plane section of $\A$ is depicted 
in  Figure \ref{fig:deletedb3}.
The characteristic variety $\VV^1(\A)\subset (\C^{\times})^8$  
was computed in \cite{Su02}: it has six  
local components of dimension $2$, one local 
component of dimension $3$, five non-local 
components of dimension $2$ corresponding to braid 
sub-arrangements, and an essential component of 
the form $\rho T$, where 
$\rho=\exp(2\pi i \lambda)$, 
with $\lambda=(1/2,0,1/2,1/2,0,1/2,0,0)$, and  
$T=\exp(\ell \otimes\C)$, with $\ell$ the line in 
$\Q^8$ spanned by $\mu=(-1,1,0,0,1,-1,-2,2)$.

Since $\dim T=1$ and $\rho$ has order $2$, the component 
$\rho T$ must arise by pullback along a pencil 
$X(\A) \to (\Sigma_{0,2}, 2)$ with a single multiple fiber, of 
multiplicity $2$.   As noted in \cite{Di07}, this pencil is 
$Q^{\mu}\colon X(\A) \to \C^{\times}$, and the multiple 
fiber is the one over $1\in \C^{\times}$. 
\end{example}

Let us summarize the above discussion, as follows. As before, 
let $\A$ be an arrangement of $n$ hyperplanes.   

\begin{theorem}
\label{thm:cv arr}
Each irreducible component of $\VV^1(\A)$ is a 
torsion-translated subtorus of $(\C^{\times})^n$. 
Moreover, each positive-dimensional, non-local 
component is of the form $\rho T$, where $\rho$ is 
a torsion character, $T=f^*(H^1(\Sigma_{0,s},\C^{\times}))$, 
for some orbifold fibration $f\colon X(\A)\to (\Sigma_{0,s},\mathbf{m})$, 
and  either 
\begin{enumerate}
\item \label{tt1} $s=2$, and $f$ has at least one multiple fiber, or 

\item \label{tt23}  $s=3$ or $4$,  and $f$ corresponds to a 
multinet on $L(\A)$. 
\end{enumerate}
\end{theorem}

\begin{remark}
\label{rem:cv comb}
As shown in \cite{Su02}, there do exist arrangements $\A$ with 
isolated torsion points in $\VV^1(\A)$. Moreover, as noted in 
Example \ref{ex:deleted B3}, there also exist arrangements 
$\A$ with $1$-dimensional translated tori (of type \eqref{tt1}) 
in $\VV^1(\A)$.  On the other hand, we 
do not know whether higher-dimensional translated tori occur (these 
would necessarily be of type \eqref{tt23}, and the corresponding 
orbifold fibrations would have multiple fibers).

Of course, the local components in $\VV^1(\A)$ are determined 
by the intersection lattice of $\A$.  Moreover, the components 
of type \eqref{tt23} corresponding to orbifold fibrations with no 
multiple fibers are also determined by $L(\A)$.  It is not known 
whether the remaining components (i.e., those of dimension $0$, 
those of type \eqref{tt1}, and those of type \eqref{tt23} having 
non-trivial translation factors) are also combinatorially determined.  
\end{remark}

\subsection{Dwyer--Fried invariants of arrangements}
\label{subsec:df arr}

Let $\A=\{H_1,\dots , H_n\}$ be a central hyperplane 
arrangement.  A connected, regular $\Z^r$-cover of 
the complement $X(\A)$ is specified by assigning to each 
hyperplane $H_i$ the $i$-th column of an integral 
$r\times n$ matrix of rank $r$.  In turn, such a matrix 
defines a point in the Grassmannian of $r$-planes 
in $\Q^n$.  The Dwyer--Fried invariants of the arrangement, 
$\Omega^i_r(\A)=\Omega^i_r(X(\A))$, record the locus of 
points in $\Grass_r(\Q^n)$ for which the first $i$ Betti 
numbers  of the corresponding cover are finite. 
From Proposition \ref{prop:df qp}, we know that 
\begin{equation}
\label{eq:om1 arr}
\Omega^1_r(\A) = \Grass_r(\Q^n) \setminus 
\bigcup_{\alpha} \bigcup_{\rho} \sigma_r(L_\alpha, \rho),
\end{equation}
where the first union runs over the set of orbifold fibrations 
$f_{\alpha}\colon X(\A)\to (\Sigma_{0,s_\alpha},\mathbf{m}_{\alpha})$ 
with either $s_{\alpha}\ge 3$, or $s_{\alpha}=2$ and 
$\abs{\mathbf{m}_{\alpha}}>0$, while the second union 
runs over all characters (respectively, non-trivial characters) 
$\rho$ in $\widehat{(f_{\alpha})_{\sharp}}(\widehat{A_{\alpha}})$, 
and where $L_{\alpha}=
f_{\alpha}^* (H^1(\Sigma_{0,s_{\alpha}},\Q))$ and 
$A_{\alpha}=\Tors(\pi_1^{\orb}(\Sigma_{0,s_{\alpha}}, 
\mathbf{m}_{\alpha}))$. 

Note that the linear subspaces $L_{\alpha}\subset \Q^n$ comprise 
the characteristic subspace arrangement $\CC_1(X(\A))$. 
Each such subspace has dimension $c_{\alpha}= s_{\alpha}-1$, 
with $1< c_{\alpha} <n$. Moreover, if $\exp(L_{\alpha}\otimes \C)$ 
is not a local component, then $c_{\alpha}=2$ or $3$.  
By Corollary \ref{cor:dfqp empty},  the set 
$\Omega^1_r(\A)$ is empty for all $r > n - c$,   
where $c = \max\, \{c_{\alpha} \}$.  In particular, 
if a generic $2$-section of $\A$ has an intersection 
point of multiplicity $s\ge 3$, then $\Omega^1_{n-s+2}(\A)=\emptyset$.

In view of Remark \ref{rem:cv comb}, it remains an open question 
whether the Dwyer--Fried sets of an arrangement are combinatorially 
determined.  Nevertheless, Proposition \ref{prop:df quasiproj} 
yields the following combinatorial upper bound for these sets:
\begin{equation}
\label{eq:omega arr bound}
\Omega^1_r(\A) \subseteq \Grass_r(\Q^n)\setminus \bigcup_{\alpha} 
\sigma_r(L_{\alpha}),
\end{equation}
where the union runs over the set of orbifold fibrations 
$f_{\alpha}\colon X(\A) \to \Sigma_{0,s_{\alpha}}$, $s_{\alpha}\ge 3$,  
corresponding to multinets on $L(\A)$. 
Of course, if either $r=1$, or there are no orbifold fibrations 
with multiple fibers, then the above inclusion holds as equality.

\begin{example}
\label{ex:df nr}
Let $\A$ be an arrangement of $n$ planes 
in $\C^{3}$, and suppose a generic $2$-section 
has one or two lines which contain all the intersection points 
of multiplicity $3$ and higher.  Then, as shown by Nazir and 
Raza in \cite{NR09}, the variety $\VV^1(\A)$ has no translated components.   
Thus, \eqref{eq:omega arr bound} holds as equality in this case.
\end{example}

Here is an example where the above condition is not satisfied, 
yet equality is still attained in \eqref{eq:omega arr bound}.

\begin{example}
\label{ex:df braid}
If $\A$ is the braid arrangement from Example \ref{ex:braid arr}, 
then all components of $\VV^1(\A)$ pass through the origin, 
and $\CC_1(X(\A))$ consists of $5$ planes in $\Q^6$. 
Thus, $\Omega^1_r(\A) =\Grass_r(\Q^6)\setminus  
\bigcup_{i=1}^{5}\sigma_r(L_i)$ for $r\le 4$, and 
$\Omega^1_r(\A) =\emptyset$ for $r> 4$. 
\end{example} 

In general, though, translated tori in the characteristic variety 
$\VV^1(\A)$ will affect the $\Omega$-invariants of an 
arrangement $\A$.

\begin{figure}[t]
\setlength{\unitlength}{0.4cm}
\begin{picture}(16,10)(2,1)
\put(0,8){\line(1,0){15.5}}
\put(14,0.8){\line(0,1){8.5}}
\put(2,8){\line(3,-1){13.5}}
\put(2,8){\line(-3,1){1.75}}
\put(2,8){\line(4,-1){13.5}}
\put(2,8){\line(-4,1){1.5}}
\put(8,8){\line(1,-1){7}} 
\put(8,8){\line(-1,1){1.3}}
\put(8,8){\line(3,-2){7.5}}
\put(8,8){\line(-3,2){1.8}}
\put(10,8){\line(2,-3){4.65}}
\put(10,8){\line(-2,3){0.9}}
\put(14,8){\line(-1,-1){5}}
\put(14,8){\line(1,1){0.9}}
\put(15.8,9.5){\makebox(0,0){$\mat{0\\1}$}}   
\put(16.5,8){\makebox(0,0){$\mat{1\\0}$}}   
\put(16.7,5){\makebox(0,0){$\mat{\!-2\\0}$}}   
\put(16.6,3.6){\makebox(0,0){$\mat{1\\0}$}}   
\put(15.6,2.3){\makebox(0,0){$\mat{0\\1}$}}   
\put(7.2,10.1){\makebox(0,0){$\mat{\!-1\\1}$}}   
\put(9.6,10.1){\makebox(0,0){$\mat{2\\0}$}}   
\put(13.5,10.1){\makebox(0,0){$\mat{\!-1\\1}$}}   
\end{picture}
\caption{Generic section of the deleted ${\rm B}_3$ 
arrangement and a $\Z^2$-cover}
\label{fig:deletedb3}
\end{figure}
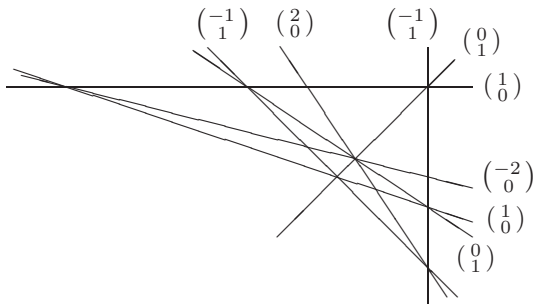

\begin{example}
\label{ex:df deleted B3}
Let $\A$ be the deleted $\operatorname{B}_3$ 
arrangement from Example \ref{ex:deleted B3}. 
We know that the characteristic subspace arrangement 
$\CC=\CC_1(X(\A))$ consists of eleven $2$-planes and 
a $3$-plane in $\Q^8$.  Moreover, we know that $\VV^1(\A)$ 
contains a $1$-dimensional translated torus.  Therefore, 
Corollary \ref{cor:tt} tells us that $\Omega^1_2(\A)$ 
is strictly contained in the set 
$U=\Grass_2(\Q^8) \setminus \bigcup_{L\in \CC} \sigma_2(L)$.  

In the notation from Example \ref{ex:deleted B3}, 
an explicit element in $U\setminus \Omega^1_2(\A)$ is the plane 
$P$ spanned by the vectors $\mu$ and $2\lambda$. 
The corresponding $\Z^2$-cover of $X(\A)$ is obtained  
by sending the meridians to the vectors indicated 
in Figure \ref{fig:deletedb3}.
\end{example}

\section{Finiteness properties of groups}
\label{sect:ff fp}

In this final section, we discuss some alternate ways to measure 
the finiteness properties of a group.  It turns out that 
the $\Omega$-invariants give subtle information about 
those properties, too.

\subsection{Properties $\FF_k$ and $\FP_k$}
\label{subsec:ff fp}

Perhaps the most basic finiteness property for discrete 
groups is the geometric finiteness condition introduced 
by C.T.C.~Wall in \cite{Wa}.  
A group $G$ is said to satisfy property $\FF_{k}$, for 
some $k \ge 1$, if $G$ admits a classifying space 
$K(G,1)$ with finite $k$-skeleton. 
Evidently, the $\FF_1$ property is equivalent to $G$ being 
finitely generated, while the $\FF_2$ property is equivalent 
to $G$ being finitely presentable.  

These notions have an algebraic analogue, introduced 
by J.P.~Serre \cite{Se} and R.~Bieri \cite{Bi}.
A group $G$ is said to satisfy property $\FP_{k}$, for some 
$k \ge 1$, if the trivial $\Z{G}$-module $\Z$ admits a partial 
resolution, $P_k \to \cdots \to P_1 \to P_0 \to \Z \to 0$, with all 
$P_i$ finitely generated, projective $\Z{G}$-modules. 

Clearly, if $G$ is of type $\FF_{k}$, then it is also of type 
$\FP_{k}$.  For $k=1$, the two conditions are equivalent.  
Furthermore, if $G$ is finitely presented and of type $\FP_k$, 
then $G$ is of type $\FF_k$.  For each $k\ge 2$,
Bestvina and Brady produced in \cite{BB} 
examples of groups of type $\FP_k$, but not of type $\FF_k$. 
Finally, if $G$ is of type $\FP_k$, then $H_i(G, \Z)$ is finitely 
generated, for all $i\le k$, and thus $b_i(G)<\infty$, for all $i\le k$. 
In general, though, none of these implications can be reversed 
(see Example \ref{ex:milnor} for one instance of this phenomenon).

\subsection{Dwyer--Fried invariants and finiteness properties}
\label{subsec:df fp}

We now present a concrete way in which the $\Omega$-sets  
can be used to inform on the aforementioned finiteness properties 
of groups.   The idea goes back to Stallings' seminal paper \cite{Sta}; 
we refer to Example \ref{ex:Stallings} and the proof of 
Theorem \ref{thm:kollar answer} for some concrete applications. 

\begin{theorem}
\label{thm:fpk}
Let $G$ be a finitely generated group, and let 
$\nu\colon G\surj \Z^r$ be an epimorphism, with 
kernel $\vG$.  Suppose that
$\Omega^{k}_{r}(G)=\emptyset$ and 
$\vG$ is of type $\FF_{k-1}$. 
Then $b_k(\vG)=\infty$.  Consequently, $H_k(\vG,\Z)$ is 
not finitely generated, and thus $\vG$ is not of type 
$\FP_k$.  
\end{theorem}

\begin{proof}
Set $X=K(G,1)$; then $X^{\nu}=K(\vG,1)$.  
Since $\vG$ is of type $\FF_{k-1}$, we have $b_i(X^{\nu})<\infty$ 
for $i\le k-1$.  On the other hand, since $\Omega^{k}_{r}(X)=\emptyset$, 
we must have $b_k(X^{\nu})=\infty$.  The conclusions follow.
\end{proof}

\begin{corollary}
\label{cor:fpk}
Let $G$ be a finitely generated group, and suppose that 
$\Omega^{3}_{1}(G)=\emptyset$.  Let $\nu\colon G\surj \Z$ 
be an epimorphism.  If the group $\vG=\ker(\nu)$ is finitely 
presented, then $H_3(\vG,\Z)$ is not finitely generated. 
\end{corollary}

\begin{example}
\label{ex:Stallings}
Let $Y_2=S^1\vee S^1$ and $X=Y_2\times Y_2 \times Y_2$. 
Clearly, $X$ is a classifying space for the group 
$G=F_2\times F_2\times F_2$.  Let $\nu\colon G\surj\Z$ 
be the homomorphism taking each standard generator to $1$.  
In \cite{Sta}, Stallings found an explicit finite presentation 
for the group $\vG=\ker(\nu)$:
\begin{equation}
\label{eq:st group}
\vG=\langle a,b,c,x,y \mid [x,a],\, [y,a],\, [x,b],\, [y,b],\, 
[a^{-1}x,c],\,  [a^{-1}y,c],\,  [b^{-1}a,c]\rangle.
\end{equation}

He then proceeded to show, via a Mayer--Vietoris 
argument, that $H_3(\vG,\Z)$ is not finitely generated.  
This last assertion can be readily explained from our point 
of view.  Indeed, by formula \eqref{eq:df prod wedge}, we have 
that $\Omega^3_1(X)=\emptyset$.  The desired conclusion 
then follows from Corollary \ref{cor:fpk}.
\end{example}

\begin{remark}
\label{rem:stallings x3}
It turns out that the Stallings group is a (non-central) 
arrangement group.  This fact was observed by the author 
during a conversation with Daniel Matei in May 2004, and announced 
in a talk given at MSRI in August 2004 \cite{Su04}. Since the explicit 
computation has not appeared in print, let us record it here. 

Consider the arrangement of $5$ lines in $\C^2$ defined 
by the vanishing of the polynomial $Q=zw(w+1)(z-1)(2z+w)$. 
Using the braid monodromy generators recorded in 
\cite[Example 10.2]{Su01}, we obtain the following 
presentation for the fundamental group of the complement 
of this arrangement:
\begin{equation}
\label{eq:x3 group}
G=\left\langle x_1,\dots,x_5\
\Bigg|  
\begin{array}{ll}%
[x_1,x_2],\: [x_1,x_4],\: [x_2,x_4],\: [x_3,x_4],\: [x_2, x_5 ], 
  \\[3pt]
x_1 x_3 x_5=x_5 x_1 x_3=x_3 x_5 x_1
\end{array}
\right\rangle,
\end{equation}
with generators corresponding to the lines, and relators 
corresponding to the intersection points.  An isomorphism 
$\Gamma \isom G$ is given by 
$a \mapsto  x_1 x_3 x_4 x_5, \: b \mapsto  x_1 x_3 x_5, \:
c \mapsto  x_2,\: x \mapsto  x_1 x_3,\: y \mapsto  x_3$.
\end{remark}

\subsection{Koll\'{a}r's question}
\label{subsec:kollar}

Two groups, $G_1$ and $G_2$, are said to be 
{\em commensurable up to finite kernels}\/
if they can be connected by a zig-zag of groups and 
homomorphisms, 
with all arrows of finite kernel and cofinite image. 
In \cite{Ko}, J.~Koll\'{a}r asked the following question: 
{\em Given a smooth, projective variety $M$, is the 
fundamental group $\vG=\pi_1(M)$ commensurable, 
up to finite kernels, with another group, $\pi$, admitting 
a $K(\pi, 1)$ which is a quasi-projective variety}? 
In \cite{DPS-crelle}, we answered this question, as follows.

\begin{theorem}[\cite{DPS-crelle}]
\label{thm:kollar answer}
For each $k\ge 3$, there is a smooth, irreducible, complex 
projective variety $M$ of complex dimension $k-1$, such 
that the group $\vG=\pi_{1}(M)$ is of type $\FF_{k-1}$, 
but not of type $\FP_k$.   
\end{theorem}

Using some classical results of R.~Bieri, it is easy to see that 
such a group $\vG$ is not commensurable, up to finite kernels, 
to any group of type $\FP_k$, and thus, to any group 
$\pi$ admitting a $K(\pi, 1)$ which is a quasi-projective variety. 
As noted in \cite{Ey}, the manifolds constructed in 
Theorem \ref{thm:kollar answer} cannot have a simply-connected 
holomorphic embedding into a compactifiable complex 
manifold with contractible universal covering space.

\subsection{Branched covers and generic fibers}
\label{subsec:br cover}

We now sketch the proof of Theorem \ref{thm:kollar answer}, 
in a streamlined version afforded by the technology developed here.  
The starting point is a classical branched covering construction.

Let $C$ be a bielliptic curve of genus $g\ge 2$.  Then $C$ 
supports an involution $\sigma$, such that the quotient 
$C/\sigma$ is an elliptic curve $E$.  (The case $g=2$ 
is depicted on the left side of Figure \ref{fig:g2g3}.) 
The projection map, $f\colon C\to E$, is a surjective 
holomorphic map, and can be viewed as a $2$-fold 
branched cover, with branch set $B\subset E$ consisting 
of $2g-2$ points. The corresponding unramified $2$-fold 
cover is classified by a homomorphism 
$\varphi \colon \pi_1(E\setminus B) \surj \Z_2$ 
taking each standard loop around a branch point to $1$. 

Now fix an integer $k\ge 3$, and set $X=C^{\times k}$.  
The group law of the elliptic curve extends by associativity 
to a map $s_k\colon E^{\times k} \to E$.  Composing this map 
with the product map  $f^{\times k}\colon C^{\times k} 
\to E^{\times k}$, we obtain a surjective holomorphic map, 
$h=s_k \circ f\colon X \to E$. 

\begin{lemma}[\cite{DPS-crelle}]
\label{lem:dps crelle}
Let $M$ be the generic fiber of $h$.  Then $M$ is a smooth, 
complex projective variety of dimension $k-1$. Moreover,
\begin{enumerate}
\item \label{m1}  $M$ is connected.  
\item \label{m2} $\pi_1(M)=\ker \big( h_{\sharp}\colon 
\pi_1(X)\surj \pi_1(E)\big)$. 
\item \label{m3} $\pi_2(M)=\cdots =\pi_{k-2}(M)=0$.  
\end{enumerate}
\end{lemma}

Assuming this lemma, it is now an easy matter to 
finish the proof of the theorem.  

\begin{proof}[Proof of Theorem \ref{thm:kollar answer}]
Set $G=\pi_1(X)$ and $\vG=\pi_1(M)$.  Identify $\pi_1(E)=\Z^2$, 
and write $\nu=h_{\sharp}$.  From Lemma \ref{lem:dps crelle}, 
parts \eqref{m1} and \eqref{m2}, we obtain a short exact sequence,
$1 \to \vG \to G \xrightarrow{\nu}  \Z^2 \to  1$. 
Since $X$ is a $k$-fold direct product of surfaces 
of genus $g\ge 2$, the space $X$ is a $K(G,1)$.  
Moreover, formula \eqref{eq:df prod surf} 
shows that $\Omega^k_2(G)=\emptyset$. 

In view of Lemma \ref{lem:dps crelle}, part \eqref{m3}, a classifying  
space $K(\vG,1)$ can be obtained from $M$ by attaching cells 
of dimension $k$ and higher.  Consequently, the group $\vG$ is of 
type $\FF_{k-1}$.  Finally, Theorem \ref{thm:fpk} shows that $\vG$ 
is not of type $\FP_{k}$. 
\end{proof}

\newcommand{\arxiv}[1]
{\texttt{\href{http://arxiv.org/abs/#1}{arXiv:#1}}}
\newcommand{\doi}[1]
{\texttt{\href{http://dx.doi.org/#1}{doi:#1}}}
\renewcommand{\MR}[1]
{\href{http://www.ams.org/mathscinet-getitem?mr=#1}{MR#1}}
\newcommand{\MRh}[2]
{\href{http://www.ams.org/mathscinet-getitem?mr=#1}{MR#1 (#2)}}


\begin{thebibliography}{00}

\bibitem{Ak} S.~Akbulut, 
{\em The Catanese--Ciliberto--Mendes Lopes surface}, 
J. G\"{o}kova Geom. Topol. \textbf{5} (2011), 86--102.
\MR{2872551}

\bibitem{AP}  A.~Akhmedov, B.~D.~Park, 
{\em Exotic smooth structures on $S^2\times S^2$}, 
\arxiv{1005.3346v5}

\bibitem{Ar} D.~Arapura, 
{\em Geometry of cohomology support loci for local systems.  
\textup{I}.}, J. Algebraic Geom. \textbf{6} (1997), no.~3, 563--597.  
\MR{1487227} 

\bibitem{ACM} E.~Artal Bartolo, J.~I.~Cogolludo-Agust\'{\i}n, D.~Matei, 
{\em Characteristic varieties of quasi-projective manifolds 
and orbifolds}, \arxiv{1005.4761v4}.

\bibitem{BCGP} I.~Bauer, F.~Catanese, F.~Grunewald, R.~Pignatelli, 
{\em Quotients of products of curves, new surfaces with $p_g=0$ 
and their fundamental groups}, Amer. J. Math. \textbf{134} (2012), 
no. 4, 993--1049.

\bibitem{Be92} A.~Beauville, 
{\em Annulation du $H\sp 1$ pour les fibr\' es en droites plats}, in: 
Complex algebraic varieties (Bayreuth, 1990), 1--15,
Lecture Notes in Math., vol.~1507, Springer, Berlin, 1992. 
\MR{1178716} 

\bibitem{BB} M.~Bestvina, N.~Brady, 
{\em Morse theory and finiteness properties of groups}, 
Invent. Math. \textbf{129} (1997), no.~3, 445--470. 
\MR{1465330}  

\bibitem{Bi}  R.~Bieri, 
{\em Homological dimension of discrete groups}, 
Second edition, Queen Mary Coll. Math. Notes, 
Queen Mary College, Dept. Pure Math., London, 1981. 
\MR{0715779} 

\bibitem{Br}  K.~S.~Brown,
{\em Cohomology of groups}, Grad. Texts in Math., 
vol.~87, Springer-Verlag, New York-Berlin, 1982.
\MR{0672956} 

\bibitem{Bu09}  N.~Budur, 
{\em Unitary local systems, multiplier ideals, and polynomial 
periodicity of Hodge numbers}, Adv. Math. \textbf{221} (2009), 
no.~1, 217--250.
\MR{2509325}  

\bibitem{Cm01} F.~Campana, 
{\em Ensembles de {G}reen-{L}azarsfeld et 
quotients resolubles des groupes de K\"{a}hler},
J. Alg. Geometry \textbf{10} (2001), no.~4, 599--622.
\MR{1838973}  

\bibitem{Cm11}  F.~Campana, 
{\em Quotients r\'{e}solubles ou nilpotents des groupes 
de K\"{a}hler orbifoldes}, Manuscripta Math. \textbf{135} 
(2011), no. 1-2, 117--150. 
\MR{2783390}  

\bibitem{Cat91} F.~Catanese, 
{\em Moduli and classification of irregular Kaehler manifolds 
(and algebraic varieties) with Albanese general type fibrations}, 
Invent. Math. \textbf{104} (1991), no.~2, 263--289. 
\MR{1098610}  

\bibitem{CCM} F.~Catanese, C.~Ciliberto, M.~Mendes Lopes, 
{\em On the classification of irregular surfaces of general type 
with nonbirational bicanonical map}, 
Trans. Amer. Math. Soc. \textbf{350} (1998), no.~1, 275--308. 
\MR{1422597}  

\bibitem{CS99} D.~Cohen, A.~Suciu, 
{\em Characteristic varieties of arrangements},
Math. Proc. Cambridge Phil. Soc. \textbf{127} (1999), 
no.~1, 33--53. 
\MR{1692519}  

\bibitem{Dz08} T.~Delzant,
{\em Trees, valuations and the {G}reen--{L}azarsfeld set},
Geom. Funct. Anal. \textbf{18} (2008), no.~4, 1236--1250. 
\MR{2465689} 

\bibitem{Di07} A.~Dimca, 
{\em Characteristic varieties and constructible sheaves}, 
Rend. Lincei Mat. Appl. \textbf{18} (2007), no.~4, 365--389. 
\MR{2349994}  

\bibitem{DPS-imrn} A.~Dimca, S.~Papadima, A.~Suciu,
{\em Alexander polynomials: {E}ssential variables and 
multiplicities}, Int. Math. Res. Notices \textbf{2008}, 
no.~3, Art. ID rnm119, 36 pp. 
\MR{2416998}  

\bibitem{DPS-crelle} A.~Dimca, \c{S}.~Papadima, A.~Suciu,
{\em Non-finiteness properties of fundamental groups of 
smooth projective varieties}, J. Reine Angew.~Math. 
\textbf{629} (2009), 89--105.  
\MR{2527414} 

\bibitem{DPS-duke} A.~Dimca, \c{S}.~Papadima, A.~Suciu,
{\em Topology and geometry of cohomology jump loci}, 
Duke Math. J. \textbf{148} (2009), no.~3, 405--457.
\MR{2527322} 

\bibitem{DF} W.~G.~Dwyer, D.~Fried,
{\em Homology of free abelian covers. \textup{I}}, Bull. 
London Math. Soc. \textbf{19} (1987), no.~4, 350--352.
\MR{0887774} 

\bibitem{EL97} L.~Ein, R.~Lazarsfeld, 
{\em Singularities of theta divisors and the birational geometry 
of irregular varieties}, J. Amer. Math. Soc. \textbf{10} (1997), 
no.~1, 243--258. 
\MR{1396893}  

\bibitem{Ey} P.~Eyssidieux,
{\em Lectures on the Shafarevich conjecture on uniformization}, 
to appear in {\em Vari\'et\'es Complexes, Feuilletages, Uniformisation},  
Panoramas et Synth\`eses,  Soci\'et\'e Math\'ematique de France, 
Paris. 

\bibitem{FY} M.~Falk, S.~Yuzvinsky,
{\em Multinets, resonance varieties, and pencils of plane curves},
Compositio Math. \textbf{143} (2007), no.~4, 1069--1088.
\MR{2339840} 

\bibitem{GL87} M.~Green, R.~Lazarsfeld, 
{\em Deformation theory, generic vanishing theorems
and some conjectures of Enriques, Catanese and Beauville}, 
Invent. Math. \textbf{90} (1987), no.~2, 389--407.
\MR{0910207}  

\bibitem{GL91} M.~Green, R.~Lazarsfeld, 
{\em Higher obstructions to deforming cohomology groups 
of line bundles}, J. Amer. Math. Soc. \textbf{4} (1991), 
no.~1, 87--103.
\MR{1076513}  

\bibitem{HP}  C.~Hacon, R.~Pardini, 
{\em Surfaces with $p_g=q=3$}, Trans. Amer. Math. Soc. 
\textbf{354} (2002), no.~7, 2631--2638. 
\MR{1895196}  

\bibitem{Har} J.~Harris, 
{\em Algebraic geometry}, Grad. Texts in Math, 
vol.~133, Springer-Verlag, New York, 1992. 
\MR{1182558} 

\bibitem{Hi} E.~Hironaka,
{\em Torsion points on an algebraic subset of an affine torus}, Int.
Math. Res. Notices \textbf{1996}, no.~19, 953--982.
\MR{1422370} 

\bibitem{Ko} J.~Koll\'{a}r,
{\em Shafarevich maps and automorphic forms}, 
Princeton Univ. Press, Princeton, NJ, 1995. 
\MR{1341589} 

\bibitem{Li02}  A.~Libgober,
{\em First order deformations for rank one local systems 
with a non-vanishing cohomology}, Topology Appl. 
\textbf{118} (2002), no.~1-2, 159--168. 
\MR{1877722} 

\bibitem{MP} A.~M\u{a}cinic,  \c{S}.~Papadima, 
{\em Characteristic varieties of nilpotent groups and 
applications}, in: Proceedings of the Sixth Congress 
of Romanian Mathematicians (Bucharest, 2007), 
pp.~57--64, vol.~1, Romanian Academy, Bucharest, 2009.
\MR{2641547}  

\bibitem{Mil} J.~Milnor, 
{\em Infinite cyclic coverings}, in: Conference on the 
Topology of Manifolds (Michigan State Univ., 1967), 
pp.~115--133, Prindle, Weber \& Schmidt, Boston, MA, 1968. 
\MR{0242163} 

\bibitem{NR09} S.~Nazir, Z.~Raza, 
{\em Admissible local systems for a class of line arrangements}, 
Proc. Amer. Math. Soc. \textbf{137} (2009), no.~4, 1307--1313. 
\MR{2465653} 

\bibitem{PS-plms} S.~Papadima, A.~Suciu,
{\em Bieri--{N}eumann--{S}trebel--{R}enz invariants and 
homology jumping loci}, Proc.~London Math.~Soc. 
\textbf{100} (2010), no.~3, 795--834.
\MR{2640291}  

\bibitem{PeY} J.~Pereira, S.~Yuzvinsky, 
{\em Completely reducible hypersurfaces in a pencil}, 
Adv.  Math.  \textbf{219} (2008), no.~2, 672--688.
\MR{2435653}  

\bibitem{Ro} D.~Rolfsen,
{\em Knots and links}, Math. Lecture Series, vol.~7,
Publish or Perish, Berkeley, CA, 1976.
\MR{0515288}  

\bibitem{Se} J.-P.~Serre,  
{\em Cohomologie des groupes discrets}, in: 
{\em Prospects in mathematics}, 77--169,  Ann. of Math. 
Studies, no.~70, Princeton Univ. Press, Princeton, NJ, 1971. 
\MR{0385006} 

\bibitem{Si93}  C. Simpson, 
{\em Subspaces of moduli spaces of rank one local systems}, 
Ann. Sci. \'{E}cole Norm. Sup. \textbf{26} (1993), no.~3, 361--401.
\MR{1222278}  

\bibitem{Sta}  J.~Stallings,  
{\em A finitely presented group whose $3$-dimensional 
integral homology is not finitely generated}, 
Amer. J. Math. \textbf{85} (1963), 541--543.
\MR{0158917}  

\bibitem{Su01} A.~Suciu,
{\em Fundamental groups of line arrangements: 
{E}numerative aspects}, in: {\em Advances in algebraic geometry 
motivated by physics (Lowell, MA, 2000)}, 43--79, Contemp. 
Math., vol.~276, Amer. Math. Soc., Providence, RI, 2001.
\MR{1837109}  

\bibitem{Su02} A.~Suciu,
{\em Translated tori in the characteristic varieties of
complex hyperplane arrangements}, Topology Appl. 
\textbf{118} (2002), no.~1-2, 209--223.  
\MR{1877726} 

\bibitem{Su04} A.~Suciu,
{\em Fundamental groups and characteristic varieties}, talk at 
the Introductory Workshop in Hyperplane Arrangements and 
Applications, MSRI, Berkeley, CA, 2004, available at 
\url{http://www.msri.org/web/msri/online-videos/-/video/showVideo/1997}.

\bibitem{Su11} A.~Suciu,
{\em Fundamental groups, Alexander invariants, and
cohomology jumping loci}, in: {\em Topology of 
algebraic varieties and singularities}, 179--223, Contemp. 
Math., vol. 538, Amer. Math. Soc., Providence, RI, 2011. 
\MR{2777821}  

\bibitem{Su-aspm} A.~Suciu,
{\em Resonance varieties and Dwyer--Fried invariants},  
in: {\em Arrangements of Hyperplanes (Sapporo 2009)}, 
359--398, Advanced Studies Pure Math., vol.~62, Kinokuniya, 
Tokyo, 2012.  
\MR{2933803}

\bibitem{Su-pisa} A.~Suciu,
{\em Geometric and homological finiteness in free 
abelian covers}, \arxiv{1112.0948v2}, to appear in 
{\em Configuration Spaces: Geometry, Combinatorics 
and Topology (Centro De Giorgi, 2010)}, 461--501, 
Edizioni della Normale, Pisa, 2012.

\bibitem{SYZ2} A.~Suciu, Y.~Yang, G.~Zhao, 
{\em  Intersections of translated algebraic subtori}, 
J. Pure Appl. Algebra \textbf{217} (2013), no.~3, 481--494.
\MR{2974227}

\bibitem{SYZ1} A.~Suciu, Y.~Yang, G.~Zhao, 
{\em  Homological finiteness of abelian covers}, 
\arxiv{1204.4873v1}.

\bibitem{Wa}  C.~T.~C.~Wall,  
{\em Finiteness conditions for ${\rm CW}$-complexes}, 
Ann. of Math. \textbf{81} (1965), 56--69.
\MR{0171284}  

\bibitem{Yu} S.~Yuzvinsky, 
{\em A new bound on the number of special fibers in a 
pencil of curves},  Proc. Amer. Math. Soc. \textbf{137} (2009), 
no.~5, 1641--1648. 
\MR{2470822}  

\end{thebibliography}
\end{document}